\documentclass[11pt]{article}
\usepackage[applemac]{inputenc}
\usepackage{latexsym,amsmath,amscd,amssymb,framed,graphics,color,mathrsfs,stmaryrd}
\usepackage{enumerate}
\usepackage{geometry}
\usepackage[colorlinks]{hyperref}
\usepackage[T1]{fontenc}
\usepackage{graphicx}
\definecolor{sand}{rgb}{0.76, 0.7, 0.5}
\definecolor{taupegray}{rgb}{0.55, 0.52, 0.54}
\usepackage[colorlinks]{hyperref}
\usepackage{subcaption}
\hypersetup{%
pdfstartview={FitH},
urlcolor=sand,
citecolor=blue,
linkcolor=taupegray,
}

\usepackage{array}
\usepackage{float}
\usepackage{multicol}
\usepackage{aeguill}
\usepackage{fancyhdr}
\usepackage[arrow, matrix, curve]{xy}
\usepackage{url}
\usepackage{colordvi}
\usepackage{framed}
\usepackage{color}
\usepackage{algorithmic}
\usepackage{algorithm}

\usepackage{manfnt}
\usepackage[bbgreekl]{mathbbol}
\usepackage{tikz}

\usepackage{enumitem}

\newcommand{\mathsym}[1]{{}}

\newtheorem{theorem}{Theorem}[section]
\newtheorem{definition}[theorem]{Definition}

\newtheorem{proposition}[theorem]{Proposition}

\newtheorem{example}[theorem]{Example}
\begin{document}
\title{Unified discrete multisymplectic Lagrangian formulation for hyperelastic solids and barotropic fluids}

\author{Fran\c{c}ois Demoures$^{1}$,  Fran\c{c}ois Gay-Balmaz$^{2}$}
\addtocounter{footnote}{1} \footnotetext{EPFL, Doc \& Postdoc Alumni. Av Druey 1, Lausanne 1018, Switzerland \\\texttt{francois.demoures@alumni.epfl.ch}}
\addtocounter{footnote}{2}
\footnotetext{Laboratoire de M\'et\'eorologie Dynamique, \'Ecole 
Normale Sup\'erieure/CNRS, Paris, France.
\texttt{francois.gay-balmaz@lmd.ens.fr}}

\maketitle

\date{}

\begin{abstract}
We present a geometric variational discretization of nonlinear elasticity in 2D and 3D in the Lagrangian description. A main step in our construction is the definition of discrete deformation gradients and discrete Cauchy-Green deformation tensors, which allows for the development of a general discrete geometric setting for frame indifferent isotropic hyperelastic models. The resulting discrete framework is in perfect adequacy with the multisymplectic discretization of fluids proposed earlier by the authors. Thanks to the unified discrete setting, a geometric variational discretization can be developed for the coupled dynamics of a fluid impacting and flowing on the surface of an hyperelastic body. The variational treatment allows for a natural inclusion of incompressibility and impenetrability constraints via appropriate penalty terms. 
We test the resulting integrators in 2D and 3D with the case of a barotropic fluid flowing on incompressible rubber-like nonlinear models.
\end{abstract}
 
 
\section{Introduction}

We present a discrete variational formulation of hyperelasticity in the Lagrangian description. This formulation is associated to the multisymplectic geometric formulation of continuum mechanics which underlies both fluid dynamics and nonlinear elasticity. From the discrete variational setting, a structure preserving numerical integrator is derived for hyperelastic solids, which is spacetime multisymplectic, symplectic in time, preserves exactly the momenta associated to symmetries, and nearly preserves total energy. Thanks to its variational character, the numerical scheme can naturally incorporate incompressibility constraints. Being in perfect adequacy with the multisymplectic discretization of fluids proposed earlier in \cite{DeGB2021}, the discrete framework can be extended to treat fluid-structure interaction dynamics by the inclusion of appropriate impenetrability constraints in the variational formulation.

The developments in this paper are founded on the multisymplectic geometric formulation of continuum mechanics, which is the spacetime version of the symplectic formulation of classical mechanics. This geometric setting allows to intrinsically formulate the multisymplectic form formula and the covariant Noether theorem, which are the field theoretic analogue to the symplecticity of the flow and to the conservation of momentum maps in classical mechanics.
The more usual infinite dimensional dynamical system formulation of continuum mechanics is recovered by space and time splitting, which also has a discrete analogue.

Among the main concepts in this paper are the discrete deformation gradients, the discrete Cauchy-Green deformation tensors, and the discrete Jacobians. These objects, directly defined from the spacetime discrete configuration map of the body, allow a systematic treatment of various models of hyperelasticity via the appropriate definition of the discrete stored energy functions. The discrete variational geometric setting also systematically guides the choice of the appropriate degree of freedom of the discrete objects, for instance, at each spatial node there is one discrete Cauchy-Green tensor associated to each of the cells sharing this node. One advantage of the discrete variational framework is the possibility to include equality and inequality constraints via appropriate modification of the discrete action functional. This is developed in this paper for the treatment of fluid-structure interaction involving a barotropic fluid and an incompressible rubber-like solid. Such a coupling can be realized thanks to the unified geometric framework underlying both the solid and fluid components.

Geometric methods for discrete elasticity have been the subject of various developments. We can mention \cite{ArOr2005} concerning crystal elasticity discretization, where the discrete notions of stress and strain in lattices are introduced, and \cite{Yavari2008}, \cite{AnYa2013}, concerning geometric discretization of elasticity from the point of view of discrete exterior calculus.

Several points of view have been developed to treat fluid-structure interaction problems, see  \cite{Lo2007}, \cite{LeBo2010}, \cite{HoWaLa2012} for the comparison of various fluid-structure interaction methods. We can mention the \textit{immersed boundary method} (see \cite{PeMcQ1989}, \cite{Peskin2002}) where the interactions between the solid boundaries and the fluid are taken into account by means of local body forces, and the \textit{fictitious domain method} (see \cite{GlPaPe1997}) where Lagrange multipliers are used instead of calculating the interaction forces. Other point of views to couple fluids and elastic bodies can be obtained, for example, via the \textit{arbitrary Lagrangian Eulerian method}, see \cite{HiAmCo1997}, techniques based on the \textit{inverse motion function}, see \cite{CoMaMil2008} and \cite{KaRyNa2012}, or the \textit{smoothed particle hydrodynamics method}, see \cite{KhGoFaSh2018}, and references therein.

Our approach is based on the geometric formulation of continuum mechanics seen as a particular classical field theory, see \cite{GiMmsy} and \cite{MaPeShWe2001}. The discrete formulation we develop follows the variational multisymplectic discretization initially developed in \cite{MaPaSh1998,LeMaOrWe2003}. Such an approach has been exploited for instance towards its application to Lie group valued field theories for the treatment of geometrically exact (Cosserat) rods \cite{DeGBKoRa2014,DeGBRa2014} and towards its application to nonsmooth mechanics \cite{DeGBRa2016,DGBDRA2017}.

\medskip

The organization of the paper is as follows. Section \ref{elast_field_th} summarizes briefly the Lagrangian variational setting for elasticity on the continuous side. The elements of the constitutive continuous theory of elasticity necessary for our development are described in rather general terms. Some aspects of the geometric formulation of field theory are also given, including the Noether theorem. In sections \ref{2D_elastic_body} and \ref{3D_elastic_body} we describe, in 2D and 3D, the discrete variational setting for elasticity. We define the discrete deformation gradients, the discrete Cauchy-Green deformation tensors, as well as the discrete Jacobians, in a parallel way with the continuous setting. Based on these definitions, a discrete Lagrangian is constructed for frame indifferent isotropic materials. It is then shown how this approach can be naturally coupled with the variational discretization of barotropic fluids thanks to the introduction of appropriate penalty functions associated to the impenetrability constraints. Illustrative numerical examples are given in 2D and 3D for a barotropic fluid flowing into an hyperelastic container described with the St. Venant-Kirchhoff or the Mooney-Rivlin model.

\section{Elasticity and field theory}\label{elast_field_th}

In this section we briefly review the variational geometric setting for elasticity in the Lagrangian (or material) description by focusing on isotropic hyperelastic materials. The description is given in terms of the (right) Cauchy-Green tensor, which plays a central role in this paper. The formulation is then recasted in a multisymplectic variational setting, which allows to formulate intrinsically the Hamilton principle, the multisymplectic property of the solutions, and the covariant Noether theorem associated to symmetries. This gives the geometric framework to be discretized in a structure preserving way later.

\subsection{Description of hyperelastic bodies}\label{description}

We assume that the elastic body is represented by a compact manifold $ \mathcal{B} \subset \mathbb{R} ^3$ with piecewise smooth boundary, and moves in the ambient space $ \mathcal{M} = \mathbb{R} ^3$. The motion of the body is described by a time dependent configuration map $ \varphi : \mathbb{R} \times \mathcal{B} \rightarrow \mathcal{M} $ which indicates the location $m= \varphi (t,X) \in \mathcal{M}$ at time $t$ of the material point $X \in \mathcal{B} $. The deformation gradient is denoted $ \mathbf{F} (t,X)= \nabla \varphi (t,X)$, given in coordinates by $ \mathbf{F} ^a{}_i= \varphi ^a{}_{,i}$, with $X^i$, $i=1,2,3$ the Cartesian coordinates on $ \mathcal{B} $ and $m^a$, $a=1,2,3$ the Cartesian coordinates on $ \mathcal{M} $. We assume that the configuration map is regular enough so that all the computations below are valid.

A material is called \textit{hyperelastic} if there is a stored energy function $\mathcal{W}$ depending on the points $X \in \mathcal{B}$ and the gradient deformation $\mathbf{F}(t,X)$ such that the first Piola-Kirchhoff stress tensor of the body in the configuration $ \varphi (t,X)$ is $\mathbf{P}= \rho  _0 \frac{\partial \mathcal{W}}{\partial \mathbf{F} }$, where $ \rho  _0$ is the mass density of the body in the Lagrangian description. Given Riemannian metrics $ \mathbf{G} $ and $\mathbf{g}$ on $ \mathcal{B} $ and $ \mathcal{M} $, the stored energy function of an hyperelastic material with configuration map $ \varphi : \mathcal{B} \rightarrow \mathcal{M} $ (at some fixed time) takes the general form
\[
\mathcal{W}( X,\mathbf{F} (X), \mathbf{G}(X),\mathbf{g}( \varphi (X))).
\]
This expression of the stored energy function explicitly written in terms of Riemannian metrics allows to naturally reformulate the axiom of material frame indifference and the isotropic property as invariances of $ \mathcal{W}$ with respect to spatial and material diffeomorphisms, i.e., diffeomorphisms of $\mathcal{M} $ and $ \mathcal{B} $, \cite{MaHu1983}.

From the \textit{axiom of material frame indifference} the stored energy function depends on the deformation gradient only through the \textit{Cauchy-Green deformation tensor} defined by
\begin{equation}\label{Cauchy_Green}
\mathbf{C} =  \varphi^*\mathbf{g}, \qquad \mathbf{C} _{ij}(X)= \mathbf{g}_{ab}( \varphi (X)) \varphi ^a{}_{,i}(X) \varphi ^b{}_{,j}  (X),\end{equation}
i.e., we have $\mathcal{W}( X,\mathbf{F} (X), \mathbf{G}(X),\mathbf{g}( \varphi (X))) = W(X,\mathbf{G}(X),\mathbf{C}(X))$, for some function $W$, where $\varphi^*$ denotes the pull-back of a tensor field by $ \varphi $.

The \textit{material (or Lagrangian) strain tensor} $\mathbf{E}$ and the \textit{second Piola-Kirchhoff stress tensor} $\mathbf{S}$ are defined respectively by 
\begin{equation}\label{material_strain}
\mathbf{E}(\mathbf{C}) =  \frac{1}{2}(\mathbf{C} - \mathbf{G}), \quad  \quad \mathbf{E}_{ij}= \frac{1}{2} \left( \mathbf{C} _{ij} - \mathbf{G} _{ij} \right),
\end{equation}
\begin{equation}\label{sec_Piola_Kir}
\mathbf{S}(\mathbf{C})=2 \rho_0 \frac{\partial W}{\partial \mathbf{C}}, \quad  \quad \mathbf{S}^{ij} = 2 \rho_0 \frac{\partial W}{\partial \mathbf{C} _{ij}}.
\end{equation}
The \textit{elasticity tensor} on the reference configuration is defined by
\begin{equation}
\boldsymbol{\mathsf{C}}= \frac{\partial \mathbf{S}}{\partial \mathbf{C}}, \quad  \quad \boldsymbol{\mathsf{C}}^{ijkl} = \frac{\partial \mathbf{S}^{ij}}{\partial \mathbf{C} _{kl}}= 2 \rho _0\frac{\partial ^2 W}{\partial \mathbf{C} _{ij} \partial \mathbf{C} _{kl}} .
\end{equation}

For simplicity, we will consider the \textit{Euclidean case}, where $\mathbf{G} _{ij} =\delta_{ij}$ and $\mathbf{g}_{ab} =\delta_{ab}$ in Cartesian coordinates, in which case we can write
\[
\mathbf{C}= \mathbf{F} ^\mathsf{T} \mathbf{F} \quad\text{and}\quad \mathbf{C} _{ij} = \varphi ^a{}_{,i} \varphi ^a{}_{,j} = \left\langle \varphi _{,i}, \varphi _{,j} \right\rangle,
\]
where $\langle \cdot  , \cdot  \rangle$ denotes the Euclidean inner product on $\mathbb{R}^n$. 
\medskip

\subsection{Euler-Lagrange equations for an hyperelastic body}

The Lagrangian density of an hyperelastic body evaluated on a configuration map $ \varphi (t,X)$ is
\begin{equation}\label{Lagrangian_hyperelastic}
\begin{aligned} 
\mathcal{L} ( \varphi , \dot \varphi , \nabla \varphi )&= L( \varphi , \dot \varphi , \nabla \varphi ) {\rm d} ^3X \wedge {\rm d} t\\
&=\Big[ \frac{1}{2} \rho  _0 \left\langle \dot \varphi , \dot \varphi \right\rangle - \rho  _0 W(\mathbf{G}, \mathbf{C}) - \rho  _0\Pi ( \varphi)\Big] {\rm d} ^3X \wedge {\rm d} t,
\end{aligned}
\end{equation}
with $ \Pi $ a potential energy term such as the gravitational potential. Note that we didn't write explicitly the dependence on $X$ and $t$. Hamilton's principle reads
\[
\delta \int_0^T\!\!\!\int_ \mathcal{B} L( \varphi , \dot \varphi , \nabla \varphi ) {\rm d} ^3X \wedge {\rm d} t=0,
\]
for variations $ \delta \varphi $ of the body configuration map with $ \delta \varphi =0$ at $t=0,T$. It yields the Euler-Lagrange equations
\begin{equation}\label{EL_general} 
\frac{d}{dt} \frac{\partial L}{\partial \dot \varphi } - \frac{\partial L}{\partial \varphi } = - \operatorname{DIV} \frac{\partial L}{\partial \mathbf{F} } 
\end{equation} 
together with the natural boundary conditions
\begin{equation}\label{NBCond_elastic} 
\frac{\partial L}{\partial \mathbf{F}  ^a{}_{i}} n_i \delta \varphi ^a =0 \quad\text{on}\quad \partial \mathcal{B} 
\end{equation} 
for the allowed variations $ \delta \varphi $ at the boundary, with $ \mathbf{F} = \nabla \varphi $. For the Lagrangian \eqref{Lagrangian_hyperelastic} we compute
\[
\frac{\partial L}{\partial  \mathbf{F}  ^a{}_{i} } =- 2 \rho  _0 \mathbf{F}  ^a{}_{j}\frac{\partial W}{\partial C_{ji}} = - \mathbf{F}  ^a{}_{j} \mathbf{S}^{ji} \quad\text{i.e.}\quad \frac{\partial L}{\partial \mathbf{F} } = - \mathbf{F} \mathbf{S},
\]
hence the Euler-Lagrange equations \eqref{EL_general} for the Lagrangian \eqref{Lagrangian_hyperelastic} gives the equations of motion for an hyperelastic body as
\begin{equation}\label{CEL_hyperelastic}
\rho_0  \ddot{\varphi} = \mathrm{DIV} \left( \mathbf{F} \mathbf{S}\right) - \rho  _0\frac{\partial \Pi }{\partial \varphi}.
\end{equation}
In this paper we focus exclusively on the variational formulation and the equations in the material description, see \cite{GBMaRa2012} for the variational formulation in the Eulerian and convective descriptions and its relation with material frame indifference and material covariance discussed above. 

\subsection{Isotropic hyperelastic materials}\label{cont_isotropic}

A frame indifferent hyperelastic material is isotropic if and only if the stored energy function can be written as $W(X, \mathbf{G} (X), \mathbf{C} (X))= \Phi (X, \lambda _1(X), \lambda _2(X), \lambda_3 (X))$, where $ \Phi $ is a symmetric function of the \textit{principal stretches} $\lambda_1, \lambda_2, \lambda_3$\footnote{Given the polar decomposition $\mathbf{F} =\mathbf{R}\mathbf{U}$ of the gradient deformation, the principal stretches $\lambda_1, \lambda_2, \lambda_3$ are the eigenvalues of $\mathbf{U}$. In particular $\lambda_1^2, \lambda_2^2, \lambda_3^2$ are the eigenvalues of $\mathbf{C} =\mathbf{U}^2$.}. This is equivalent to state that $W$ is a function of the invariants $I_1$, $I_2$, $I_3$ of $\mathbf{C}$ defined by\footnote{For simplicity we didn't write explicitly the dependence of the invariants on $ \mathbf{G} $ and assumed $\mathbf{G}_{ij}= \delta _{ij}$.}
\begin{equation}\label{invariants_C}
\begin{aligned} 
I_1(\mathbf{C}) = \mathrm{Tr}(\mathbf{C}), \qquad I_2(\mathbf{C}) &= \frac{1}{2} \left[(I_1(\mathbf{C}))^2 - \mathrm{Tr}(\mathbf{C}^2) \right] \qquad I_3(\mathbf{C}) = \mathrm{det}(\mathbf{C}).\\
&= \mathrm{det}(\mathbf{C})\mathrm{Tr}(\mathbf{C} ^{-1} ),
\end{aligned}
\end{equation}
The derivative of the invariants with respect to $\mathbf{C}$ are given by
\begin{equation}\label{deriv_invariant}
\frac{\partial I_1}{\partial \mathbf{C}} = \mathbf{G}^{-1} , \qquad  \frac{\partial I_2}{\partial \mathbf{C}} = I_2 \mathbf{C}^{-1} -  I_3 \mathbf{C}^{-2}, \qquad  \frac{\partial I_3}{\partial \mathbf{C}} =I_3  \mathbf{C}^{-1},
\end{equation}
from which the second Piola-Kirchhoff stress tensor \eqref{sec_Piola_Kir} reads
\begin{equation} \label{sec_Piola_Kir2}
\mathbf{S}=2 \rho_0 \left[ \frac{\partial W}{\partial I_1} \mathbf{G}^{-1}  + \left(\frac{\partial W}{\partial I_2}I_2 + \frac{\partial W}{\partial I_3}I_3 \right) \mathbf{C}^{-1} - \frac{\partial W}{\partial I_2}I_3 \mathbf{C}^{-2} \right].
\end{equation}
We refer to \cite[\S3.5]{MaHu1983} for more details.

The Lagrangian density of a frame indifferent isotropic hyperelastic material is therefore
\begin{equation}\label{Lagrangian_hyperelastic_isotropic}
\mathcal{L} ( \varphi , \dot \varphi , \nabla \varphi )=\Big[ \frac{1}{2} \rho  _0 \left\langle \dot \varphi , \dot \varphi \right\rangle - \rho  _0 W(I_1,I_2,I_3) - \rho  _0\Pi ( \varphi)\Big] {\rm d} ^3X \wedge {\rm d} t,
\end{equation} 
while the equations of motion take the form
\begin{equation}\label{CEL_hyperelastic_isotropic}
\rho_0  \ddot{\varphi} = \mathrm{DIV} \left[ 2 \rho_0 \mathbf{F} \left(    \frac{\partial W}{\partial I_1} \mathbf{G} ^{-1}  + \left(\frac{\partial W}{\partial I_2}I_2 + \frac{\partial W}{\partial I_3}I_3 \right) \mathbf{C}^{-1} - \frac{\partial W}{\partial I_2}I_3 \mathbf{C}^{-2}  \right) \right] - \rho  _0\frac{\partial \Pi }{\partial \varphi}.
\end{equation}

Incompressible materials can be treated by inserting the constraint $I_3=1$ in the Hamilton principle as
\[
\delta \int_0^T\!\!\!\int_ \mathcal{B} \big(L( \varphi , \dot \varphi , \nabla \varphi ) + \mu ( I_3 -1 ) \big) {\rm d} ^3X \wedge {\rm d} t=0,
\]
where $ \mu $ is a Lagrange multiplier.

\begin{example}{\rm A well-known example of isotropic hyperelastic model is the \textit{Moonley-Rivlin} model (\cite{Mooney1940}, \cite{Rivlin1948,Rivlin1949a,Rivlin1949b}), which is suitable for the description of certain incompressible rubber-like materials. Its stored energy function is defined by the simple expression
\begin{equation} \label{Mooney_Rivlin}
W(I_1, I_2) = C_1(I_1-3) + C_2(I_2-3),
\end{equation}
with positive material constants $C_1,C_2$.}
\end{example}

\subsection{Multisymplectic variational continuum mechanics}

We briefly recall here the geometric variational framework of classical field theory, as it applies to continuum mechanics (fluid and elasticity), see \cite{MaPeShWe2001} and \cite{LeMaOrWe2003}. In \S\ref{2D_elastic_body} and \S\ref{3D_elastic_body}  we  shall carry out a structure preserving discretization of this setting which allows the identification of the notion of discrete multisymplecticity, discrete momentum map, and discrete Noether theorems.

\paragraph{Configuration bundle, jet bundle, and Lagrangian density.}
The geometric formulation of classical field theories starts with the identification of the \textit{configuration bundle} of the theory, denoted $ \pi _{ \mathcal{Y} , \mathcal{X} }: \mathcal{Y} \rightarrow \mathcal{X} $. The configuration fields $ \varphi $ of the theory are sections of this fiber bundle, i.e., they are smooth maps $ \varphi : \mathcal{X} \rightarrow \mathcal{Y} $ such that $ \pi _{ \mathcal{Y} , \mathcal{X} } \circ \varphi = id_ \mathcal{X} $,  where $ id_ \mathcal{X} $ denotes the identity map on $ \mathcal{X} $.

The \textit{first jet bundle} of the configuration bundle $\pi _{ \mathcal{Y} , \mathcal{X} }: \mathcal{Y} \rightarrow \mathcal{X} $ is the field theoretic analogue of the tangent bundle of classical mechanics. It is defined as the fiber bundle $ \pi _{J^1 \mathcal{Y},\mathcal{Y}  }:  J^1 \mathcal{Y}  \rightarrow \mathcal{Y} $ whose fiber at $y \in \mathcal{Y} $ consists of linear maps $\gamma :T_x \mathcal{X} \rightarrow T_y \mathcal{Y} $ satisfying $T \pi _{ \mathcal{Y} , \mathcal{X} } \circ \gamma = id_{T_x \mathcal{X} }$, where $x= \pi _{ \mathcal{Y} , \mathcal{X} }(y)$.
The derivative of a field $ \varphi $ can be regarded as a section of the fiber bundle $ \pi _{J^1 \mathcal{Y} , \mathcal{X} }= \pi _{ \mathcal{Y} , \mathcal{X} } \circ \pi _{J^1 \mathcal{Y} , \mathcal{Y} }: J^1 \mathcal{Y} \rightarrow \mathcal{X} $, by considering its \textit{first jet extension} $ x \in \mathcal{X} \rightarrow j^1 \varphi (x)= T_ x \varphi \in J^1_{ \varphi (x)} \mathcal{Y} $, with $T_x \varphi :T_x \mathcal{X} \rightarrow T_{ \varphi (x)} \mathcal{Y} $ the tangent map of $ \varphi $.

A \textit{Lagrangian density} is a smooth bundle map $ \mathcal{L} : J^1 \mathcal{Y} \rightarrow \Lambda ^{n+1} \mathcal{X} $ over $ \mathcal{X} $, where $ \Lambda ^{n+1} \mathcal{X} \rightarrow \mathcal{X} $ is the vector bundle of $(n + 1)$-forms on $ \mathcal{X} $, with $ \operatorname{dim} \mathcal{X} =n+1$. The associated \textit{action functional} is given by
\begin{equation}\label{action_funct}
S( \varphi ):=\int_\mathcal{X} \mathcal{L} ( j ^1 \varphi (x)),
\end{equation}
and the Euler-Lagrange equations follow from the stationary condition for $ \delta S( \varphi)=0$ for appropriate boundary conditions.

For continuum mechanics, the configuration bundle is the trivial fiber bundle
\[
\mathcal{Y} = \mathcal{M} \times \mathcal{X} \rightarrow \mathcal{X} , \quad\text{with}\quad \mathcal{X} = \mathbb{R} \times \mathcal{B} ,
\]
where $ \mathcal{B} $ is the reference configuration of the elastic body and $ \mathcal{M} $ is the ambient space. We have the variables $x=(t,X) \in\mathcal{X} $ and $y=(x,m)=(t,X,m) \in \mathcal{Y}$. A section of this bundle is a map $ \varphi : \mathcal{X} \rightarrow \mathcal{X} \times \mathcal{M} $, whose first component is $ id_ \mathcal{X} $. It is canonically identified with a map $ \varphi : \mathcal{X} = \mathbb{R} \times \mathcal{B} \rightarrow \mathcal{M} $ referred to as the elastic body configuration in \S\ref{description}.
The first jet extension of $ \varphi $ is $j^1 \varphi (t,X)= ( \varphi (t,X), \dot \varphi (t,X), \nabla \varphi (t,X))$ and the Lagrangian density reads
\[
\mathcal{L} ( \varphi , \dot \varphi , \nabla \varphi )= L( \varphi , \dot \varphi , \nabla \varphi ) {\rm d} ^3 X \wedge {\rm d} t.
\]

\paragraph{Multisymplectic form formula and Noether theorem.} A solution $ \varphi $ of the Euler-Lagrange equation associated to $ \mathcal{L} :J^1 \mathcal{Y} \rightarrow \Lambda ^{n+1} \mathcal{X} $ satisfies the \textit{multisymplectic form formula} given by
\begin{equation}\label{MSFF} 
\int_{ \partial U} ( j ^1 \varphi ) ^* \mathbf{i} _{j^1V} \mathbf{i} _{j^1W} \Omega _ \mathcal{L} =0,
\end{equation} 
for all open subset $U \subset \mathcal{X} $ with piecewise smooth boundary, and for all vector fields $V,W$ solutions of the first  variation of the Euler-Lagrange equations at $ \varphi $, see \cite{MaPaSh1998}. In this formula, $ \Omega _ \mathcal{L}$ denotes the Cartan $(n+2)$-form on $J^1 \mathcal{Y} $ associated to $ \mathcal{L} $ and  $ j ^1 V,j^1W$ are the first jet extension of $V,W$ to $J^1 \mathcal{Y} $, see \cite{GiMmsy}. Property \eqref{MSFF} is an extension to field theory of the symplectic property of the solution of the Euler-Lagrange equations of classical mechanics.

Consider the action of a Lie group $G$ on $ \mathcal{Y} $, assume that this action covers a diffeomorphism of $ \mathcal{X} $, and assume that the Lagrangian density $ \mathcal{L} $ is $G$-equivariant. Then a solution $ \varphi $ of the Euler-Lagrange equation satisfies the \textit{covariant Noether theorem}
\begin{equation}\label{CNT}
\int_{\partial U} (j^1\varphi )^* J^ \mathcal{L} ( \xi )=0,
\end{equation} 
for all open subset $U \subset \mathcal{X} $ with piecewise smooth boundary and for all $ \xi \in \mathfrak{g} $, the Lie algebra of $G$, with $J^ \mathcal{L} :J^1 \mathcal{Y} \rightarrow \mathfrak{g} ^* \otimes \Lambda ^n J^1 \mathcal{Y} $ the covariant momentum map associated to $ \mathcal{L} $ and to the Lie group action, see \cite{GiMmsy} for details.

The properties \eqref{MSFF} and  \eqref{CNT} admit discrete versions that will be considered later.

\section{2D discrete fluid - elastic body interactions}  \label{2D_elastic_body}

In this section we develop a multisymplectic variational discretization of two dimensional hyperelasticity. A main step in our approach is the definition of the discrete Cauchy-Green tensor associated to a discrete configuration map. We then show how this approach can be naturally coupled with the multisymplectic discretization of barotropic fluids. This is then used for the development of a variational discretization of the coupled dynamics of a barotropic fluid flowing on a hyperelastic body.

\subsection{Multisymplectic discretization and the discrete Cauchy-Green tensor} \label{2D_mult_discret}

\paragraph{Discrete Lagrangian setting.} Consider the configuration bundle of continuum mechanics $ \mathcal{Y} = \mathcal{X} \times \mathcal{M} \rightarrow \mathcal{X} = \mathbb{R} \times \mathcal{B} $. In this section we assume that $ \mathcal{B}$ is a domain in $ \mathbb{R} ^2 $ and take $ \mathcal{M} = \mathbb{R} ^2 $. To discretize the geometric setting, one first considers a \textit{discrete parameter space} $\mathcal{U} _d$ and a \textit{discrete base-space configuration}, which is a one-to-one map
\[
\phi _{ \mathcal{X} _d}: \mathcal{U} _d \rightarrow \phi _{ \mathcal{X} _d}( \mathcal{U} _d) = \mathcal{X} _d \subset \mathcal{X} 
\]
whose image $ \mathcal{X} _d$ is the discrete spacetime. The \textit{discrete fields} are the sections of the \textit{discrete configuration bundle} $ \pi _d: \mathcal{Y} _d= \mathcal{X} _d \times \mathcal{M} \rightarrow \mathcal{X} _d$, which can be identified with maps $ \varphi _d: \mathcal{X} _d \rightarrow \mathcal{M} $. The discrete field and the discrete base-space configuration can be simultaneously described by introducing the \textit{discrete configuration}, which is a map $ \phi _d: \mathcal{U} _d \rightarrow \mathcal{Y} $. From  $\phi _d$ we obtain the discrete base-space configuration and the discrete field as $ \phi _{ \mathcal{X} _d}= \pi _d \circ \phi _d$ and $ \varphi _d= \phi _d \circ \phi _{ \mathcal{X} _d}^{-1} $, see Fig.\,\ref{D_conf_bundle}.

\begin{figure}[H] {
\begin{displaymath}
\begin{xy}
\xymatrix{   &    \mathcal{Y} _d=  \mathcal{X} _d \times \mathcal{M}     \ar@<2pt>[dd]^{\pi_d}  \\ 
&\\
\mathcal{U} _d  \ar[ruu]^-{\phi_d}  \ar[r]_-{\phi_{ \mathcal{X} _d}}  &   \phi _{\mathcal{X} _d}(\mathcal{U} _d) = \mathcal{X} _d\subset \mathcal{X} \ar@<2pt>[uu]^{\varphi_d }  }
\end{xy}
\end{displaymath} }
\caption{\footnotesize Discrete configuration and discrete configuration bundle}\label{D_conf_bundle} 
\end{figure}

\begin{figure}[H] \centering 
\includegraphics[width=2.2 in]{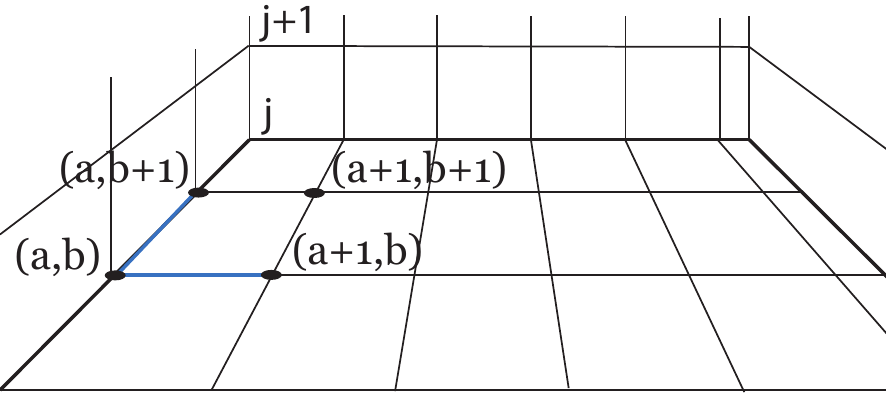} 
\caption{\footnotesize Discrete spacetime domain $\mathcal{U}_d$.} \label{pyramid} 
\end{figure}

We choose the discrete parameter space $\mathcal{U} _d= \{0,...,j,...,N\} \times \mathbb{B}_d$. Here 
$ \{0,...,j,...,N\}$ describes an increasing sequence of time and $\mathbb{B}_d$ encodes the nodes and simplexes of the discretization of $ \mathcal{B} $. Assuming for simplicity that $ \mathcal{B} $ is a rectangle, we consider $\mathbb{B}_d =   \{0, \ldots, A\}\times \{0, \ldots, B\}$ and denote $(j,a,b) \in \mathcal{U} _d$ the element of the discrete parameter space. The latter determines a set of parallelepipeds, denoted $\mbox{\mancube}_{a,b}^j$, defined by the following eight 
pairs of indices
\begin{equation}\label{cube}
\begin{aligned}
&\mbox{\mancube}_{a,b}^j = \big\{ (j,a,b), (j+1,a,b),(j,a+1,b),(j,a, b+1), (j,a+1,b+1), \\
&\hspace{3cm}(j+1,a+1, b), (j+1,a,b+1),  (j+1,a+1,b+1) \big\},
\end{aligned} 
\end{equation}
with $j= 0, ...,N-1$, $a=0, ..., A-1$, $b=0, ..., B-1$, see Fig.\,\ref{pyramid}. We assume that the discrete base-space configuration is of the form
\begin{equation}\label{disc_base_space_conf}
\phi_{\mathcal{X}_d}:\mathcal{U}_d \ni (j,a,b) \mapsto s^j_{a,b} = (t^j,s^j_a,s^j_b) \in \mathcal{X}_d
\end{equation}
and denote by $ \varphi ^j_{a,b}:= \varphi _d(s^j_{a,b})$ the value of the discrete field  at $ s^j_{a,b}$.

Denoting by $\mathcal{U}_d^{\,\mbox{\mancube}}$ the set of all parallelepipeds in the discrete parameter space $ \mathcal{U} _d$, we write $\mathcal{X} _d^{\,\mbox{\mancube}}:=\phi _{ \mathcal{X} _d}\big( \mathcal{U} _d^{\,\mbox{\mancube}}\big)$ and recall that the \textit{discrete version of the first jet bundle} is given by
\begin{equation}\label{discrete_J1} 
J^1 \mathcal{Y} _d:= \mathcal{X} _d^{\,\mbox{\mancube}} \times \underbrace{ \mathcal{M} \times ... \times \mathcal{M}}_{ \text{$8$ times}} \rightarrow \mathcal{X} _d^{\,\mbox{\mancube}}.
\end{equation} 
Given a discrete field $\varphi _d $, its \textit{first jet extension} is the section of \eqref{discrete_J1} defined by
\begin{equation}\label{discrete_jet_extension} 
j^1 \varphi _d( \mbox{\mancube}_{a,b}^j) = \big( \varphi _{a,b}^j, \varphi _{a,b}^{j+1},  \varphi _{a+1,b}^j, \varphi _{a+1,b}^{j+1},  \varphi _{a,b+1}^j ,  \varphi _{a,b+1}^{j+1}, \varphi _{a+1,b+1}^j ,  \varphi _{a+1,b+1}^{j+1} \big),
\end{equation} 
which associates to each parallelepiped, the values of the field at its nodes.
Finally, we recall that a discrete Lagrangian is a map
\[
\mathcal{L} _d : J^1 \mathcal{Y} _d \rightarrow \mathbb{R},
\]
from which the discrete action functional is constructed as
\[
S_d( \varphi _d) = \sum_{\mbox{\mancube}\, \in \mathcal{X} _d^{\mbox{\mancube}}} \mathcal{L} _d\big( j^1 \varphi _d( \mbox{\mancube})\big).
\]

\begin{figure}[H] \centering 
\includegraphics[width=3.8 in]{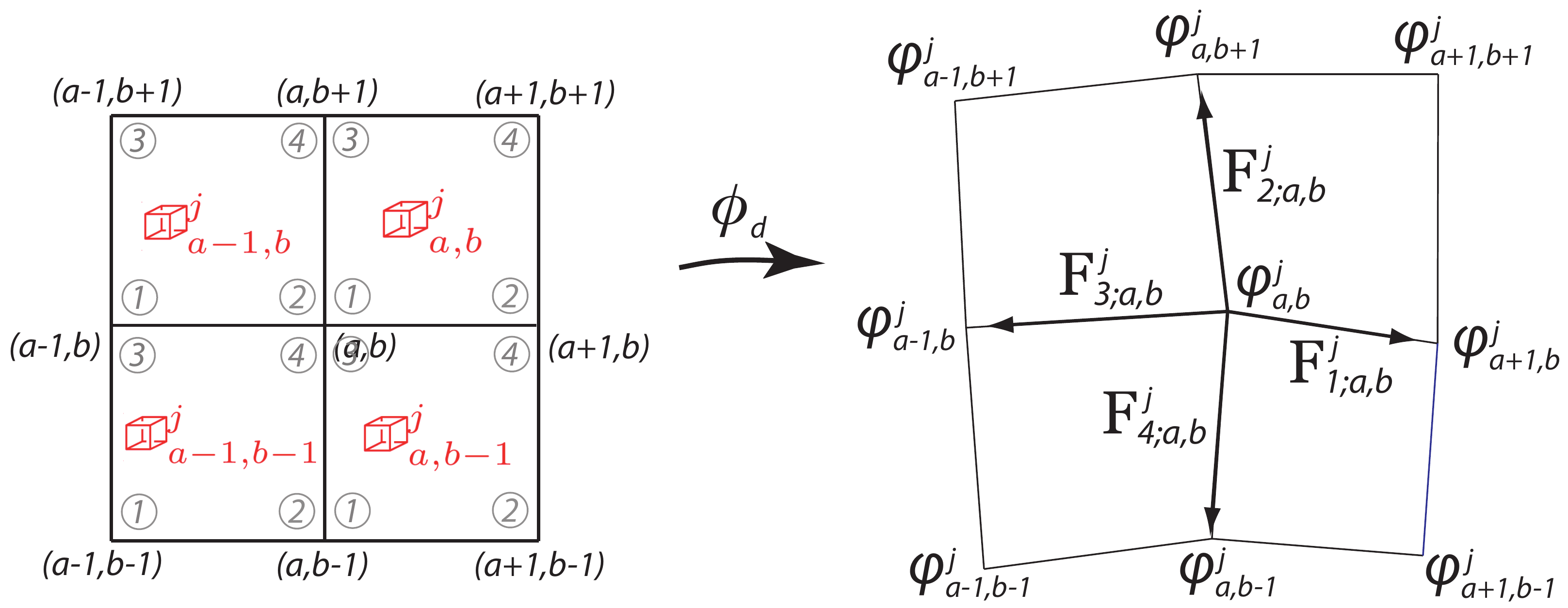} \vspace{-3pt}  
\caption{\footnotesize Discrete field $\phi_d=\varphi_d\circ  \phi _{ \mathcal{X} _d}$ evaluated on $\mbox{\mancube}_{a,b}^j $, $\mbox{\mancube}_{a,b}^j $, $\mbox{\mancube}_{a,b}^j $, $\mbox{\mancube}_{a,b}^j $ at time $t^j$ and the associated vectors $ \mathbf{F} _{\ell;a,b}$, $\ell=1,2,3,4$.} \label{2D_tessellation} 
\end{figure}

\paragraph{Discrete Cauchy-Green deformation tensor.} Given a discrete base-space configuration $ \phi _{ \mathcal{X} _d}$ and a discrete field $ \varphi _d$, the following four vectors $\mathbf{F}_{\ell;a,b}^j \in \mathbb{R} ^2$, $\ell =1,2,3,4$ are defined at each node $(j,a,b)\in \mathcal{U}_d$, see Fig.\,\ref{2D_tessellation} on the right:  
\begin{equation}\label{basis_vectors}
\begin{aligned}
\mathbf{F}_{1;a,b}^j= \frac{\varphi_{a+1,b}^j - \varphi_{a,b}^j}{ |s^j_{a+1} - s^j_{a}| } \quad & \text{and} \quad  \mathbf{F}_{2;a,b}^j = \frac{\varphi_{a,b+1}^j- \varphi_{a,b}^j}{|s^j_{b+1} - s^j_{b}| }
\\
\mathbf{F}_{3;a,b}^j= \frac{\varphi_{a-1,b}^j - \varphi_{a,b}^j}{|s^j_{a} - s^j_{a-1}| } = -  \mathbf{F}_{1;a-1,b}^j \quad &\text{and} \quad \mathbf{F}_{4;a,b}^j = \frac{\varphi_{a,b-1}^j- \varphi_{a,b}^j}{ |s^j_{b} - s^j_{b-1}|} = - \mathbf{F}_{2;a,b-1}^j .
\end{aligned}
\end{equation}
  
\begin{definition}\label{gradient_definition}
The \textbf{discrete deformation gradients} of a discrete field $ \varphi _d$ at the parallelepiped $\mbox{\mancube}_{a,b}^j$ are the four $2 \times 2$ matrices $\mathbf{F}_\ell(\mbox{\mancube}_{a,b}^j)$, $\ell =1,2,3,4$, defined by
\begin{equation}\label{gradient_def}
\begin{aligned}
\mathbf{F}_1(\mbox{\mancube}_{a,b}^j)&= \left[\mathbf{F}_{1;a,b}^j \; \; \mathbf{F}_{2;a,b}^j \right], \qquad & \mathbf{F}_2(\mbox{\mancube}_{a,b}^j) &= \left[\mathbf{F}_{2;a+1,b}^j  \; \; \mathbf{F}_{3;a+1,b}^j \right],\\
\mathbf{F}_3(\mbox{\mancube}_{a,b}^j) &= \left[\mathbf{F}_{4;a,b+1}^j \; \; \mathbf{F}_{1;a,b+1}^j \right], \qquad & \mathbf{F}_4(\mbox{\mancube}_{a,b}^j) &= \left[\mathbf{F}_{3;a+1,b+1}^j  \; \; \mathbf{F}_{4;a+1,b+1}^j \right].
\end{aligned}
\end{equation}
\end{definition}

\medskip 

Note that the discrete deformation gradients are defined at each of the four nodes of $\mbox{\mancube}^j_{a,b}$ that are associated to time $t^j$, see Fig.\,\ref{2D_tessellation} on the left for the ordering $\ell=1,2,3,4$ with respect to the nodes. Also, for each $ \mbox{\mancube}$, the discrete deformation gradients $ \mathbf{F} _\ell( \mbox{\mancube})$, $\ell=1,2,3,4$, depend on the discrete field $ \varphi _d$ only through its first jet extension $j^ 1\varphi _d(\mbox{\mancube})$ at the given $ \mbox{\mancube}$, see \eqref{discrete_jet_extension}.

It is assumed that the discrete field $\varphi_d$ is such that the determinant of the discrete deformation gradient are positive. Using Definition \ref{gradient_definition}, the discrete Cauchy-Green deformation tensors can be defined exactly as in the continuous case.

\begin{definition}\label{CG_definition}
The \textbf{discrete Cauchy-Green deformation tensors} of a discrete field $ \varphi _d$ at the parallelepiped $\,\mbox{\mancube}_{a,b}^j$ are the four $2 \times 2$ symmetric matrices
\begin{equation}\label{def_CG} 
\mathbf{C} _\ell( \mbox{\mancube}_{a,b}^j) = \mathbf{F} _\ell( \mbox{\mancube}_{a,b}^j) ^\mathsf{T}\mathbf{F} _\ell( \mbox{\mancube}_{a,b}^j).
\end{equation} 
\end{definition}
\medskip 

With the assumption on the discrete field mentioned above, the discrete deformation gradients and the discrete Cauchy-Green deformation tensors are maps
\[
\mathbf{F} _\ell: \mathcal{X} _d^{\,\mbox{\mancube}} \rightarrow GL(2), \qquad \mathbf{C} _\ell: \mathcal{X} _d^{\,\mbox{\mancube}} \rightarrow  \operatorname{Sym}_+(2) , \qquad \ell=1,2,3,4,
\]
defined on $ \mathcal{X} _d^{\,\mbox{\mancube}}$, with $ \operatorname{Sym}_+(2)$ the set of symmetric positive definite $ 2 \times 2$ matrices. By using the definitions \eqref{gradient_def} and \eqref{def_CG}, the discrete Cauchy-Green deformation tensors can be explicitly written in terms of the vectors defined in \eqref{basis_vectors} as follows, which is a discrete counterpart of the continuous relation $\mathbf{C}_{ij} = \langle \varphi_{,i} , \varphi_{,j} \rangle$.

\begin{proposition}\label{2D_DCG}
Given a discrete field $ \varphi _d: \mathcal{X}_d \rightarrow \mathcal{M} $, the discrete Cauchy-Green deformation tensors are explicitly given as
\begin{align}
& \mathbf{C}_1(\mbox{\mancube}_{a,b}^j) =
\begin{bmatrix}
\langle\mathbf{F}_{1;a,b}^j , \mathbf{F}_{1;a,b}^j \rangle & \; \langle \mathbf{F}_{1;a,b}^j , \mathbf{F}_{2;a,b}^j\rangle \\[6pt]
\langle \mathbf{F}_{2;a,b}^j , \mathbf{F}_{1;a,b}^j\rangle  & \; \langle\mathbf{F}_{2;a,b}^j,  \mathbf{F}_{2;a,b}^j\rangle \end{bmatrix},\\
&\mathbf{C}_2(\mbox{\mancube}_{a,b}^j) =
\begin{bmatrix} 
\langle\mathbf{F}_{3;a+1,b}^j ,  \mathbf{F}_{3;a+1,b}^j \rangle  & \; \langle\mathbf{F}_{3;a+1,b}^j , \mathbf{F}_{2;a+1,b}^j \rangle  \\[6pt]
\langle\mathbf{F}_{2;a+1,b}^j , \mathbf{F}_{3;a+1,b}^j\rangle &  \; \langle\mathbf{F}_{2;a+1,b}^j, \mathbf{F}_{2;a+1,b}^j \rangle \end{bmatrix}, \nonumber\\
& \mathbf{C}_3(\mbox{\mancube}_{a,b}^j) =
\begin{bmatrix}  
\langle \mathbf{F}_{1;a,b+1}^j ,  \mathbf{F}_{1;a,b+1}^j \rangle & \; \langle \mathbf{F}_{1;a,b+1}^j  , \mathbf{F}_{4;a,b+1}^j \rangle  \\[6pt]
\langle \mathbf{F}_{4;a,b+1}^j , \mathbf{F}_{1;a,b+1}^j\rangle  & \; \langle \mathbf{F}_{4;a,b+1}^j , \mathbf{F}_{4;a,b+1}^j\rangle \end{bmatrix}, \label{2D_disc_Cauchy_Green}\\
&\mathbf{C}_4(\mbox{\mancube}_{a,b}^j) =
\begin{bmatrix}
\langle\mathbf{F}_{3;a+1,b+1}^j , \mathbf{F}_{3;a+1,b+1}^j\rangle & \; \langle\mathbf{F}_{3;a+1,b+1}^j , \mathbf{F}_{4;a+1,b+1}^j \rangle\\[6pt]
\langle\mathbf{F}_{4;a+1,b+1}^j , \mathbf{F}_{3;a+1,b+1}^j \rangle & \; \langle\mathbf{F}_{4;a+1,b+1}^j ,  \mathbf{F}_{4;a+1,b+1}^j\rangle \end{bmatrix}, \nonumber
\end{align}
in terms of the vector $ \mathbf{F}^j_{\ell;a,b}$ defined in \eqref{basis_vectors}. 
\end{proposition}

\paragraph{Discrete Jacobian.} We recall from \cite{DeGB2021} that the discrete Jacobian of a discrete field $ \varphi _d$ at the parallelepiped $\,\mbox{\mancube}_{a,b}^j$ are the four numbers $J_\ell(\mbox{\mancube}_{a,b}^j)$, $\ell =1,2,3,4$, defined by
\begin{equation}\label{disc_rlat_surf_jacob}
\begin{aligned}
J_1(\mbox{\mancube}_{a,b}^j) &=  | \mathbf{F}_{1;a,b}^j \times \mathbf{F}_{2;a,b}^j |,\qquad  & J_2(\mbox{\mancube}_{a,b}^j) &=  | \mathbf{F}_{2;a+1,b}^j \times \mathbf{F}_{3;a+1,b}^j |,\\
J_3(\mbox{\mancube}_{a,b}^j) &=  | \mathbf{F}_{4;a,b+1}^j \times \mathbf{F}_{1;a,b+1}^j |, \qquad & J_4(\mbox{\mancube}_{a,b}^j) &=  | \mathbf{F}_{3;a+1,b+1}^j \times \mathbf{F}_{4;a+1,b+1}^j | .
\end{aligned}
\end{equation}

The discrete Jacobian can be written in terms of the discrete deformation gradient and the discrete Cauchy-Green tensors exactly as in the continuous case, namely, we have the following relations.

\begin{proposition}\label{discrete_2D_Jac} The discrete Jacobian is related to the discrete deformation gradients and discrete Cauchy-Green tensors as follows
\begin{equation}\label{Jacobian_rel_2D}
\mathrm{det} \big( \mathbf{F}_\ell(\mbox{\mancube}_{a,b}^j)\big)=J_\ell(\mbox{\mancube}_{a,b}^j), \qquad \mathrm{det} \big( \mathbf{C}_\ell(\mbox{\mancube}_{a,b}^j)\big) = J_\ell(\mbox{\mancube}_{a,b}^j)^2,
\end{equation}
for each $\mbox{\mancube}_{a,b}^j \in \mathcal{X} _d^{\,\mbox{\mancube}}$ and all $\ell=1,2,3,4$.
\end{proposition}

\subsection{Discrete frame indifferent linear hyperelastic model} \label{disc_hyperelastic}

In the 2D case, we shall focus on the St. Venant-Kirchhoff model
\begin{equation}\label{StVK} 
\rho  _0 W( \mathbf{C} )= \lambda \operatorname{Tr}( \mathbf{E}) ^2 + 2 \mu \operatorname{Tr}( \mathbf{E} ^2 )= \mathsf{E}^ \mathsf{T} \mathsf{C}\,\mathsf{E},
\end{equation}
where $ \lambda $ and $ \mu $ are the first and second Lam\'e coefficients and $ \mathbf{E}$ is the material strain tensor defined in \eqref{material_strain} with $ \mathbf{G} _{ij}= \delta _{ij}$. In the last term we rewrote the stored energy function in a vectorial form by introducing the vector $\mathsf{E}$ and the symmetric matrix $\mathsf{C}$ defined by
\[
\mathsf{E}=( \mathbf{E}_{11}\;\; \mathbf{E}_{22}\;\;2\mathbf{E}_{12})^\mathsf{T} \in \mathbb{R} ^3  
\qquad\text{and}\qquad 
\mathsf{C}=\frac{E}{1-\nu^2}\left[
\begin{matrix} 
1 &\nu & 0 \\
\nu & 1 & 0 \\
0 & 0 & \frac{(1-\nu)}{2} \\
\end{matrix} \right]
\]
with $ \nu=\frac{ \lambda }{ \lambda + 2 \mu } $ the Poisson ratio and $E=\frac{4 \mu ( \lambda + \mu )}{ \lambda + 2 \mu }$ the Young modulus in 2D. 

\paragraph{Discrete stored energy function.} The discrete stored energy function $W_d:J^1 \mathcal{Y} _d \rightarrow \mathbb{R}$ is defined in terms of the continuous one as follows
\begin{equation}\label{2D_stored_energy}
\rho_0 W_d\big(j^1 \varphi _d(\mbox{\mancube})\big) := \frac{\rho_0}{4}  \sum_{\ell=1}^4 W\big( \mathbf{C}_\ell(\mbox{\mancube})\big)=\frac{1}{4}  \sum_{\ell=1}^4 \mathsf{E}_\ell (\mbox{\mancube})^ \mathsf{T} \mathsf{C}\,\mathsf{E}_\ell(\mbox{\mancube}),
\end{equation} 
with $\mathsf{E}_\ell( \mbox{\mancube}) = ( \mathbf{E}_{\ell\, 11}( \mbox{\mancube})\;\; \mathbf{E}_{\ell\, 22}( \mbox{\mancube})\;\;2\mathbf{E}_{\ell\, 12}( \mbox{\mancube}))^\mathsf{T}$. For instance, for $\ell=1$, we have
\[
\mathsf{E}_1( \mbox{\mancube}^j_{a,b})= \Big(\frac{1}{2}( \langle \mathbf{F}_{1;a,b} ,  \mathbf{F}_{1;a,b} \rangle -1), \frac{1}{2}( \langle \mathbf{F}_{2;a,b} ,  \mathbf{F}_{2;a,b} \rangle -1),  \langle \mathbf{F}_{1;a,b} ,  \mathbf{F}_{2;a,b} \rangle\Big).
\]

\paragraph{Discrete second Piola-Kirchhoff stress.}
The discrete second Piola-Kirchhoff stress at time $t^j$ and spatial position $\ell$, with the ordering $\ell= 1$ to $\ell = 4$,  respectively associated to the nodes $(j,a,b)$, $(j,a+1,b)$, $(j,a,b+1)$, $(j,a+1,b+1)$, is defined as
\begin{equation}\label{2D_SPiolaK}
\mathbf{S}_\ell( \mbox{\mancube}^j_{a,b}) = 2 \rho_0 \frac{\partial W}{\partial \mathbf{C}}\big( \mathbf{C}_\ell(\mbox{\mancube}^j_{a,b})\big),
\end{equation}
for all $\mbox{\mancube}^j_{a,b} \in \mathcal{X} _d^{\,\mbox{\mancube}}$.


\subsection{2D Barotropic fluid - elastic body interactions}

We describe here the discrete variational setting for the coupled dynamics of  a barotropic fluid flowing on a hyperelastic body. 
The Lagrangian density of a barotropic fluid is of the form
\begin{equation}\label{Lagr_fluid}
\begin{aligned} 
\mathcal{L} ^\mathsf{f} ( \varphi , \dot \varphi , \nabla \varphi )&= L^\mathsf{f}( \varphi , \dot \varphi , \nabla \varphi ) {\rm d} ^2X \wedge {\rm d} t\\
&=\Big[ \frac{1}{2} \rho  _0^\mathsf{f} \left\langle \dot \varphi , \dot \varphi \right\rangle - \rho  _0^\mathsf{f}  W^\mathsf{f}( \rho  _0^\mathsf{f} , J) - \rho  _0^\mathsf{f} \Pi ( \varphi)\Big] {\rm d} ^2X \wedge {\rm d} t,
\end{aligned}
\end{equation}
with $ \rho  _0^\mathsf{f} $ the fluid mass density and $J= \operatorname{det} \mathbf{F} = \sqrt{ \operatorname{det} \mathbf{C} }$, where we assume $ \mathbf{G} _{ij}= \delta _{ij}$ and $ g_{ab}= \delta _{ab}$. Here $ W^\mathsf{f}( \rho  _0^\mathsf{f} , J) = w ( \rho  _0^\mathsf{f} /J)$, with $w( \rho  )$ the specific internal energy barotropic fluid. We recall that Hamilton's principle yields the barotropic fluid equations in the Lagrangian (or material) description as
\begin{equation}\label{CEL_barotropic}
\rho_0^\mathsf{f}   \ddot{\varphi} = - \frac{\partial }{\partial x^i} \big(P_W J \mathbf{F} ^{-1} \big) ^i - \rho  _0^\mathsf{f} \frac{\partial \Pi }{\partial \varphi},
\end{equation}
with $P_W( \rho_0  ^\mathsf{f} , J)= - \rho_0^\mathsf{f}   \frac{\partial W^\mathsf{f}}{\partial J} ( \rho  _0^\mathsf{f} ,J)$ the pressure in the material description.

\paragraph{Discrete Lagrangians for the fluid and the elastic body.} For the elastic body, we consider discrete Lagrangians of the form
\begin{equation}\label{Discrete_Lagrangian_2D_wave}
\mathcal{L}_{d}^\mathsf{e}\big(j^1 \varphi _d(\mbox{\mancube})\big) = \operatorname{vol}(\mbox{\mancube}) \Big( \rho  _0^\mathsf{e} K_d\big(j^1 \varphi _d(\mbox{\mancube})\big)  - \rho_0^\mathsf{e} W_d^\mathsf{e}\big(j^1 \varphi _d(\mbox{\mancube})\big) - \rho  _0^\mathsf{e}\Pi_d\big(j^1 \varphi _d(\mbox{\mancube})\big) \Big),
\end{equation}
where $\operatorname{vol}(\mbox{\mancube}) $ is the volume of the parallelepiped $\mbox{\mancube} \in \mathcal{X} _d^{\,\mbox{\mancube}}$. For simplicity, we have assume that $ \rho  _0^\mathsf{e}$ is constant in space.
The discrete kinetic energy $K_d: J^1 \mathcal{Y} _d \rightarrow  \mathbb{R}$ is defined by
\begin{equation}\label{kinetic_energy}
 K_d\big(j^1 \varphi _d(\mbox{\mancube}_{a,b}^j)\big)  := \frac{1}{4} \sum_{\alpha=a}^{a+1} \sum_{\beta=b}^{b+1} \frac{1}{2} \left| \mathbf{v}_{\alpha,\beta}^j\right|^2 ,
\end{equation}
with $\mathbf{v}_{\alpha,\beta}^j= (\varphi_{\alpha,\beta}^{j+1}- \varphi_{\alpha,\beta}^j)/\Delta t^j$.
The discrete stored energy function $W_d^\mathsf{e}:J^1\mathcal{Y}_d \rightarrow \mathbb{R}$ is defined by
\begin{equation}\label{2D_stored_energy}
W_d^\mathsf{e}\big(j^1 \varphi _d(\mbox{\mancube}_{a,b}^j)\big) := \frac{1}{4}  \sum_{\ell=1}^4 W^\mathsf{e}\big( \mathbf{C}_\ell(\mbox{\mancube}_{a,b}^j)\big),
\end{equation}
with $W^\mathsf{e}$ the stored energy function of the hyperelastic body.
Finally, the discrete potential energy $\Pi_d:J^1\mathcal{Y}_d \rightarrow \mathbb{R}$ is defined by
\begin{equation}\label{grav_potential}
\Pi_d\big(j^1 \varphi _d(\mbox{\mancube}_{a,b}^j)\big):= \frac{1}{4}  \sum_{\alpha=a}^{a+1} \sum_{\beta=b}^{b+1}  \Pi ( \varphi ^j_{ \alpha , \beta }),
\end{equation} 
where $ \Pi   $ is the potential energy of the continuous model. We shall focus on the
gravitation potential $ \Pi ( \varphi )= \left\langle  \mathbf{g} , \varphi \right\rangle $, with $ \mathbf{g}$ the gravitational acceleration vector.

For the barotropic fluid, we shall consider discrete versions of \eqref{Lagr_fluid} given in a similar way with the elasticity case \eqref{Discrete_Lagrangian_2D_wave} by 
\begin{equation}\label{Discrete_Lagrangian_2D_fluid}
\mathcal{L}_{d}^\mathsf{f}\big(j^1 \varphi _d(\mbox{\mancube})\big) = \operatorname{vol}(\mbox{\mancube}) \Big( \rho  _0 ^\mathsf{f}K_d\big(j^1 \varphi _d(\mbox{\mancube})\big)  - \rho_0 ^\mathsf{f}W_d^\mathsf{f}\big(\rho  _0^\mathsf{f},j^1 \varphi _d(\mbox{\mancube})\big) - \rho  _0^\mathsf{f}\Pi_d\big(j^1 \varphi _d(\mbox{\mancube})\big) \Big),
\end{equation}
where now $W^\mathsf{f}_d: J^1 \mathcal{Y} _d \rightarrow \mathbb{R} $ is the discrete internal energy of the fluid given by
\begin{equation}\label{2D_internal_energy}
W_d^\mathsf{f}\big(\rho  _0^\mathsf{f}, j^1 \varphi _d(\mbox{\mancube}_{c,d}^j)\big) := \frac{1}{4}  \sum_{\ell=1}^4 W^\mathsf{f}\big( \rho  _0^\mathsf{f}, J_\ell(\mbox{\mancube}_{c,d}^j)\big),
\end{equation}
with $W^\mathsf{f}(\varrho _0^\mathsf{f}, J)$ the internal energy of the barotropic fluid and $J_\ell(\mbox{\mancube}_{c,d}^j)$ the discrete Jacobians associated to $\mbox{\mancube}_{c,d}^j$. We assumed that $ \rho  _0^\mathsf{f}$ is constant in space for simplicity.

\paragraph{Discrete action and variations.} To simplify the exposition, we assume that the discrete base space configuration for the elastic body is fixed and given by $ \phi _ { \mathcal{X} _d}(j,a,b)= (j \Delta  t, a \Delta s_1, b \Delta s_2)$, for given $ \Delta t$, $ \Delta s_1$, $ \Delta s_2$, similarly for the fluid. In this case, we have $ \operatorname{vol}(\mbox{\mancube}) = \Delta t \Delta s_1 \Delta s_2$ in the discrete Lagrangian and the mass of each 2D cell in $\phi _ { \mathcal{X} _d}(\mathbb{B}_d)$ is $M^\mathsf{e} = \rho  _0^\mathsf{e} \Delta s_1 \Delta s_2$.

We consider the discrete action functional
\begin{equation}\label{Disc_act_sum_2D}
\begin{aligned}
S_d( \varphi _d)&= \sum_{\mbox{\mancube} \,\in \mathcal{X} _d^{\,\mbox{\mancube}}} \mathcal{L} _d^\mathsf{e}\big(j^1 \varphi _d( \mbox{\mancube})\big) +\mathcal{L} _d^\mathsf{f}\big(j^1 \varphi _d( \mbox{\mancube})\big) \\
&= \sum_{j=0}^{N-1} \left( \sum_{a=0}^{A-1} \sum_{b=0}^{B-1} \mathcal{L}^\mathsf{e}_{d}\big(j^1 \varphi _d(\mbox{\mancube}_{a,b}^j)\big) + \sum_{c=0}^{C-1} \sum_{d=0}^{D-1} \mathcal{L}^\mathsf{f}_{d}\big(j^1 \varphi _d(\mbox{\mancube}_{c,d}^j)\big)\right),
\end{aligned} 
\end{equation}
where $\mathcal{L}^\mathsf{e}_{d}$ and $\mathcal{L}^\mathsf{f}_{d}$ are the discrete Lagrangians associated to the elastic body and to the fluid.
For notational simplicity we have used the same notation $ \varphi _d$ for the both the discrete configuration of the fluid and the elastic body.

Let us denote by $D_k$, $k=1,...,8$, the partial derivative of the discrete Lagrangians with respect to the $k$-th component of $j^1 \varphi _d (\mbox{\mancube} )$, in the order given in \eqref{discrete_jet_extension}. From the discrete Hamilton principle $ \delta S_d( \varphi _d)=0$, we get the discrete Euler-Lagrange equations in the general form
\begin{equation}\label{discrete_Euler_Lagrange} 
\begin{aligned}
&D _1 \mathcal{L} ^{j}_{a,b} + D _2 \mathcal{L} ^{j-1}_{a,b} + D _3 \mathcal{L} ^{j}_{a-1,b} + D _4\mathcal{L} ^{j-1}_{a-1,b} \\
& \qquad \qquad \qquad + D _5 \mathcal{L} ^{j}_{a,b-1} + D _6 \mathcal{L} ^{j-1}_{a,b-1} + D _7 \mathcal{L} ^{j}_{a-1,b-1} + D _8 \mathcal{L} ^{j-1}_{a-1,b-1}  =0
\end{aligned}
\end{equation}
for both $\mathcal{L}^\mathsf{e}_{d}$ and $\mathcal{L}^\mathsf{f}_{d}$, where we used the abbreviate notation $ \mathcal{L} ^j_{a,b}:= \mathcal{L} (j^1 \varphi _d (\mbox{\mancube}^j_{a,b}))$.
By using the expressions of the discrete Lagrangians given in \eqref{Discrete_Lagrangian_2D_wave} and \eqref{Discrete_Lagrangian_2D_fluid}, the general form \eqref{discrete_Euler_Lagrange} yields
\begin{equation}\label{concrete_EL}
\begin{aligned} 
& \rho  _0^\mathsf{n} \frac{\mathbf{v}^{j-1}_{a,b}- \mathbf{v}^j_{a,b}}{ \Delta t} -  \frac{1}{4}\sum_{\ell=1}^4 \Big[ \mathfrak{D}_1^\ell   \left( \mathbf{F} _\ell(\mbox{\mancube}_{a,b}^j) \mathbf{S}_\ell(\mbox{\mancube}_{a,b}^j)\right)+\mathfrak{D}_3^\ell   \left( \mathbf{F} _\ell(\mbox{\mancube}_{a-1,b}^j) \mathbf{S}_\ell(\mbox{\mancube}_{a-1,b}^j)\right)\\
& \hspace{2cm} +\mathfrak{D}_5^\ell   \left( \mathbf{F} _\ell(\mbox{\mancube}_{a,b-1}^j) \mathbf{S}_\ell(\mbox{\mancube}_{a,b-1}^j)\right)+\mathfrak{D}_7^\ell   \left( \mathbf{F} _\ell(\mbox{\mancube}_{a-1,b-1}^j) \mathbf{S}_\ell(\mbox{\mancube}_{a-1,b-1}^j)\right) \Big] \\
& \hspace{2cm}- \rho  _0^\mathsf{n} \frac{\partial \Pi}{\partial \varphi }    ( \varphi ^j _{a,b})=0,
\end{aligned}
\end{equation}  
$\mathsf{n}=\mathsf{e},\mathsf{f}$, see Appendix \ref{discrete_EL_derivation} for  the derivation, as well as for the definition of the operators $\mathfrak{D}_k$, $k=1,3,5,7$. Note that these are the discrete equations of motion corresponding to \eqref{CEL_hyperelastic}. In \eqref{concrete_EL}, $\mathbf{S}_\ell$ is the discrete second Piola-Kirchhoff tensor of the model considered, see \eqref{2D_SPiolaK}.  For the St. Venant-Kirchhoff model one has $\mathbf{S}=\lambda \mathrm{Tr}( \mathbf{E}) \mathbf{l} + 2\mu \mathbf{E}$ while for the barotropic fluid $\mathbf{S}= - P_WJ \mathbf{C} ^{-1} $.

Implicit versions of these Euler-Lagrange equations can be derived by considering the appropriate modifications in the definition of the discrete deformation gradients, see \cite[\S3.2.1]{DeGB2021}.  We also refer to that paper for the discrete boundary conditions emerging from the discrete Hamilton principle.

\paragraph{Discrete multisymplectic form formula.} By evaluating the differential $ \mathbf{d} S_d$ of the discrete action functional, restricted to an arbitrary subdomain $ \mathcal{U}_d'$ of $ \mathcal{U} _d$, on a solution $ \varphi _d$ of the discrete Euler-Lagrange equation and by taking the exterior derivative of the resulting expression along the first variations of this solution, one obtains a discrete analogue of the multisymplectic form formula \eqref{MSFF}, see \cite{MaPaSh1998} for a detailed treatment. This formula extends to spacetime discretization, the symplectic property of variational time integrators \cite{MaWe2001}. It also encodes a discrete version of the reciprocity theorem of continuum mechanics, as well as discrete time symplecticity of the solution, see  \cite{LeMaOrWe2003}. We refer to \cite{DeGB2021} for application to the case of the barotropic fluid, the case of the elastic body being similar. In the present situation, the solution of the discrete Euler-Lagrange equations \eqref{concrete_EL} satisfies the discrete multisymplectic form formula on each subdomain $ \mathcal{U} _d'$ of the discrete spacetime associated to either the fluid motion or the elastic body motion.

\paragraph{Symmetries and discrete Noether theorems.}  Consider the action of the special Euclidean group $SE(2)\ni ( R_ \theta , u)$ on a discrete field  given by $\varphi _d \mapsto R_ \theta \varphi _d + u $ with $R_ \theta $ the rotation of angle $ \theta $. The \textit{discrete covariant momentum maps} $J^\mathsf{p}_{ \mathcal{L} _d}: J^1 \mathcal{Y} _d \rightarrow \mathfrak{se}(2) ^* $, $\mathsf{p}=1,...,8$, associated to the $SE(2)$ action and to a discrete Lagrangian $ \mathcal{L} _d: J^1 \mathcal{Y} _d \rightarrow \mathbb{R} $, are found as
\begin{equation}\label{Momap_SE2} 
J^\mathsf{p}_{ \mathcal{L} _d}\big( j^1 \varphi _d(\mbox{\mancube}_{a,b}^j)\big) = \big( \varphi ^{(\mathsf{p})} \times D_ \mathsf{p} \mathcal{L} ^j_{a,b}, D_ \mathsf{p} \mathcal{L} ^j_{a,b}\big) \in \mathfrak{se}(2) ^* .
\end{equation} 
In \eqref{Momap_SE2}, $\varphi ^{(\mathsf{p})}$ is the $\mathsf{p}$-th component of $j^1 \varphi _d$, see \eqref{discrete_jet_extension}, $D_\mathsf{p}$ denotes the partial derivative with respect to the $\mathsf{p}$-th argument of $ \mathcal{L} _d$, and $\mathcal{L} _{a,b}^j:= \mathcal{L} _d\big( j^1 \varphi _d(\mbox{\mancube}_{a,b}^j)\big)$.
From the discrete covariant momentum maps \eqref{Momap_SE2}, the \textit{discrete classical momentum map} is given by
\begin{align*} 
\mathbf{J} _d ( \boldsymbol{\varphi} ^j, \boldsymbol{\varphi} ^{j+1}) &= \sum_{a=0}^{A-1}\sum_{b=0}^{B-1} \left( \mathcal{J} ^2_{ \mathcal{L} _d} +\mathcal{J} ^4_{ \mathcal{L} _d} +\mathcal{J} ^6_{ \mathcal{L} _d}+\mathcal{J} ^8_{ \mathcal{L} _d}\right)\big( j^1 \varphi _d(\mbox{\mancube}_{a,b}^j)\big) \\
&= -  \sum_{a=0}^{A-1}\sum_{b=0}^{B-1} \left( \mathcal{J} ^1_{ \mathcal{L} _d} +\mathcal{J} ^3_{ \mathcal{L} _d} +\mathcal{J} ^5_{ \mathcal{L} _d}+\mathcal{J} ^7_{ \mathcal{L} _d}\right)\big( j^1 \varphi _d(\mbox{\mancube}_{a,b}^j)\big) \in \mathfrak{se}(2) ^* 
\end{align*} 
with $ \boldsymbol{\varphi} ^j$ the collection of all positions at time $t^j$, see \cite{DeGBRa2014}. For the discrete Lagrangians \eqref{Discrete_Lagrangian_2D_wave} and \eqref{Discrete_Lagrangian_2D_fluid}, 
we get the expressions
\begin{equation}\label{2D_momentum_map_fluid}
\mathbf{J} _d ( \boldsymbol{\varphi} ^j, \boldsymbol{\varphi} ^{j+1}) = \sum_{a=0}^{A-1}\sum_{b=0}^{B-1} \mathbf{J} ^\mathsf{e}\big(j^1 \varphi _d (\mbox{\mancube}_{a,b}^j)\big) +\sum_{c=0}^{C-1}\sum_{d=0}^{D-1} \mathbf{J} ^\mathsf{f}\big(j^1 \varphi _d (\mbox{\mancube}_{c,d}^j)\big)
\end{equation}
with
\begin{equation}
\mathbf{J} ^\mathsf{e}\big(j^1 \varphi _d (\mbox{\mancube}_{a,b}^j)\big)= \sum_{\alpha =a}^{a+1} \sum_{\beta =b}^{b+1} \left( \varphi_{\alpha,\beta}^j \times \frac{M^\mathsf{e}}{4} v_{\alpha,\beta}^j ,  \frac{M^\mathsf{e}}{4} v_{\alpha,\beta}^j\right) 
\end{equation}
similarly for the fluid.

In absence of gravitational potential, the discrete Lagrangians \eqref{Discrete_Lagrangian_2D_wave} and \eqref{Discrete_Lagrangian_2D_fluid} are $SE(2)$ invariant. Indeed, when $\varphi _d \mapsto R_ \theta \varphi _d + u $ the associated discrete deformation gradients \eqref{gradient_def} transform as $ \mathbf{F} _\ell \mapsto R_ \theta \mathbf{F} _\ell$, so that the discrete Cauchy-Green deformation tensors \eqref{def_CG} and the discrete Jacobians \eqref{disc_rlat_surf_jacob} are $SE(2)$-invariant. The discrete covariant Noether theorem for the discrete covariant momentum maps $J^\mathsf{p}_{ \mathcal{L} _d}$, given as
\[
\sum_{\mbox{\mancube}\, \in {\mathcal{U} '}_d^{\,\mbox{\mancube}}}\, \sum_{\mathsf{p}; \,\mbox{\mancube}^{(\mathsf{p})} \in \partial \mathcal{U} _d'}J^\mathsf{p}_{ \mathcal{L} _d}(\mbox{\mancube})=0,
\]
thus holds on arbitrary subdomain $\mathcal{U} _d'$ of the discrete spacetime associated to either the fluid motion or the elastic body motion. Furthermore, the classical discrete Noether theorem holds for the discrete classical momentum map \eqref{3D_momentum_map_fluid} as
\[
\mathbf{J} _d ( \boldsymbol{\varphi} ^j, \boldsymbol{\varphi} ^{j+1}) = \mathbf{J} _d ( \boldsymbol{\varphi} ^{j-1}, \boldsymbol{\varphi} ^j), 
\] 
for all $t^j$.

When the gravitational potential is taken into account, only the component of the discrete momentum maps associated to horizontal translations satisfies the discrete Noether theorems above.

\paragraph{Incompressibility and impenetrability conditions.} The constraints will be included in the discrete variational formulation by the addition of appropriate penalty terms.

The first constraint is the incompressibility of the elastic body which is imposed via the condition $J_\ell(\mbox{\mancube})- 1= 0$, for all $\mbox{\mancube} \in \mathcal{X} _d^{\, \mbox{\mancube}}$ and all $\ell=1,2,3,4$, which mimics the properties of rubber. The associated quadratic penalty function is given by
\begin{equation}\label{2D_disc_penalty_function}
\Phi_{\rm in}\big(j^1 \varphi _d (\mbox{\mancube})\big) := \frac{1}{4}  \sum_{\ell=1}^4 \frac{r}{2} \big(J_\ell(\mbox{\mancube})- 1  \big)^2,
\end{equation} 
where $r$ is the penalty parameter.

The second constraint is the impenetrability condition between the fluid node $\varphi_{c,d}^j$ flowing above the elastic body nodes $\varphi_{a,b}^j$, $\varphi_{a+1,b}^j$,  which is defined as follows
\begin{equation}\label{imp_cond}
\Psi_{\rm im}(\varphi_{a,b}^j, \varphi_{a+1,b}^j, \varphi_{c,d}^j) =  \langle (\varphi_{c,d}^j -\varphi_{a+1,b}^j), \mathrm{R}_{\pi/2} (\varphi_{a+1,b}^j -\varphi_{a,b}^j)  \rangle \geq 0,
\end{equation}
where $\mathrm{R}_{\pi/2}$ is a rotation of $+\pi/2$. In this expression $\mathrm{R}_{\pi/2} (\varphi_{a+1,b}^j -\varphi_{a,b}^j) $ represents the outward pointing normal vector to the body. See, e.g., \cite{CiWe2005} for  the use of such constraint functions in impact problems. The function $ \Psi  _{\rm im}$ in \eqref{imp_cond} is evaluated on all the boundary nodes $\varphi_{c,d}^j$ of the fluid and on the two boundary nodes of the elastic body that are the closest to $ \varphi_{c,d}^j$.
The associated quadratic penalty term is given by
\begin{equation}\label{penalty_2D} 
\Phi_{\rm im}(\varphi_{a,b}^j, \varphi_{a+1,b}^j, \varphi_{c,d}^j) = \frac{1}{2} K \Psi_{\rm im}(\varphi_{a,b}^j, \varphi_{a+1,b}^j, \varphi_{c,d}^j)^2,
\end{equation} 
with $K \in \; ]0, \infty[\;$ if $\Phi_{\rm im}(\varphi_{a,b}^j, \varphi_{a+1,b}^j, \varphi_{c,d}^j)\leq 0$ and $K=0$ if $\Phi_{\rm im}(\varphi_{a,b}^j, \varphi_{a+1,b}^j, \varphi_{c,d}^j)>0$.
The derivatives of $\Phi_{\rm im}$ with respect to the nodes $\varphi_{a,b}^j, \varphi_{a+1,b}^j, \varphi_{c,d}^j$ are given by
\begin{equation}\label{deriv_2D_constraint}
\begin{aligned}
&D_1\Phi_{\rm im}(\varphi_{a,b}^j, \varphi_{a+1,b}^j, \varphi_{c,d}^j)
 =-K\Psi_{\rm im}(\varphi_{a,b}^j, \varphi_{a+1,b}^j, \varphi_{c,d}^j)  \mathrm{R}^\top_{\pi/2} (\varphi_{c,d}^j -\varphi_{a+1,b}^j)
\\
&D_2\Phi_{\rm im}(\varphi_{a,b}^j, \varphi_{a+1,b}^j, \varphi_{c,d}^j)=K\Psi_{\rm im}(\varphi_{a,b}^j, \varphi_{a+1,b}^j, \varphi_{c,d}^j) \mathrm{R}_{\pi/2}( \varphi ^j_{a,b}- \varphi ^j_{c,d})\\
&D_3\Phi_{\rm im}(\varphi_{a,b}^j, \varphi_{a+1,b}^j, \varphi_{c,d}^j)
 =K\Psi_{\rm im}(\varphi_{a,b}^j, \varphi_{a+1,b}^j, \varphi_{c,d}^j) \mathrm{R}_{\pi/2} (\varphi_{a+1,b}^j -\varphi_{a,b}^j),
\end{aligned}
\end{equation}
which express the directions of the reactions forces.

The resulting discrete action functional is found as
\begin{equation}\label{tilde_S}
\begin{aligned} 
\widetilde{S}_d( \varphi _d ) &= S_d( \varphi  _d ) - \sum_{j=0}^{N-1}\sum_{a=0}^{A-1}\sum_{j=0}^{B-1} \Delta t \Delta s_1 \Delta s_2\,\Phi_{\rm in}\big(j^1 \varphi _d (\mbox{\mancube}_{a,b}^j)\big) \\
&\hspace{4cm}  - \sum_{j=0}^{N-1}\sum_{(c,d)\, \in \,\mathbb{B}^\mathsf{f}}\Phi_{\rm im}(\varphi_{a,b}^j, \varphi_{a+1,b}^j, \varphi_{c,d}^j),
\end{aligned}
\end{equation} 
where $S_d$ is associated to the discrete Lagrangians of the fluid and the elastic body, see \eqref{Disc_act_sum_2D}. In the third term the second sum is taken over the boundary nodes of the fluid and, given such a boundary node $ \varphi _{c,d}^j$, the nodes $ \varphi_{a,b}^j, \varphi_{a+1,b}^j$ are chosen as explained above.
For the boundary nodes $ \varphi _{c,d}^j$ that are not in contact, the corresponding term vanishes by definition of $\Phi _{\rm im}$. Away from the interface between the barotropic fluid and the elastic body, the boundary nodes are subject to the ambient pressure. This boundary condition directly arises from the variations of $S_d$ as explained below. Parts of the boundary of the body can also be prescribed, see below.

In order to avoid the possible instabilities related to the presence of a large penalty parameter, one could use the methods of augmented Lagrangians, in which both a Lagrange multiplier term and a penalty term are used, see
\cite{Ro1993}. From a practical point of view, see \cite{NoWr2006}, this approach is however much more expensive as it needs to determine the Lagrange multiplier at each contact point.

\subsection{Barotropic fluid flowing over a St. Venant-Kirchhoff hyperelastic container}\label{2D_num_simul}

We consider an incompressible elastic body with stored energy given by the St. Venant-Kirchhoff model \eqref{StVK} and a compressible barotropic fluid described by the Tait equation
\begin{equation}\label{Dis_Tait_eq} 
W_d^\mathsf{f}\big( \rho  _0, j^1 \varphi _d( \mbox{\mancube})\big) =  \frac{1}{4}  \sum_{\ell=1}^4\left[ \frac{A}{ \gamma -1} \left(\frac{J_\ell (\mbox{\mancube})}{ \rho  _0}  \right)    ^{1- \gamma } + B  \left(\frac{J_\ell (\mbox{\mancube})}{ \rho  _0}  \right) \right].
\end{equation}
We shall use the expression \eqref{Dis_Tait_eq}  for the treatment of an isentropic perfect fluid, where the value of the constant $B\neq 0$ does not affect the dynamics, while it allows to naturally impose the boundary condition $P|_{ \partial \mathcal{B} }=B$, with $P=A \left( \frac{\rho  _0}{J} \right) ^ \gamma  $ and $B$ the external pressure. This is crucial for the discretization, since it allows to find the appropriate discretization of the boundary condition directly from the boundary terms of the discrete variational principle, see \cite{DeGB2021}.

\begin{figure}[H] \centering 
\includegraphics[width=3.2 in]{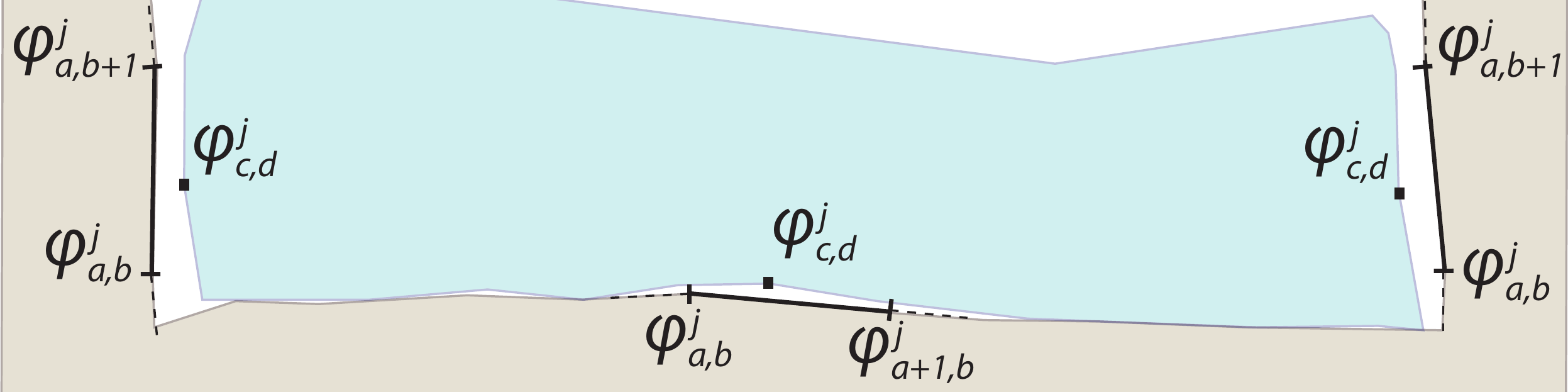} \vspace{-3pt}  
\caption{\footnotesize Elastic body and fluid configurations respectively in nodes $\varphi_{a,b}^j$, $\varphi_{a+1,b}^j$, $\varphi_{a,b+1}^j$, and $\varphi_{c,d}^j$ at time $t^j$.} \label{impenetrability} 
\end{figure}

\begin{figure}[H] \centering 
\includegraphics[width=2.5 in]{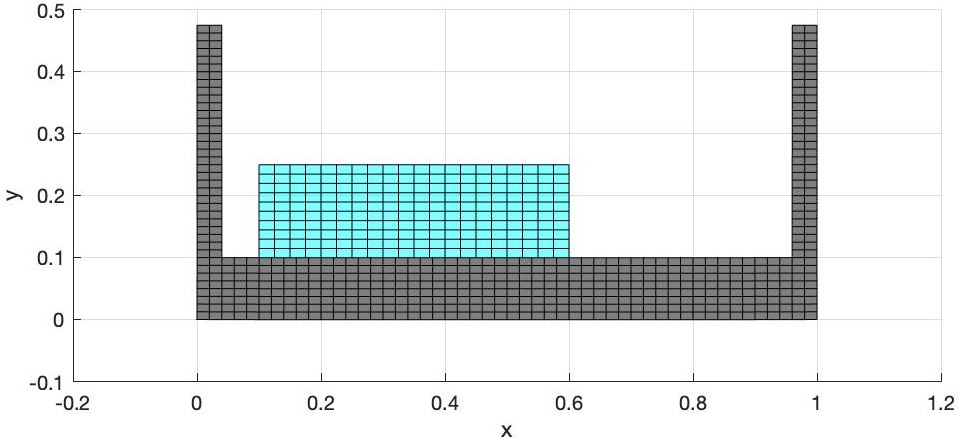} \hspace{0.15cm} \qquad  \qquad \includegraphics[width=1.65 in]{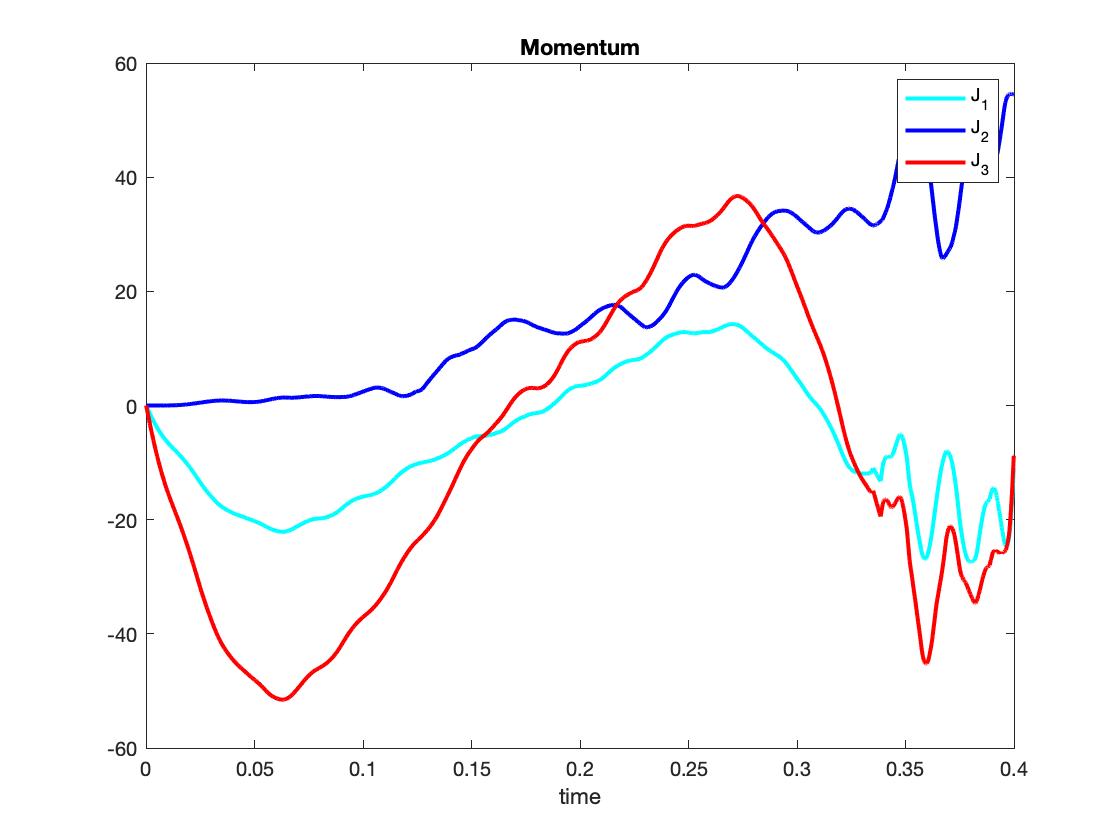} \vspace{-3pt}\\
\hspace{0.6cm}  \includegraphics[width=2.2 in]{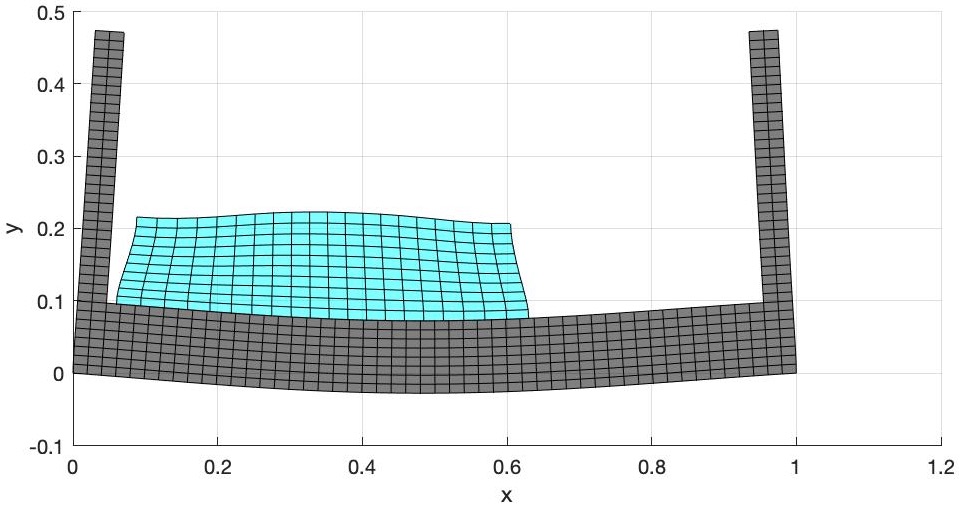}  \quad \qquad \qquad \includegraphics[width=1.65 in]{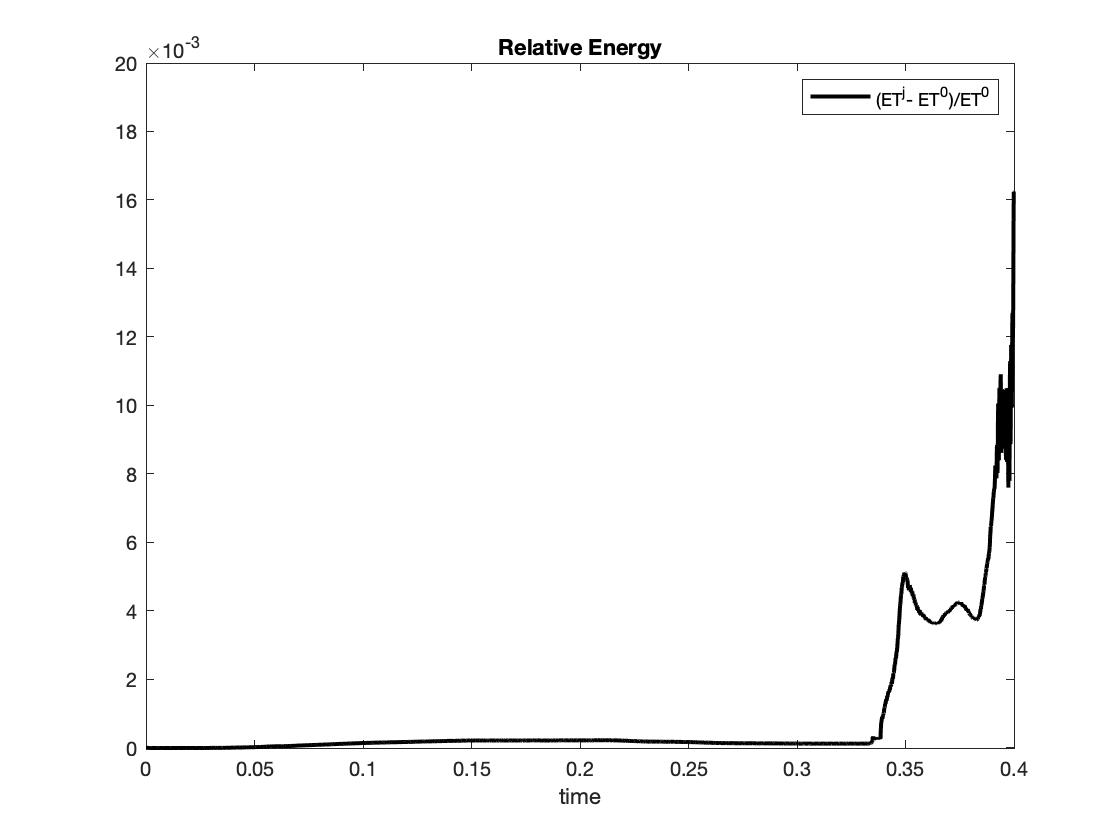} \vspace{-3pt}\\
\hspace{0.6cm} \includegraphics[width=2.2 in]{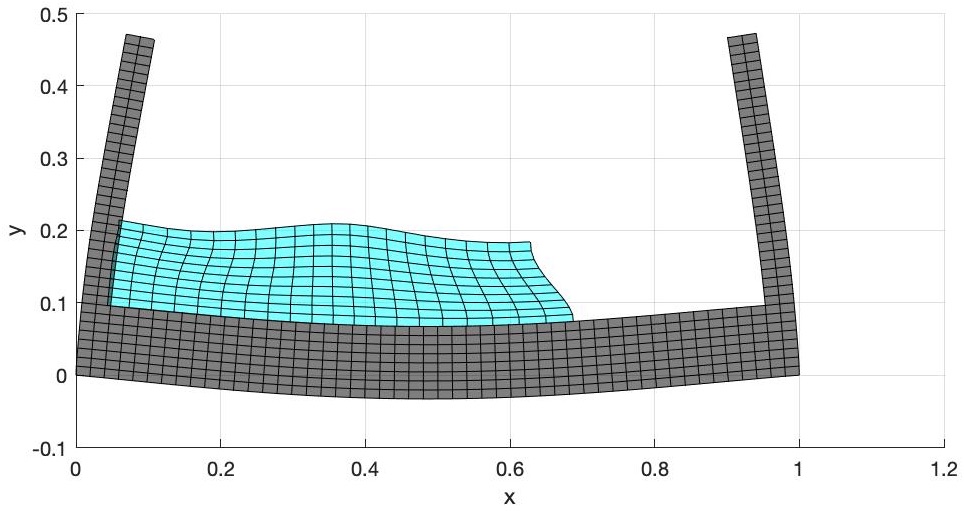}  \qquad \qquad  \quad \includegraphics[width=1.65 in]{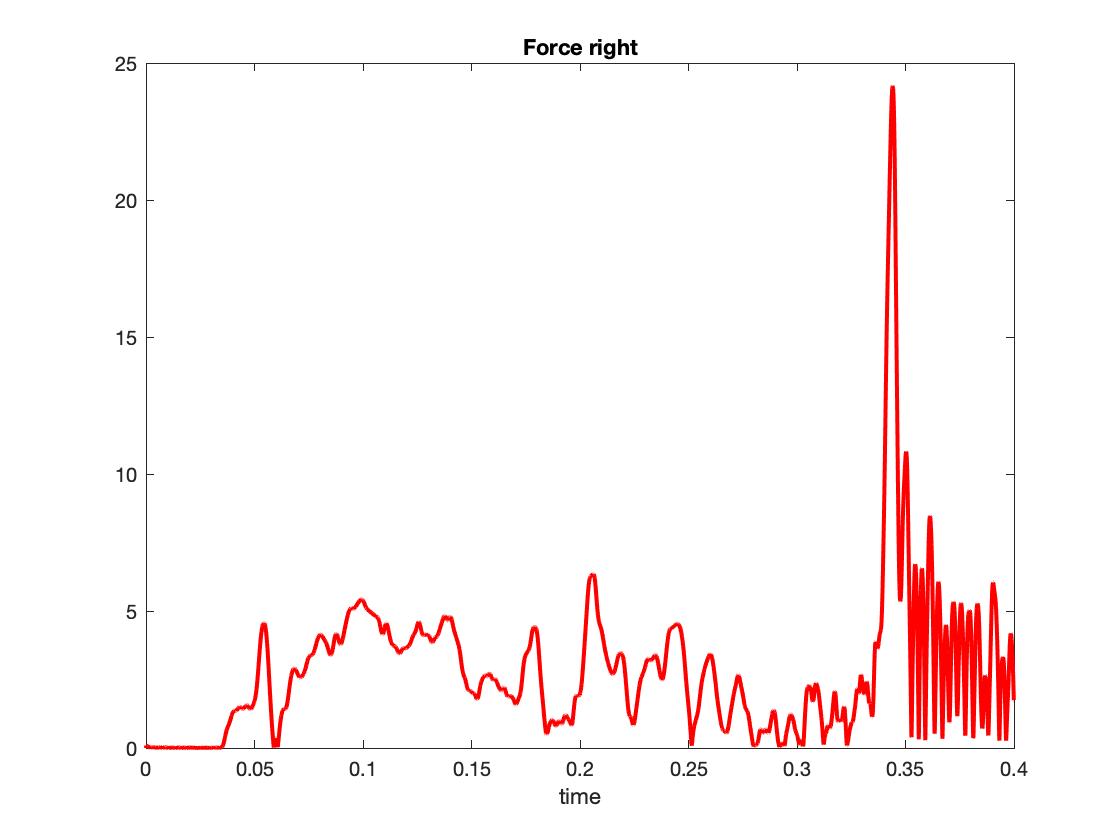} \vspace{-3pt}\\
\;  \includegraphics[width=2.6 in]{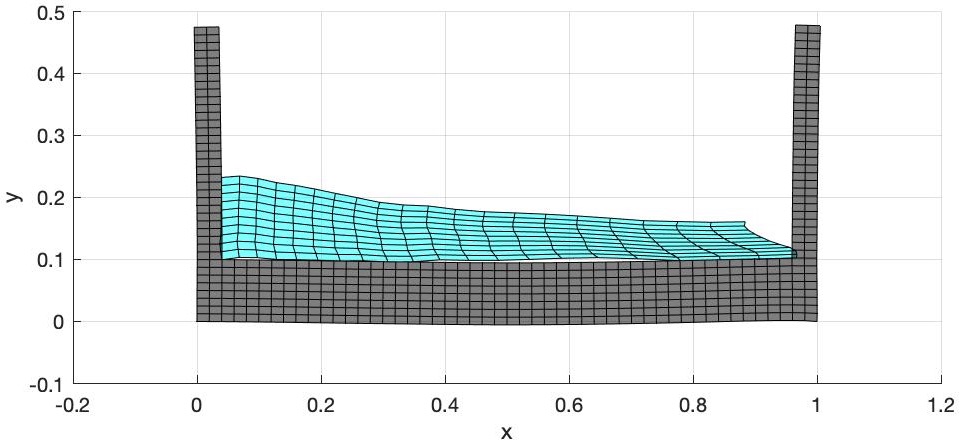}  \hspace{1.7cm} \includegraphics[width=1.65 in]{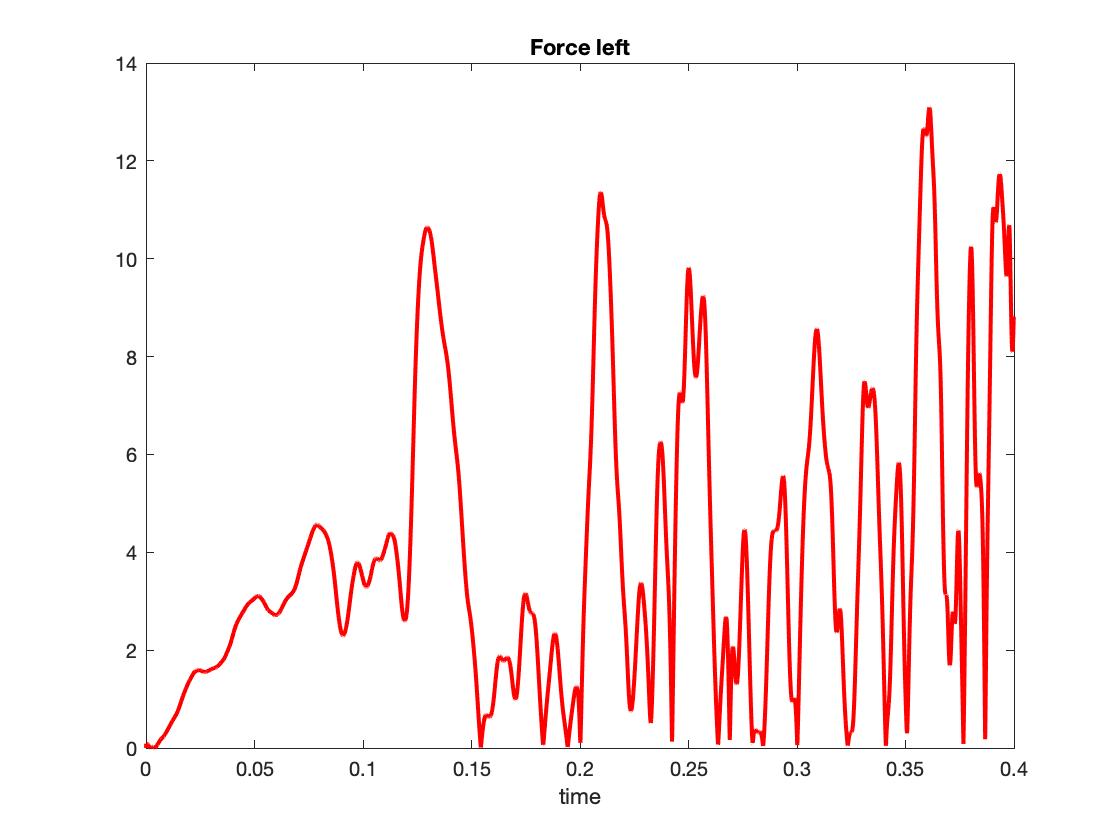} \quad \vspace{-3pt}\\
 \includegraphics[width=2.27 in]{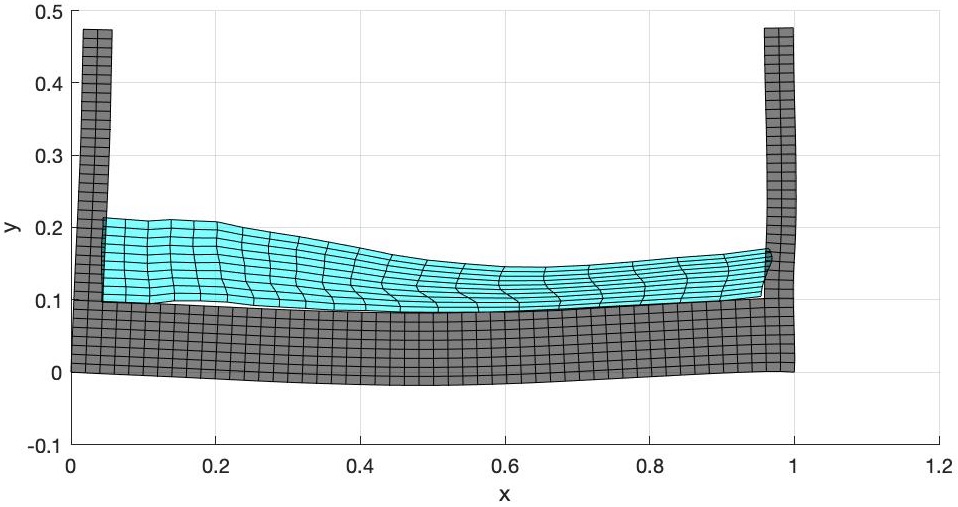} \hspace{0.55cm} \qquad \qquad \hspace{1.25 in} \vspace{-3pt}  
\caption{\footnotesize  \textit{Top to bottom on the left}: Configuration after $0.01$s, $0.1$s, $0.15$s, $0.35$s, $0.4$s. \textit{Top to bottom on the right}: Evolution of momentum maps, relative energy, and norm of the resultant of forces acting on the right and on the left of the fluid, associated respectively to pressures $P_2 ( \mbox{\mancube}_{C-1,d}^j)$ $\&$ $P_4 ( \mbox{\mancube}_{C-1,d-1}^j)$ and $P_1 ( \mbox{\mancube}_{0,d}^j)$ $\&$ $P_3 ( \mbox{\mancube}_{0,d-1}^j)$ during $0.4$s.}\label{fluid_elastic_body_2D}
\end{figure}

\begin{figure}[H] \centering 
\includegraphics[width=2.5 in]{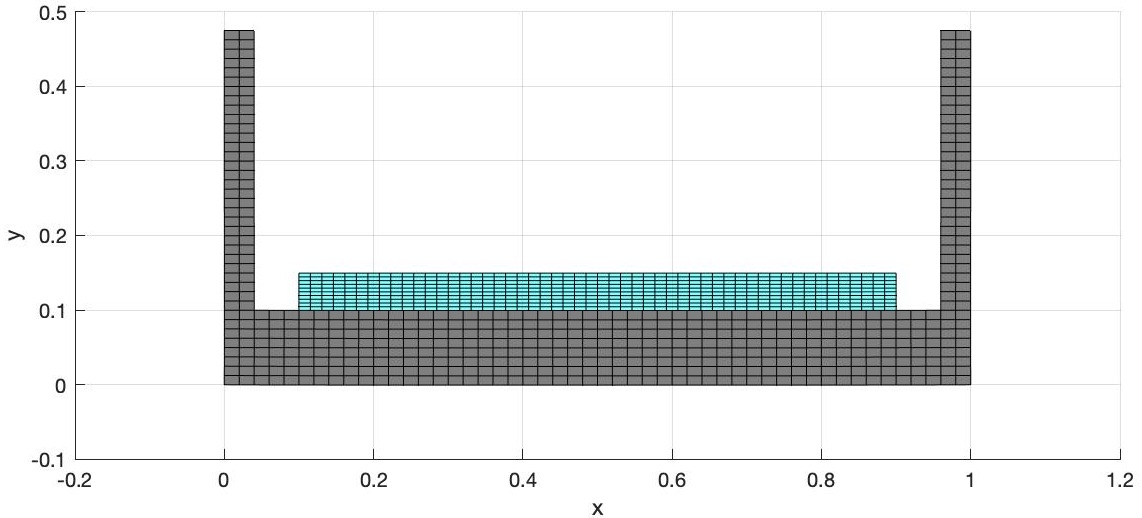} \hspace{0.15cm} \qquad  \qquad \includegraphics[width=1.65 in]{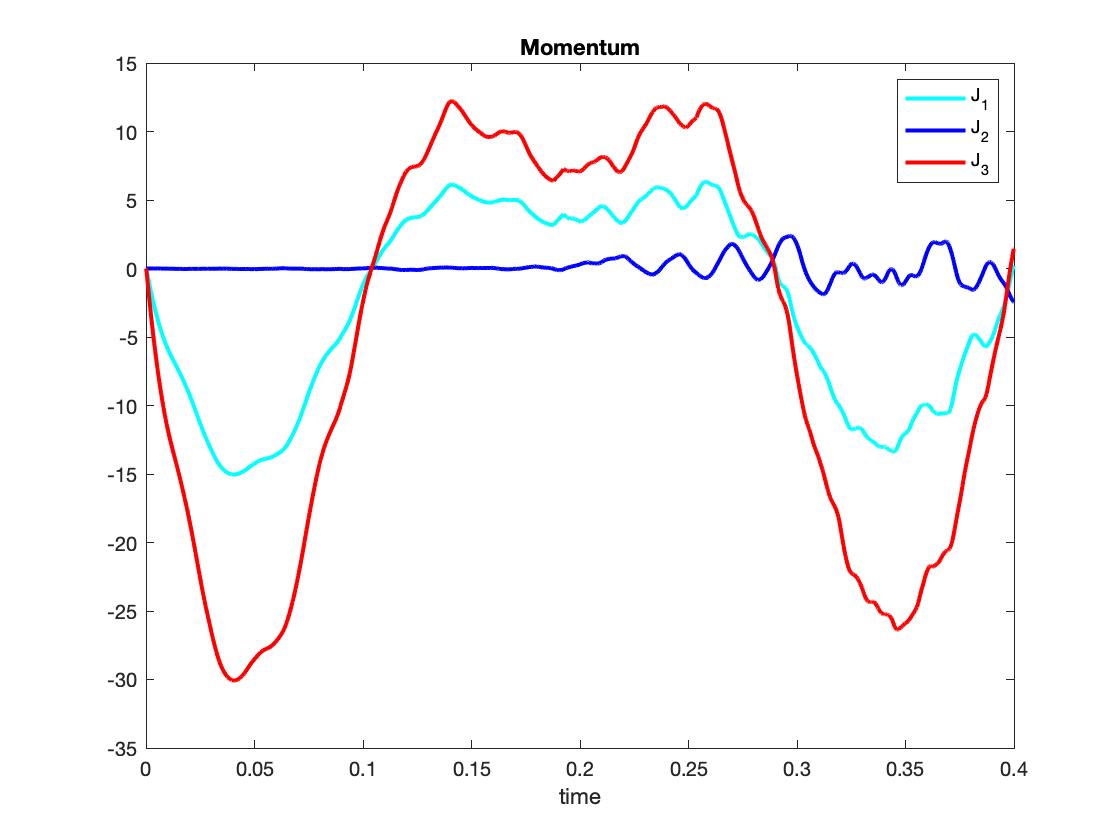} \vspace{-3pt}\\
\hspace{0.7cm}  \includegraphics[width=2.2 in]{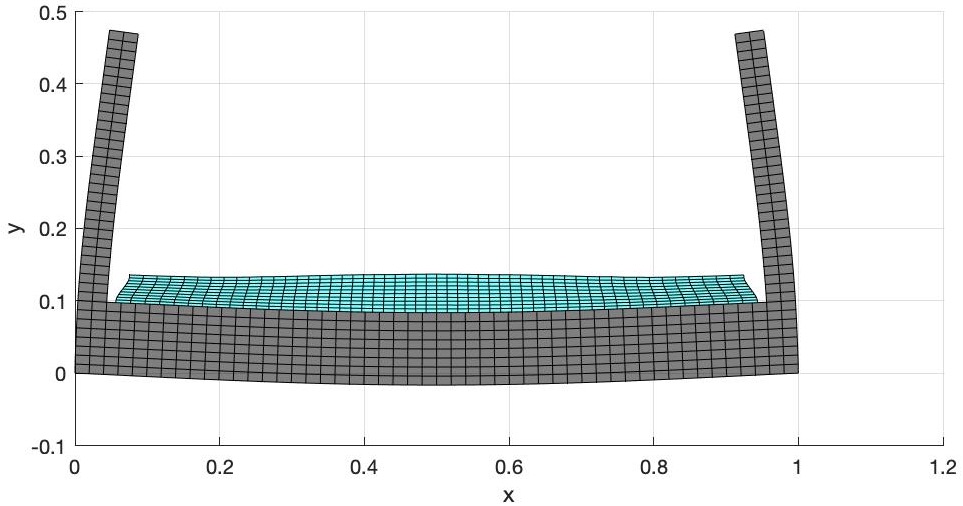}  \quad \qquad \qquad \includegraphics[width=1.65 in]{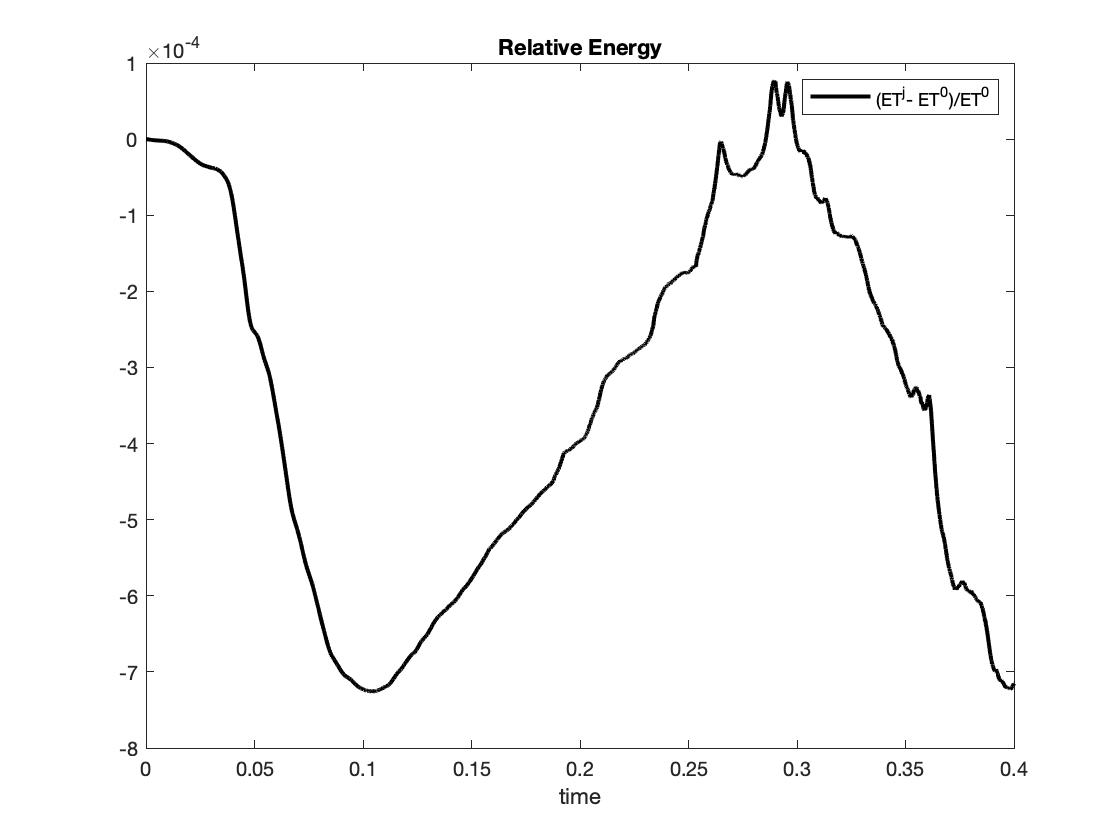} \vspace{-3pt}\\
 \includegraphics[width=2.59 in]{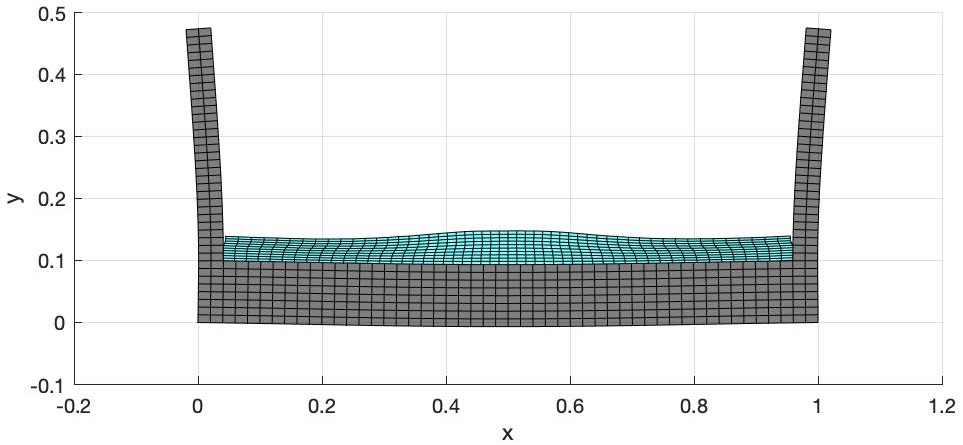}  \qquad \qquad  \quad \includegraphics[width=1.65 in]{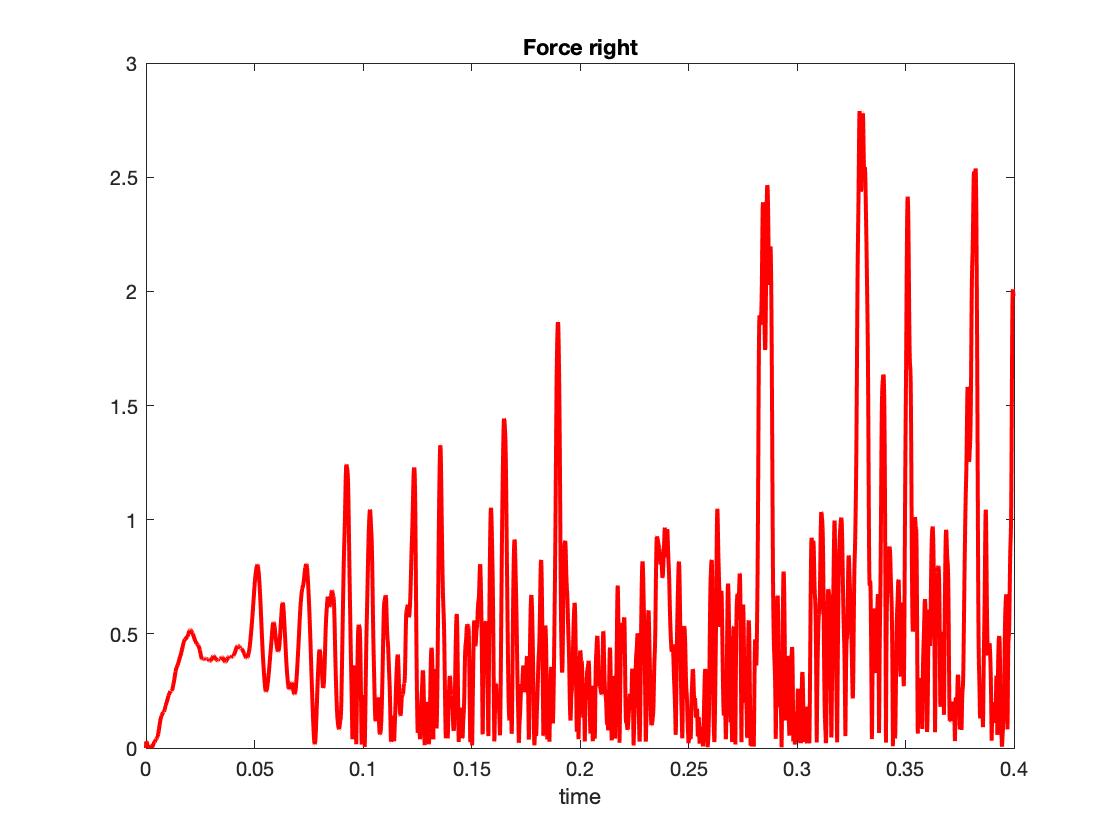} \vspace{-3pt}\\
\hspace{1.2 cm}  \includegraphics[width=2.22 in]{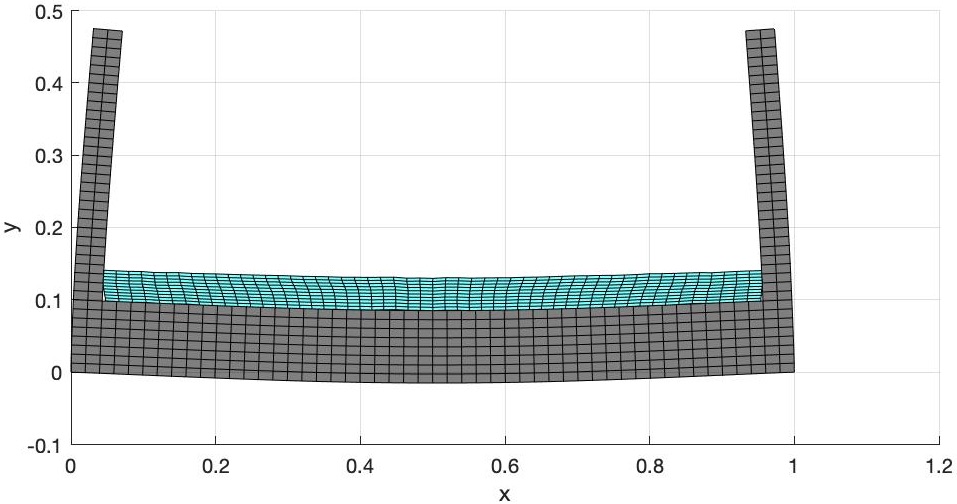}  \hspace{1.8cm} \includegraphics[width=1.65 in]{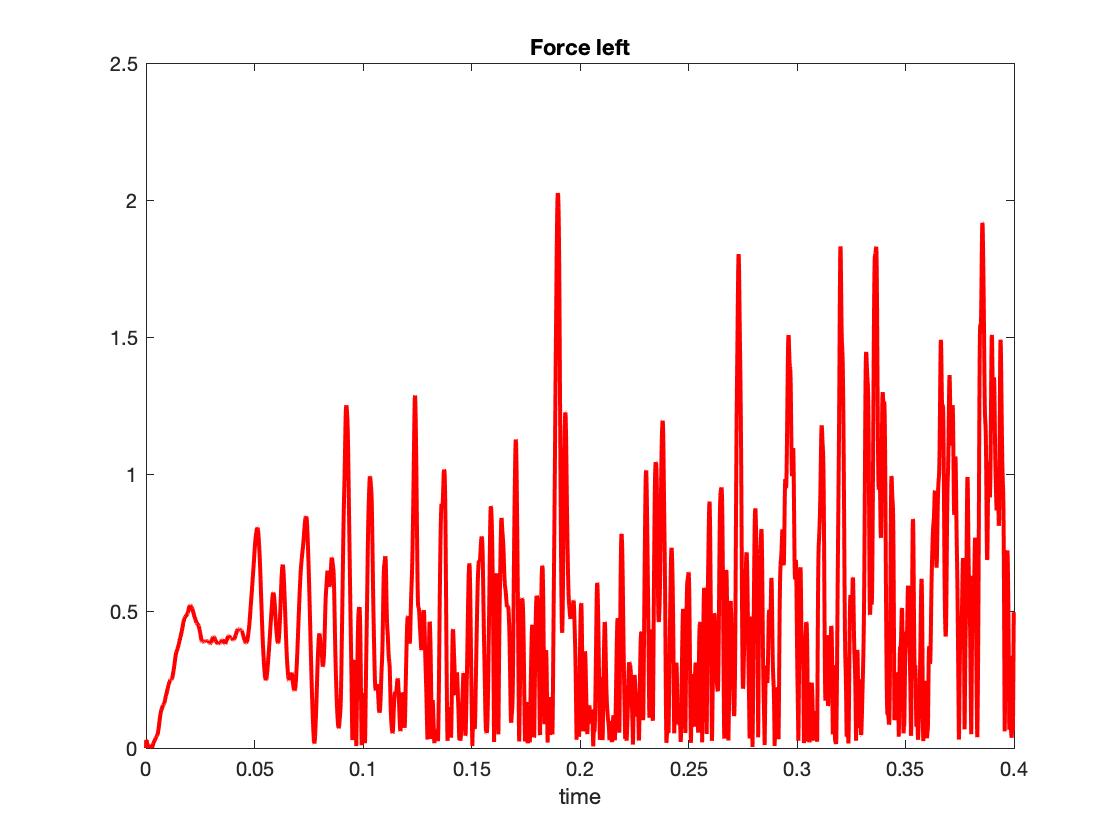} \quad \vspace{-3pt}\\
\hspace{0.142cm}\includegraphics[width=2.26 in]{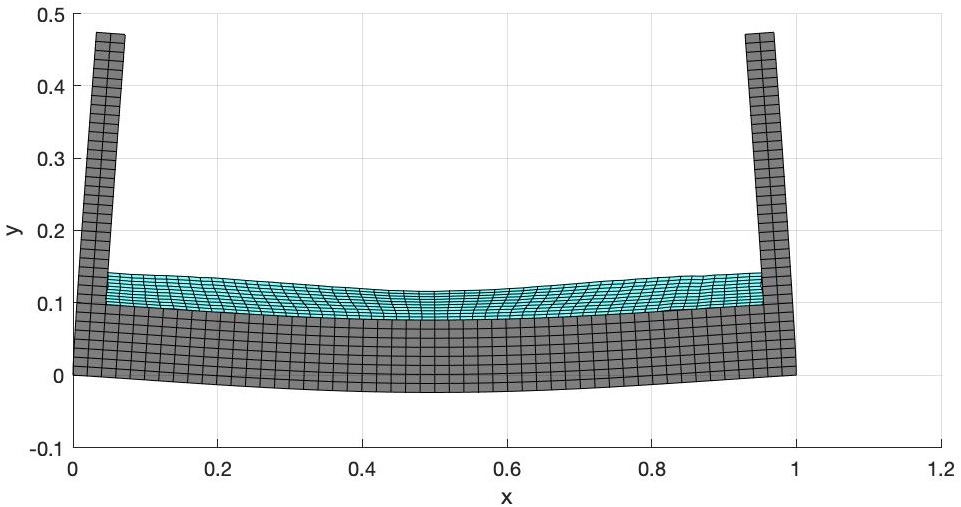} \hspace{0.55cm} \qquad \qquad \hspace{1.22 in} \vspace{-3pt}  
\caption{\footnotesize \textit{Top to bottom on the left}: Configuration after $0.01$s, $0.15$s, $0.25$s, $0.35$s, $0.4$s. \textit{Top to bottom on the right}: Evolution of momentum maps, relative energy, and norm of the resultant of forces acting on the right and on the left of the fluid, associated respectively to pressures $P_2 ( \mbox{\mancube}_{C-1,d}^j)$ $\&$ $P_4 ( \mbox{\mancube}_{C-1,d-1}^j)$ and $P_1 ( \mbox{\mancube}_{0,d}^j)$ $\&$ $P_3 ( \mbox{\mancube}_{0,d-1}^j)$ during $0.5$s.}\label{fluid_elastic_body_2D_bis}
\end{figure}

The barotropic fluid is subject to gravity and is flowing without friction inside an elastic container.
In this case, see Fig.\,\ref{impenetrability}, the impenetrability constraints take the following form along the three edges
\begin{equation} \label{2D_impenetrability_constraints}
\begin{aligned}
\Psi_{\rm im_1}(\varphi_{a,b}^j, \varphi_{a+1,b}^j, \varphi_{c,d}^j) &=  \langle  (\varphi_{c,d}^j -\varphi_{a+1,b}^j), \mathrm{R}_{\pi/2} (\varphi_{a+1,b}^j -\varphi_{a,b}^j) \rangle \geq 0
\\
\Psi_{\rm im_2}(\varphi_{a,b}^j, \varphi_{a,b+1}^j, \varphi_{c,d}^j) &= \langle (\varphi_{c,d}^j -\varphi_{a,b}^j), \mathrm{R}_{\pi/2} (\varphi_{a,b}^j -\varphi_{a,b+1}^j) \rangle \geq 0
\\
\Psi_{\rm im_3}(\varphi_{a,b}^j, \varphi_{a,b+1}^j, \varphi_{c,d}^j) &= \langle (\varphi_{c,d}^j -\varphi_{a,b+1}^j), \mathrm{R}_{\pi/2} (\varphi_{a,b+1}^j -\varphi_{a,b}^j) \rangle \geq 0.
\end{aligned}
\end{equation}
The integrator is found by computing the criticality condition for the discrete action functional \eqref{tilde_S} with incompressibility and impenetrability penalty terms. We simulate the following situation:
\begin{enumerate}
\item The barotropic fluid has the properties $\rho_0= 997 \, \mathrm{kg/m}^3$, $\gamma =6 $, $A = \tilde{A} \rho_0^{-\gamma}$ with $\tilde A = 3.041\times 10^4$ Pa, and $B = 3.0397\times 10^4$ Pa. The size of the discrete reference configuration at time $t^0$ is $0.5\mathrm{m} \times 0.15 \mathrm{m}$, with space-steps $\Delta s_1=0.025$, $\Delta s_2= 0.015$m. 

\item The hyperelastic St. Venant-Kirchhoff model has the properties $\rho_0= 945 \, \mathrm{kg/m}^3$, Poisson ratio $\nu= 0.4999$, and Young modulus $E = 4.5\times  10^6$. The incompressibility parameter is $r=10^4$ and the space-steps are $\Delta s_1=0.02$, $\Delta s_2= 0.0125$m.
\item The two corners at the bottom of the hyperelastic body are fixed.
\item The time-step is $\Delta t=10^{-4}$ and the test is carried out for $0.4s$. After this time, the impact between the fluid and the obstacle requires taking into account the physical phenomena of breaking waves.
\end{enumerate}
The results are reproduced in Fig.\,\ref{fluid_elastic_body_2D} and Fig.\,\ref{fluid_elastic_body_2D_bis}, which also show the evolution of the total momentum maps (fluid+solid), the relative total energy behavior (fluid+solid), and the norm of the resultant of forces acting on the right and on the left of the fluid. We used the notation $\mathbf{J}_d=(\mathsf{J}_1,\mathsf{J}_2, \mathsf{J}_3) \in \mathfrak{se}(2)^*$ for the angular, horizontal, and vertical component, respectively.
Note that only the component $\mathsf{J}_2$ corresponds to a symmetry of the discrete Lagrangian in presence of gravity. Its conservation is however broken, even before the fluid impacts the elastic body, see Fig.\,\ref{fluid_elastic_body_2D}, due to the fixing of the two corners at the bottom of the hyperelastic body. Its exact conservation in Fig.\,\ref{fluid_elastic_body_2D_bis} is due to the symmetric character of the initial configuration. The graph of the resultant force also illustrates the times of impact of the fluid on the right and left parts of the hyperelastic container.

The forces acting on the right at the spacetime nodes $(j,C-1,d)$, $d=1,...,D-1$,  of the fluid associated to pressures $P_2 ( \mbox{\mancube}_{C-1,d}^j)$ $\&$ $P_4 ( \mbox{\mancube}_{C-1,d-1}^j)$ are given by
\begin{equation*}{\footnotesize
\begin{aligned}
 \frac{P_2 ( \mbox{\mancube}_{C-1,d}^j)}{4}\left(\varphi_{C,d}^j - \varphi_{C,d+1}^j \right)\times \frac{\mathbf{n}_2( \mbox{\mancube}_{C-1,d}^j)}{|\mathbf{n}_2( \mbox{\mancube}_{C-1,d}^j) |} + \frac{P_4 ( \mbox{\mancube}_{C-1,d-1}^j)}{4} \left(\varphi_{C,d-1}^j - \varphi_{C,d}^j \right)\times \frac{\mathbf{n}_4( \mbox{\mancube}_{C-1,d-1}^j)}{|\mathbf{n}_4( \mbox{\mancube}_{C-1,d-1}^j) |},
\end{aligned} }
\end{equation*}
see Appendix A.2 in \cite{DeGB2021} for details. Similarly, the forces acting on the left at the spacetime nodes $(j,0,d)$, $d=1,...,D-1$, of the fluid associated to pressures $P_1 ( \mbox{\mancube}_{0,d}^j)$ $\&$ $P_3 ( \mbox{\mancube}_{0,d-1}^j)$ are 
\begin{equation*}{\footnotesize
\begin{aligned}
 \frac{P_1 ( \mbox{\mancube}_{0,d}^j)}{4} \left(\varphi_{0,d+1}^j - \varphi_{0,d}^j \right)\times \frac{\mathbf{n}_1( \mbox{\mancube}_{0,d}^j)}{|\mathbf{n}_1( \mbox{\mancube}_{0,d}^j) |}  +  \frac{P_3 ( \mbox{\mancube}_{0,d-1}^j)}{4}  \left(\varphi_{0,d}^j - \varphi_{0,d-1}^j \right)\times \frac{\mathbf{n}_3( \mbox{\mancube}_{0,d-1}^j)}{|\mathbf{n}_3( \mbox{\mancube}_{0,d-1}^j) |} .
\end{aligned} }
\end{equation*}

\subsection{Convergence tests for the hyperelastic model}

\begin{figure}[H] \centering 
\includegraphics[width=2.2 in]{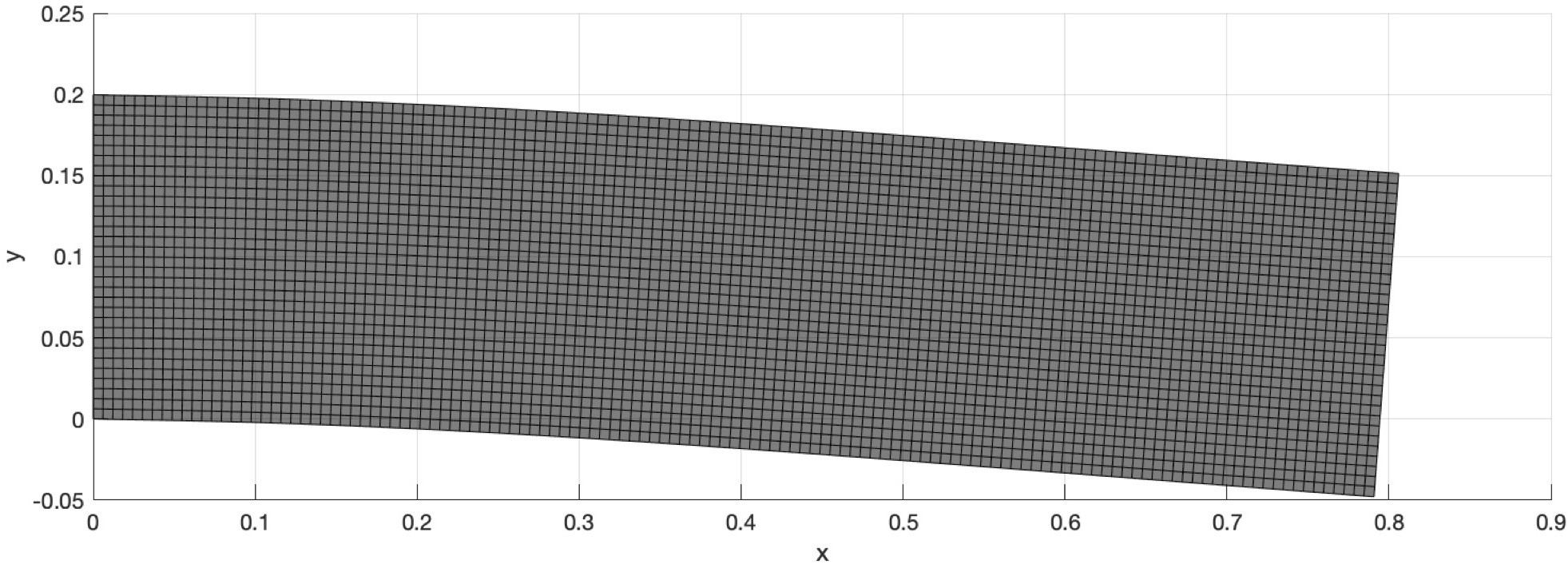} \vspace{-3pt}   \includegraphics[width=1.5 in]{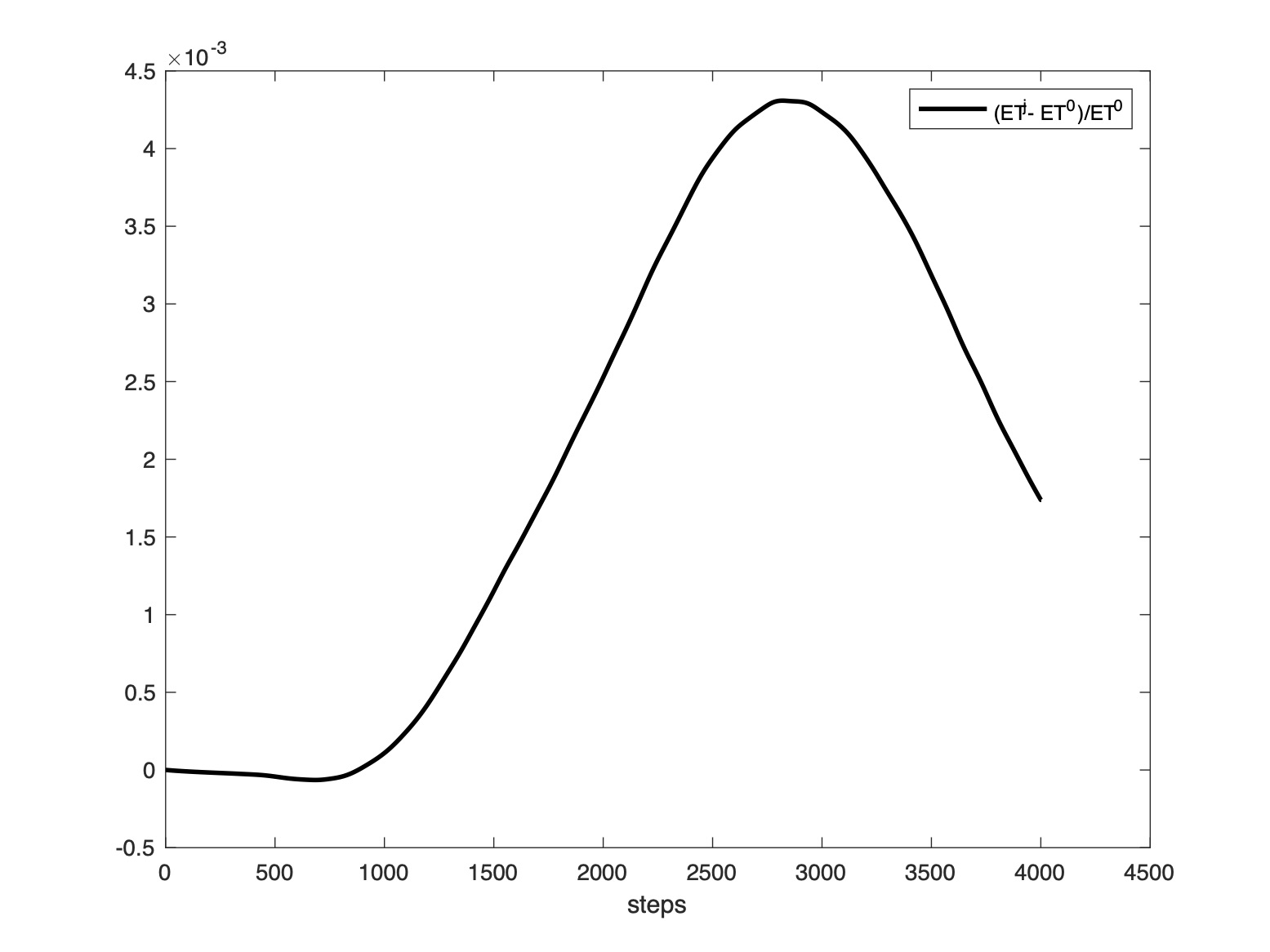} \vspace{-3pt}   \includegraphics[width=1.5 in]{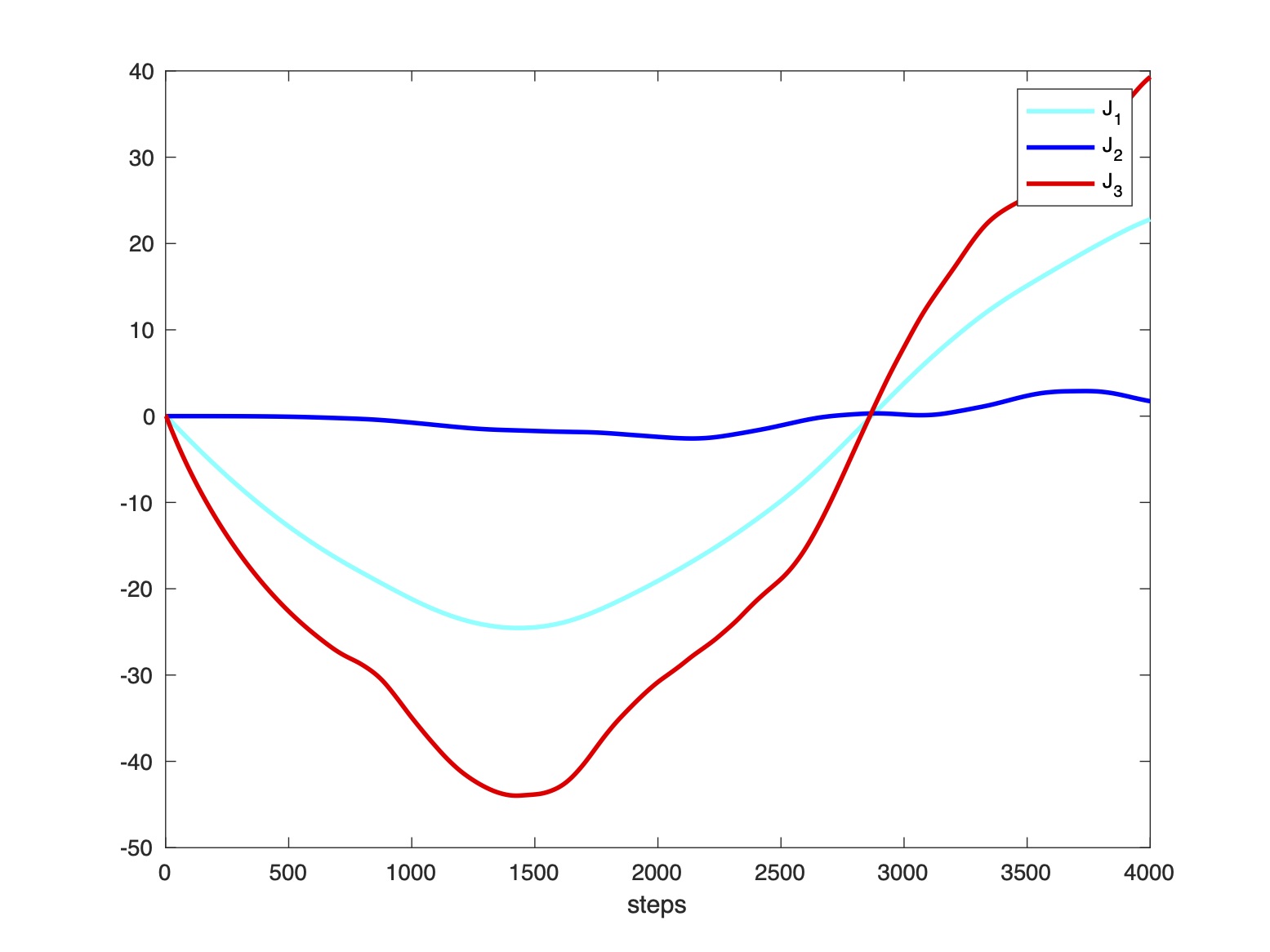} \vspace{-3pt}  
\caption{\footnotesize \textit{From left to right}: Cantilever configuration, relative energy, and momentum maps after $0.2$s, when $\Delta t=5\times 10^{-5}$ and $\Delta s_1=\Delta s_2=6.25\times 10^{-3}$. }\label{convergence_2D} 
\end{figure}

We consider a St. Venant-Kirchhoff hyperelastic cantilever with mass density $\rho_0= 945 \, \mathrm{kg/m}^3$, Poisson ratio $\nu= 0.4999$, and Young modulus $E= 2.5\times  10^6$. The model is incompressible, with coefficient $r=10^4$ in \eqref{2D_disc_penalty_function}.
The size of the rectangular reference cantilever configuration at time $t^0$ is $ L_1\times L_2 = 0.8 \mathrm{m} \times 0.2 \mathrm{m}$, see Fig.\,\ref{convergence_2D}. 

We consider the explicit integrator described in \eqref{concrete_EL}, for the elastic body only, and study the convergence with respect to $\Delta t$ and $\Delta s_i$, $i=1,2$.

\begin{enumerate}
\item
 Given a fixed mesh, with values $\Delta s_1=\Delta s_2= 0.025$m, we vary the time-steps as $\Delta t \in \{ 2 \times 10^{-4}, \, 10^{-4}, \, 5 \times 10^{-5}, \, 2.5 \times 10^{-5} \}$.
We compute the $L^2$-errors in the position $\varphi_d$ at time $t^N=0.5$s, by comparing $\varphi_d$ with an ``exact solution'' obtained with the time-step $\Delta t_{\rm ref}=6.25 \times 10^{-6}$s.  That is, for each value of $\Delta t$ we calculate
\begin{equation}\label{L_2norm_3D}
\| \varphi_d - \varphi_{\rm ref} \|_{L^2} = \left( \sum_a \sum_b  \| \varphi_{a,b}^N - \varphi_{{\rm ref};a,b}^N \|^2 \right)^{1/2}.
\end{equation}
This yields the following convergence with respect to $\Delta t$ 
\begin{figure}[H] \centering 
\begin{tabular}{| c | c | c | c | c |}
\hline
$\Delta t$ & $2 \times 10^{-4}$ & $ 10^{-4}$ & $ 5 \times 10^{-5}$ & $2.5 \times 10^{-5}$  \\
\hline
$\| \varphi_d - \varphi_{\rm ref} \|_{L^2}$ &  $1.3 \times 10^{-3}$  & $6.47\times 10^{-4}$  & $3.02\times 10^{-4}$  &  $1.29 \times 10^{-4}$ \\
\hline
$ \text{rate} $  &    & 1.007 & 1.1  & 1.23  \\
\hline
\end{tabular}
\end{figure}

\item Given a fixed time-step $\Delta t=5 \times 10^{-5}$s, we vary the space-steps as $\Delta s_1 =\Delta s_2$ $\in \{0.1, 0.05,\, 0.025, \,0.0125 \}$. The ``exact solution'' is chosen with $\Delta s_{1;\rm ref}=\Delta s_{2;\rm ref}  = 0.00625$m. We compute the $L^2$-errors in the position $\varphi_d$ at time $t^N=0.1$s. We get the following convergence with respect to $\Delta s_1=\Delta s_2$.
\begin{figure}[H] \centering 
\begin{tabular}{| c | c | c | c |c |c|c|}
\hline
$\Delta s_1 =\Delta s_2  $ & $0.1$  & $0.05$  & $0.025$    & $0.0125$\\
\hline
$\| \varphi_d - \varphi_{\rm ref} \|_{L^2}$ & $0.0481$   &  $0.0151 $ &  $0.0044$ &  $0.0012$ \\
\hline
$ \text{rate} $  &  & 1.6715  & 1.779 & 1.8745 \\
\hline
\end{tabular}
\end{figure}
The results indicate that the method is second-order accurate with respect to $\Delta s_1=\Delta s_2$.
\end{enumerate}

\section{3D Discrete fluid - elastic body interactions} \label{3D_elastic_body}

In this section we show how the setting developed previously can be extended to the three dimensional case. The main step is to find an appropriate definition of the 3D discrete deformation gradients.

\subsection{3D discrete setting and the discrete Cauchy-Green tensor} \label{3D_elast_disc}

\begin{figure}[H] \centering 
\includegraphics[width=4.9 in]{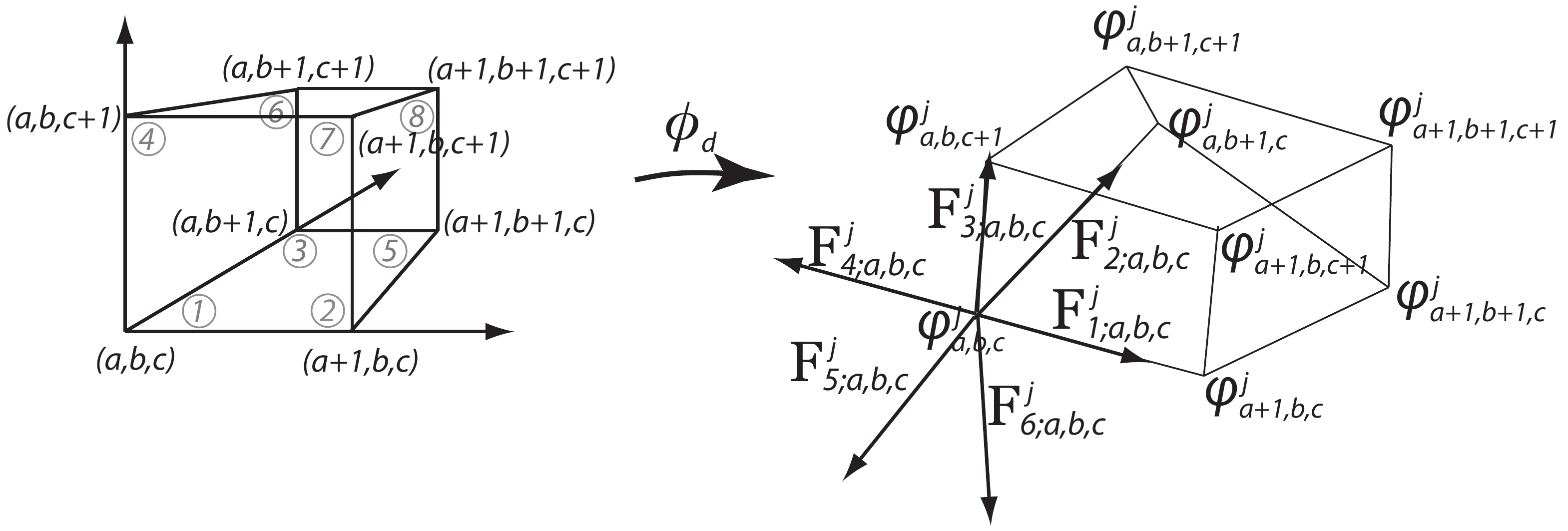} \vspace{-3pt}  
\caption{\footnotesize Discrete field $\phi_d=\varphi_d\circ  \phi _{ \mathcal{X} _d}$ evaluated on $\mbox{\mancube}_{a,b,c}^j$ at time $t^j$. }\label{phi_evalutation_3D} 
\end{figure}

\paragraph{Discrete Lagrangian setting.} With $ \mathcal{B} $ a domain in $ \mathbb{R} ^3  $ with piecewise smooth boundary and $ \mathcal{M} =\mathbb{R}  ^3  $, the discrete configuration bundle is defined as in the 2D case, see the beginning of \S\ref{2D_mult_discret} and Fig.\,\ref{D_conf_bundle}. The discrete parameter space is now $\mathcal{U}_d= \{0,...,j,...,N\} \times \mathbb{B}_d$, where we consider $ \mathbb{B}_d=\{0,...,A\} \times \{0,...,B\} \times \{0,..., C\}$ and we denote $(j,a,b,c) \in \mathcal{U} _d$. It determines a set $ \mathcal{U} _d^{\,\mbox{\mancube}}$ of parallelepipeds denoted $\mbox{\mancube}_{a,b,c}^{\,j} $,  defined by 16 pairs of indices, see Fig.\,\ref{phi_evalutation_3D} for the eight pairs of indices in $\mbox{\mancube}_{a,b,c}^j$ at time $t^j$. We assume that the discrete base-space configuration is of the form
\begin{equation}\label{3D_disc_base_space_conf}
\phi_{\mathcal{X}_d}:\mathcal{U}_d \ni (j,a,b,c) \mapsto s^j_{a,b,c} = (t^j,s^j_a,s^j_b,s^j_c) \in \mathcal{X}_d
\end{equation} 
and denote by $ \varphi ^j_{a,b,c}:= \varphi _d(s^j_{a,b,c})$ the value of the discrete field  at $ s^j_{a,b,c}$. The discrete first jet bundle and discrete first jet extensions are defined as in the 2D case earlier. In particular, we write the first jet extension as follows{ \fontsize{9pt}{13pt}\selectfont
\begin{equation}\label{discrete_jet_extension_3D} 
\begin{aligned} 
&j^1 \varphi _d( \mbox{\mancube}_{a,b,c}^j)= \big( \varphi _{a,b,c}^j, \varphi _{a,b,c}^{j+1},  \varphi _{a+1,b,c}^j, \varphi _{a+1,b,c}^{j+1},  \varphi _{a,b+1,c}^j ,  \varphi _{a,b+1,c}^{j+1}, \varphi _{a,b,c+1}^j, \varphi _{a,b,c+1}^{j+1},\\
&\quad \varphi _{a+1,b+1,c}^j , \varphi _{a+1,b+1,c}^{j+1}, \varphi _{a,b+1,c+1}^j, \varphi _{a,b+1,c+1}^{j+1} ,\varphi _{a+1,b,c+1}^j, \varphi _{a+1,b,c+1}^{j+1},\varphi _{a+1,b+1,c+1}^j, \varphi _{a+1,b+1,c+1}^{j+1}\big),
\end{aligned} 
\end{equation} }

\paragraph{Discrete Cauchy-Green deformation tensor in 3D.} Given a discrete base-space configuration $ \phi _{ \mathcal{X} _d}$ and a discrete field $ \varphi _d$, the following six vectors $\mathbf{F}_{\ell;a,b,c}^j \in \mathbb{R} ^3$, $\ell =1,...,6$ are defined at each node $(j,a,b,c)\in \mathcal{U}_d$, see Fig.\,\ref{phi_evalutation_3D} on the right:
\begin{equation*}
\begin{aligned}
\mathbf{F}_{1;a,b,c}^j & = \frac{\varphi_{a+1,b,c}^j - \varphi_{a,b,c}^j}{ |s^j_{a+1} - s^j_{a}| } , \quad  \mathbf{F}_{2;a,b,c}^j = \frac{\varphi_{a,b+1,c}^j- \varphi_{a,b,c}^j}{|s^j_{b+1} - s^j_{b}| }, \quad \mathbf{F}_{3;a,b,c}^j = \frac{\varphi_{a,b,c+1}^j- \varphi_{a,b,c}^j}{|s^j_{c+1} - s^j_{c}| },\\
\mathbf{F}_{4;a,b,c}^j &=  \frac{\varphi_{a-1,b,c}^j - \varphi_{a,b,c}^j}{|s^j_{a} - s^j_{a-1}| }, \quad\mathbf{F}_{5;a,b,c}^j = \frac{\varphi_{a,b-1,c}^j- \varphi_{a,b,c}^j}{ |s^j_{b} - s^j_{b-1}|}, \quad \mathbf{F}_{6;a,b,c}^j= \frac{\varphi_{a,b,c-1}^j - \varphi_{a,b,c}^j}{|s^j_{c} - s^j_{c-1}| }.
\end{aligned}
\end{equation*}

\medskip

\begin{definition}  \label{3D_gradient_definition}
The \textbf{discrete deformation gradients} of a discrete field $ \varphi _d$ at the parallelepiped $\mbox{\mancube}_{a,b}^j$ are the eight $3 \times 3$ matrices $\mathbf{F}_\ell(\mbox{\mancube}_{a,b,c}^j)$, $\ell =1,...,8$, defined by
\begin{equation}\label{3D_gradient_def}
\begin{aligned}
& \mathbf{F}_1(\mbox{\mancube}_{a,b,c}^j) = \left[\mathbf{F}_{1;a,b,c}^j \; \; \mathbf{F}_{2;a,b,c}^j \; \; \mathbf{F}_{3;a,b,c}^j \right],\\ 
&\mathbf{F}_2(\mbox{\mancube}_{a,b,c}^j) = \left[\mathbf{F}_{2;a+1,b,c}^j \;\; \mathbf{F}_{4;a+1,b,c}^j  \; \;  \mathbf{F}_{3;a+1,b,c}^j \right] ,\\
&  \mathbf{F}_3(\mbox{\mancube}_{a,b,c}^j) = \left[\mathbf{F}_{5;a,b+1,c}^j  \; \; \mathbf{F}_{1;a,b+1,c}^j \; \; \mathbf{F}_{3;a,b+1,c}^j \right],\\
&   \mathbf{F}_4(\mbox{\mancube}_{a,b,c}^j) = \left[ \mathbf{F}_{2;a,b,c+1}^j  \; \; \mathbf{F}_{1;a,b,c+1}^j \; \; \mathbf{F}_{6;a,b,c+1}^j \right],\\
&  \mathbf{F}_5(\mbox{\mancube}_{a,b,c}^j) = \left[\mathbf{F}_{4;a+1,b+1,c}^j  \; \; \mathbf{F}_{5;a+1,b+1,c}^j  \; \; \mathbf{F}_{3;a+1,b+1,c}^j \right],\\
& \mathbf{F}_6(\mbox{\mancube}_{a,b,c}^j) = \left[\mathbf{F}_{1;a,b+1,c+1}^j  \; \; \mathbf{F}_{5;a,b+1,c+1}^j \; \; \mathbf{F}_{6;a,b+1,c+1}^j \right],\\
& \mathbf{F}_7(\mbox{\mancube}_{a,b,c}^j) = \left[\mathbf{F}_{4;a+1,b,c+1}^j \; \; \mathbf{F}_{2;a+1,b,c+1}^j  \; \; \mathbf{F}_{6;a+1,b,c+1}^j \right], \\
& \mathbf{F}_8(\mbox{\mancube}_{a,b,c}^j) = \left[\mathbf{F}_{5;a+1,b+1,c+1}^j  \; \; \mathbf{F}_{4;a+1,b+1,c+1}^j \; \; \mathbf{F}_{6;a+1,b,c+1}^j \right].
\end{aligned}
\end{equation}
\end{definition}

\medskip

We note that the discrete deformation gradients are defined at each of the eight nodes of $\mbox{\mancube}_{a,b,c}^j$ that are associated to time $t^j$. We assume that the determinant of these matrices is always positive. The ordering $\ell=1$ to $\ell=8$ is respectively associated to the nodes $(a,b,c)$, $(a+1,b,c)$, $(a,b+1,c)$, $(a,b,c+1)$, $(a+1,b+1,c)$, $(a,b+1,c+1)$, $(a+1,b,c+1)$, $(a+1,b+1,c+1)$ see Fig.\,\ref{phi_evalutation_3D} on the left. Note that with this concept of discrete deformation gradient, one can consider the right and left polar decompositions $ \mathbf{F} = \mathbf{R} \mathbf{U}= \mathbf{V} \mathbf{R} $ exactly as in the continuous case.

\begin{definition}\label{CG_definition}
The \textbf{3D discrete Cauchy-Green deformation tensors} of a discrete field $ \varphi _d$ at the parallelepiped $\,\mbox{\mancube}_{a,b,c}^j$ are the eight $3 \times 3$ symmetric matrices
\begin{equation}\label{def_CG_3D} 
\mathbf{C} _\ell( \mbox{\mancube}_{a,b,c}^j) = \mathbf{F} _\ell( \mbox{\mancube}_{a,b,c}^j) ^\mathsf{T}\mathbf{F} _\ell( \mbox{\mancube}_{a,b,c}^j), \quad \ell=1,...,8.
\end{equation} 
\end{definition}
\medskip 

As in Proposition \ref{2D_DCG}, we have following result.

\begin{proposition}
Given a discrete field $ \varphi _d: \mathcal{X}_d \rightarrow \mathcal{M} $, the 3D discrete Cauchy-Green deformation tensor associated to $\mbox{\mancube}_{a,b,c}^j$ at the node $(a,b,c)$ at time $t^j$, for $\ell=1$, is 
 \begin{equation}\label{3D_Cauch_Green}
 \mathbf{C}_1(\mbox{\mancube}_{a,b,c}^j) =
 \begin{bmatrix}
\langle\mathbf{F}_{1;a,b,c}^j , \mathbf{F}_{1;a,b,c}^j \rangle & \; \langle \mathbf{F}_{1;a,b,c}^j , \mathbf{F}_{2;a,b,c}^j\rangle & \; \langle \mathbf{F}_{1;a,b,c}^j , \mathbf{F}_{3;a,b,c}^j\rangle \\[6pt]
\langle \mathbf{F}_{2;a,b,c}^j , \mathbf{F}_{1;a,b,c}^j\rangle  & \; \langle\mathbf{F}_{2;a,b,c}^j,  \mathbf{F}_{2;a,b,c}^j\rangle & \; \langle\mathbf{F}_{2;a,b,c}^j,  \mathbf{F}_{3;a,b,c}^j\rangle \\[6pt]
\langle \mathbf{F}_{3;a,b,c}^j , \mathbf{F}_{1;a,b,c}^j\rangle  & \; \langle\mathbf{F}_{3;a,b,c}^j,  \mathbf{F}_{2;a,b,c}^j\rangle & \; \langle\mathbf{F}_{3;a,b,c}^j,  \mathbf{F}_{3;a,b,c}^j\rangle   \end{bmatrix}.
 \end{equation}
See Appendix \ref{3D_Cauchy_Green} for the discrete Cauchy-Green deformation tensor at the nodes $(a+1,b,c)$, $(a,b+1,c)$, $(a,b,c+1)$, $(a+1,b+1,c)$, $(a,b+1,c+1)$, $(a+1,b,c+1)$, $(a+1,b+1,c+1)$.
\end{proposition}

\paragraph{Discrete Jacobian.}  The 3D discrete Jacobian of $ \varphi _d$ associated to $\mbox{\mancube}_{a,b,c}^{\,j}$ at the node $(j,a,b,c)$ is defined by
\begin{equation} \label{3D_Jacobian}
 J_1(\mbox{\mancube}_{a,b,c}^j) =  (\mathbf{F}_{1;a,b,c}^j \times \mathbf{F}_{2;a,b,c}^j)\cdot \mathbf{F}_{3;a,b,c}^j.
\end{equation}
We refer to Appendix \ref{3D_Jacobian's} for the expression of the 3D discrete Jacobians at the other nodes. The relations with the 3D discrete deformation gradients and discrete Cauchy-Green tensors hold exactly as in the 2D case in Proposition \ref{discrete_2D_Jac}.

\subsection{Discrete frame indifferent and isotropic nonlinear hyperelastic models}\label{isotropic_hyperelastic}

\paragraph{Invariants of the discrete Cauchy-Green deformation tensor.} Our definition of the discrete Cauchy-Green deformation tensor as a symmetric matrix allows the definition of its invariants exactly as in the continuous case. Note that in the discrete case each node has several associated discrete Cauchy-Green deformation tensors, namely, one for each parallelepiped sharing this node. This follows directly from the discrete field theoretic Lagrangian setting that we follow. We have the following expressions.

\begin{definition}\label{def_discrete_invariants_C}
The invariants of the symmetric matrix $\mathbf{C}_\ell (\mbox{\mancube}_{a,b,c}^j)$, $\ell=1,...,8$, are defined by
\begin{equation}\label{discrete_invariants_C_d}
\begin{aligned}
I_1\big(\mathbf{C}_\ell (\mbox{\mancube}_{a,b,c}^j)\big) &= \mathrm{Tr}\big(\mathbf{C}_\ell (\mbox{\mancube}_{a,b,c}^j)\big)\\
I_2\big(\mathbf{C}_\ell (\mbox{\mancube}_{a,b,c}^j)\big) &=\frac{1}{2} \left[I_1\big(\mathbf{C}_\ell (\mbox{\mancube}_{a,b,c}^j)\big)^2 - \mathrm{Tr}\big(\mathbf{C}_\ell (\mbox{\mancube}_{a,b,c}^j)^2\big) \right]\\
&= \mathrm{det}\big(\mathbf{C}_\ell (\mbox{\mancube}_{a,b,c}^j)\big)\mathrm{Tr}\big(\mathbf{C}_\ell (\mbox{\mancube}_{a,b,c}^j )^{-1}\big)\\
I_3\big(\mathbf{C}_\ell (\mbox{\mancube}_{a,b,c}^j)\big)&= \mathrm{det}\big(\mathbf{C}_\ell (\mbox{\mancube}_{a,b,c}^j)\big).
\end{aligned}
\end{equation}
\end{definition}

For each fixed $\mbox{\mancube}_{a,b,c}^j$, these invariants depend only of the first jet extension $j^ 1 \varphi _d( \mbox{\mancube}_{a,b,c}^j)$ of the discrete field at the given $\mbox{\mancube}_{a,b,c}^j$.
Recall that these invariants are related to the coefficients in the characteristic polynomial $P(\nu)$ of the $3\times 3$ symmetric matrix $\mathbf{C}_\ell \big(\mbox{\mancube}_{a,b,c}^j\big)$ as follows
\[
P(\nu)=\nu^3 - I_1\big(\mathbf{C}_\ell (\mbox{\mancube}_{a,b,c}^j)\big) \nu^2 + I_2\big(\mathbf{C}_\ell (\mbox{\mancube}_{a,b,c}^j)\big) \nu - I_3\big(\mathbf{C}_\ell (\mbox{\mancube}_{a,b,c}^j)\big).
\]
In particular, the eigenvalues $\nu _1, \nu _2, \nu _3$ are defined at each node with respect to each of the eight neighbouring parallelepipeds. For each $\mbox{\mancube}_{a,b,c}^j \in \mathcal{X} _d^{\, \mbox{\mancube}}$ and $\ell=1,...,8$, we have
\begin{align*}
&I_1\big(\mathbf{C}_\ell (\mbox{\mancube}_{a,b,c}^j)\big) = \nu_1+\nu_2 + \nu_3, \qquad I_2\big(\mathbf{C}_\ell (\mbox{\mancube}_{a,b,c}^j)\big) = \nu_1\nu_2 + \nu_1\nu_3 + \nu_2\nu_3,\\
&I_3\big(\mathbf{C}_\ell (\mbox{\mancube}_{a,b,c}^j)\big) = \nu_1\nu_2\nu_3,
\end{align*}
and $\nu_1=\lambda_1^2$, $\nu_2=\lambda_2^2$, $\nu_3=\lambda_3^2$, with $\lambda_1$, $\lambda_2$, $\lambda_3$ the discrete principal stretches. At each node, there are thus several classes of discrete principal stretches, namely, one for each parallelepiped sharing this nodes. We can also recall the following standard result, which holds here in the discrete case and can be regarded as a discrete analogue of the characterization of the stored energy function for frame indifferent isotropic hyperelastic materials.

\begin{proposition}\label{isotropic_prop}
The following are equivalent:
\begin{enumerate}              [label=\alph*),itemsep=0 pt]
\item[\rm (i)] A scalar function $f$ of $\mathbf{C}_\ell (\mbox{\mancube}_{a,b,c}^j)$ is invariant under orthogonal transformations.
\item[\rm (ii)] $f$ is a function of the invariants of $\mathbf{C}_\ell (\mbox{\mancube}_{a,b,c}^j)$.
\item[\rm (iii)] $f$ is a symmetric function of the principal stretches.
\end{enumerate}
\end{proposition}

\paragraph{3D Discrete stored energy function.} Given the stored energy function $W(I_1,I_2,I_3)$ of a frame indifferent isotropic hyperelastic material, the associated discrete stored energy function is defined as $W_d:J^1 \mathcal{Y} _d \rightarrow \mathbb{R}$ with \begin{equation} \label{energy_iso_hyperelast}
W_d\big(j^1 \varphi _d(\mbox{\mancube})\big) =  \frac{1}{8}  \sum_{\ell=1}^{8}  W\big(I_1(\mathbf{C}_\ell (\mbox{\mancube})),I_2(\mathbf{C}_\ell (\mbox{\mancube})), I_3(\mathbf{C}_\ell (\mbox{\mancube}))\big).
\end{equation}
For each $\mbox{\mancube} \in \mathcal{X} _d^{\,\mbox{\mancube}}$ the value of $W_d$ is found as the averaged value of the stored energy function evaluated on the discrete Cauchy-Green tensor at the nodes of $\mbox{\mancube}$ at time $t^j$.

\paragraph{3D Discrete second Piola-Kirchhoff stress.} It is defined at each node of each $\mbox{\mancube} \in \mathcal{X} _d^{\,\mbox{\mancube}}$ as
\begin{equation}\label{3D_SPiolaK}
\mathbf{S}_\ell( \mbox{\mancube}) = 2 \rho_0 \frac{\partial W}{\partial \mathbf{C}}\big( \mathbf{C}_\ell(\mbox{\mancube})\big)= \mathbf{S}\big( \mathbf{C}_\ell(\mbox{\mancube})\big),
\end{equation}
$\ell=1,...,8$. In particular, for an isotropic material, it takes the form
\begin{equation} \label{sec_Piola_Kir2}
\mathbf{S}_\ell( \mbox{\mancube})=2 \rho_0 \left[ \frac{\partial W}{\partial I_1} \mathbf{G}^{-1}  + \left(\frac{\partial W}{\partial I_2}I_2 + \frac{\partial W}{\partial I_3}I_3 \right) \big(\mathbf{C}_\ell( \mbox{\mancube})\big)^{-1} - \frac{\partial W}{\partial I_2}I_3 \big(\mathbf{C}_\ell (\mbox{\mancube})\big)^{-2} \right].
\end{equation}

\medskip

\begin{example}[Nonlinear Mooney-Rivlin incompressible model]{\rm For this model, see \eqref{Mooney_Rivlin}, the discrete stored energy becomes
\begin{equation} \label{3D_disc_Mooney_Rivlin}
W_d\big(j^1 \varphi _d(\mbox{\mancube}_{a,b,c}^j)\big) = \frac{1}{8}  \sum_{\ell=1}^{8}  \Big[ C_1 \big( I_1\big(\mathbf{C}_\ell (\mbox{\mancube}_{a,b,c}^j)\big) -3 \big) + C_2 \big( I_2\big(\mathbf{C}_\ell (\mbox{\mancube}_{a,b,c}^j)\big) -3 \big) \Big]
\end{equation}
and the 3D discrete second Piola-Kirchhoff stress is
\[
\mathbf{S}_\ell( \mbox{\mancube}) = 2 \rho  _0 \Big(C_1 \mathbf{G} ^{-1} + C_2 I_2\big(\mathbf{C}_\ell( \mbox{\mancube}) \big) \mathbf{C}_\ell ( \mbox{\mancube})^{-1} - C_2 I_3\big(\mathbf{C}_\ell ( \mbox{\mancube})\big) \mathbf{C} _\ell( \mbox{\mancube})^{-2}\Big).
\]
Incompressibility is imposed via a Lagrange multiplier exactly as in the 2D case earlier.}
\end{example}

\subsection{3D Barotropic fluid - elastic body interactions}

\paragraph{Discrete Lagrangians for the fluid and the elastic body.} The discrete Lagrangian $ \mathcal{L} ^\mathsf{e}_d$ of the elastic body is defined exactly as in the 2D case in \eqref{Discrete_Lagrangian_2D_wave} with the obvious modifications, namely we take the discrete kinetic energy 
\begin{equation}\label{3D_kinetic_energy}
K_d\big(j^1 \varphi _d(\mbox{\mancube}_{a,b,c}^j)\big) : = \frac{1}{8} \sum_{\alpha=a}^{a+1} \sum_{\beta=b}^{b+1} \sum_{\gamma=c}^{c+1} \frac{1}{2} \left| \mathbf{v}_{\alpha,\beta, \gamma}^j\right|^2 ,
\end{equation}
where $\mathbf{v}_{\alpha,\beta,\gamma}^j= (\varphi_{a,b,c}^{j+1}- \varphi_{a,b,c}^j)/\Delta t$, we take the discrete stored energy function
\begin{equation}\label{3D_stor_ener}
W_d^\mathsf{e}\big(j^1 \varphi _d(\mbox{\mancube}_{a,b,c}^j)\big) := \frac{1}{8}  \sum_{\ell=1}^{8}  W^\mathsf{e}\big(\mathbf{C}_\ell(\mbox{\mancube}_{a,b,c}^j)\big)
\end{equation}
and the gravitational potential energy density
\begin{equation}\label{3D_grav_potential}
\Pi_d(\mbox{\mancube}_{a,b,c}^j) := \frac{1}{8}  \sum_{\alpha=a}^{a+1} \sum_{\beta=b}^{b+1} \sum_{\gamma=c}^{c+1} \Pi ( \varphi_{\alpha,\beta, \gamma}^j).
\end{equation}
We proceed similarly for the discrete Lagrangian $ \mathcal{L} _d^\mathsf{f}$ of the fluid, see \eqref{Discrete_Lagrangian_2D_fluid}.

\paragraph{Discrete action and discrete Euler-Lagrange equations.} As earlier, we assume $ \phi _ { \mathcal{X} _d}(j,a,b,c)= (j \Delta  t, a \Delta s_1, b \Delta s_2, c \Delta s_3)$ and the mass of each 3D cell in $ \Phi _{ \mathcal{X} _d}( \mathbb{B}_d)$ is $M^\mathsf{e} = \rho  _0^\mathsf{e} \Delta s_1 \Delta s_2 \Delta s_3$. The discrete action functional reads
 \begin{equation}\label{3D_discrete_action_sum}
 S_d(\varphi_d) = \sum_{j=0}^{N-1} \left( \sum_{a=0}^{A-1} \sum_{b=0}^{B-1} \sum_{c=0}^{C-1}\mathcal{L}^\mathsf{e}_{d}\big(j^1 \varphi _d(\mbox{\mancube}_{a,b,c}^j)\big) + \sum_{d=0}^{D-1} \sum_{e=0}^{E-1} \sum_{f=0}^{F-1}\mathcal{L}^\mathsf{f}_{d} \big(j^1 \varphi _d(\mbox{\mancube}_{d,e,f}^j)\big)\right).
\end{equation}

Let us denote by $D_k$, $k=1,...,16$, the partial derivative of the discrete Lagrangians with respect to the $k$-th component of $j^1 \varphi _d (\mbox{\mancube} )$, in the order given in \eqref{discrete_jet_extension_3D}. From the discrete Hamilton principle $ \delta S_d( \varphi _d)=0$, we get the discrete Euler-Lagrange equations in the general form
\begin{equation}\label{discrete_Euler_Lagrange} { \fontsize{10pt}{13pt}\selectfont
\begin{aligned}
&D _1 \mathcal{L} ^{j}_{a,b,c} + D _2 \mathcal{L} ^{j-1}_{a,b,c} + D _3 \mathcal{L} ^{j}_{a-1,b,c} + D _4\mathcal{L} ^{j-1}_{a-1,b,c} +D _5 \mathcal{L} ^{j}_{a,b-1,c} + D _6 \mathcal{L} ^{j-1}_{a,b-1,c}+ D _7 \mathcal{L} ^{j}_{a,b,c-1}\\
& \qquad     + D _8 \mathcal{L} ^{j-1}_{a,b,c-1} + D _9 \mathcal{L} ^{j}_{a-1,b-1,c} + D _{10} \mathcal{L} ^{j-1}_{a-1,b-1,c} + D _{11} \mathcal{L} ^{j}_{a,b-1,c-1} + D _{12} \mathcal{L} ^{j-1}_{a,b-1,c-1}\\
&  \qquad  +  D _{13} \mathcal{L} ^{j}_{a-1,b,c-1} + D _{14} \mathcal{L} ^{j-1}_{a-1,b,c-1}  + D _{15} \mathcal{L} ^{j}_{a-1,b-1,c-1} + D _{16} \mathcal{L} ^{j-1}_{a-1,b-1,c-1}  =0
\end{aligned}}
\end{equation}
for both $\mathcal{L}^\mathsf{e}_{d}$ and $\mathcal{L}^\mathsf{f}_{d}$, where we used the abbreviate notation $ \mathcal{L} ^j_{a,b,c}:= \mathcal{L} (j^1 \varphi _d (\mbox{\mancube}^j_{a,b,c}))$. For the discrete Lagrangians based on the discretizations given in \eqref{3D_kinetic_energy}--\eqref{3D_grav_potential}, the discrete Euler-Lagrange equations \eqref{discrete_Euler_Lagrange} yield
\begin{equation}\label{concrete_EL_3D}{\fontsize{10pt}{13pt}\selectfont
\hspace{-0.4cm}\begin{aligned} 
& \rho  _0^\mathsf{n} \frac{\mathbf{v}^{j-1}_{a,b,c}- \mathbf{v}^j_{a,b,c}}{ \Delta t} -  \frac{1}{8}\sum_{\ell=1}^4 \Big[ \mathfrak{D}_1^\ell   \left( \mathbf{F} _\ell(\mbox{\mancube}_{a,b,c}^j) \mathbf{S}_\ell(\mbox{\mancube}_{a,b,c}^j)\right)+\mathfrak{D}_3^\ell   \left( \mathbf{F} _\ell(\mbox{\mancube}_{a-1,b,c}^j) \mathbf{S}_\ell(\mbox{\mancube}_{a-1,b,c}^j)\right)\\
 &  +\mathfrak{D}_5^\ell   \left( \mathbf{F} _\ell(\mbox{\mancube}_{a,b-1,c}^j) \mathbf{S}_\ell(\mbox{\mancube}_{a,b-1,c}^j)\right)+ ... + \mathfrak{D}_{15}^\ell   \left( \mathbf{F} _\ell(\mbox{\mancube}_{a-1,b-1,c-1}^j) \mathbf{S}_\ell(\mbox{\mancube}_{a-1,b-1,c-1}^j)\right) \Big]\\
 &  - \rho  _0^\mathsf{n} \frac{\partial \Pi}{\partial \varphi }    ( \varphi ^j _{a,b})=0,
\end{aligned}}
\end{equation}  
$\mathsf{n}=\mathsf{e},\mathsf{f}$. The derivation of \eqref{concrete_EL_3D}  is similar to the 2D case, see Appendix \ref{discrete_EL_derivation}. These equations give the dynamics at the interior of the elastic body and the fluid.

\begin{figure}[H] \centering 
\includegraphics[width=3.2 in]{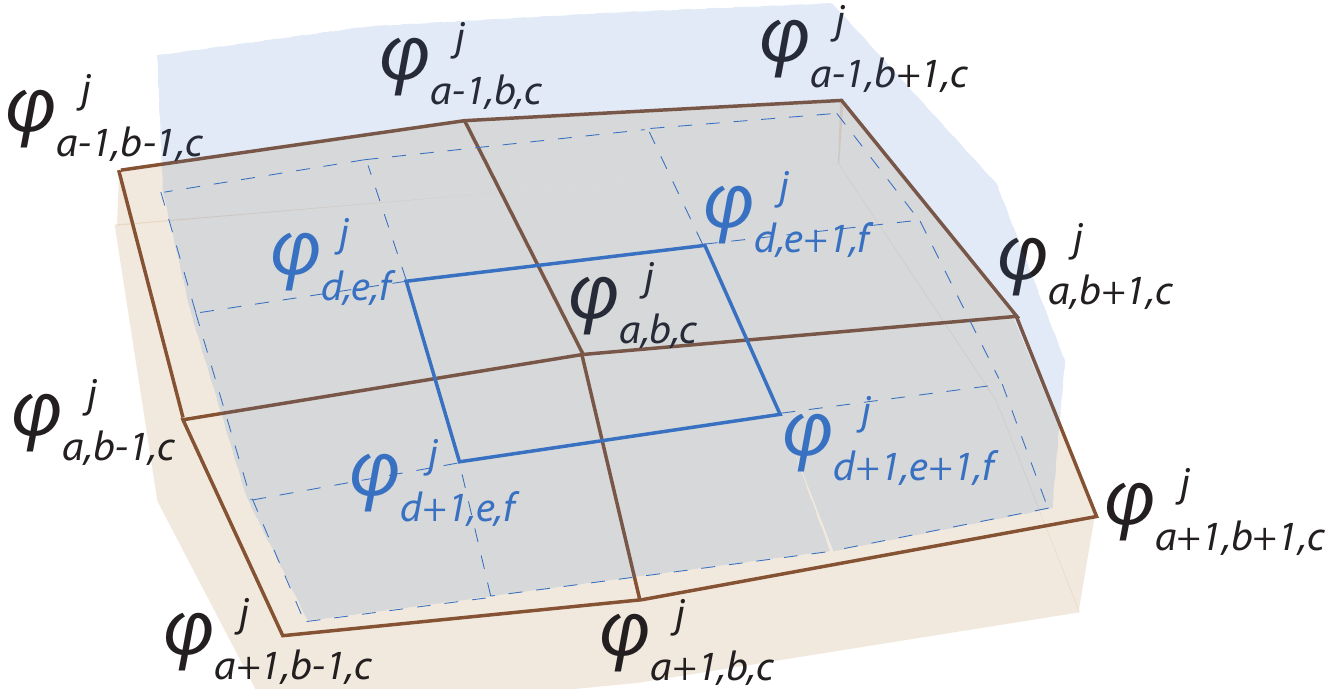} \vspace{-3pt}  
\caption{\footnotesize Fluid flowing above the hyperelastic body.} \label{3D_impenetrability} 
\end{figure}

\paragraph{Impenetrability conditions.} The impenetrability conditions between the fluid nodes $\varphi_{d,e,f}^j$, $\varphi_{d+1,e,f}^j$, $\varphi_{d,e+1,f}^j$, $\varphi_{d+1,e+1,f}^j$ and the elastic body nodes are defined as follows,  see Fig.\,\ref{3D_impenetrability},
\begin{equation} \label{3D_impenetrability_constraints}
\begin{aligned}
&\Psi_{{\rm im}_1}(\varphi_{a,b,c}^j, \varphi_{a-1,b,c}^j, \varphi_{a,b-1,c}^j, \varphi_{d,e,f}^j) 
\\
& \qquad = \big\langle (\varphi_{d,e,f}^j - \varphi_{a,b,c}^j ), (\varphi_{a-1,b,c}^j -\varphi_{a,b,c}^j) \times (\varphi_{a,b-1,c}^j -\varphi_{a,b,c}^j) \big\rangle \geq 0
\\
&\Psi_{{\rm im}_2}(\varphi_{a,b,c}^j, \varphi_{a+1,b,c}^j, \varphi_{a,b-1,c}^j, \varphi_{d+1,e,f}^j) 
\\
& \qquad = \big\langle (\varphi_{d+1,e,f}^j - \varphi_{a,b,c}^j ), (\varphi_{a,b-1,c}^j -\varphi_{a,b,c}^j) \times (\varphi_{a+1,b,c}^j -\varphi_{a,b,c}^j) \big\rangle \geq 0
\\
&\Psi_{\rm im_3}(\varphi_{a,b,c}^j, \varphi_{a-1,b,c}^j, \varphi_{a,b+1,c}^j, \varphi_{d,e+1,f}^j) 
\\
& \qquad = \big\langle (\varphi_{d,e+1,f}^j - \varphi_{a,b,c}^j ), (\varphi_{a,b+1,c}^j -\varphi_{a,b,c}^j) \times (\varphi_{a-1,b,c}^j -\varphi_{a,b,c}^j) \big\rangle \geq 0
\\
&\Psi_{\rm im_4}(\varphi_{a,b,c}^j, \varphi_{a+1,b,c}^j, \varphi_{a,b+1,c}^j, \varphi_{d+1,e+1,f}^j) 
\\
& \qquad = \big\langle (\varphi_{d+1,e+1,f}^j -\varphi_{a,b,c}^j ), (\varphi_{a+1,b,c}^j -\varphi_{a,b,c}^j) \times (\varphi_{a,b+1,c}^j -\varphi_{a,b,c}^j) \big\rangle \geq 0.
\end{aligned}
\end{equation}
The expression $(\varphi_{a-1,b,c}^j -\varphi_{a,b,c}^j) \times (\varphi_{a,b-1,c}^j -\varphi_{a,b,c}^j)$ represents the outward pointing normal vector field to the body, similarly in the other three conditions.

The associated penalty terms $\Phi_{{\rm im}_i}$, $i=1,...,4$ are defined as in \eqref{penalty_2D}. 
The directions of the reaction forces are given by the derivative with respect to the positions of the nodes, i.e.,
\begin{align*}
&D_1\Psi_{{\rm im}_1}(\varphi_{a,b,c}^j, \varphi_{a-1,b,c}^j, \varphi_{a,b-1,c}^j, \varphi_{d,e,f}^j), \quad D_2\Psi_{{\rm im}_1}(\varphi_{a,b,c}^j, \varphi_{a-1,b,c}^j, \varphi_{a,b-1,c}^j, \varphi_{d,e,f}^j),\\
&D_3 \Psi_{{\rm im}_1}(\varphi_{a,b,c}^j, \varphi_{a-1,b,c}^j, \varphi_{a,b-1,c}^j, \varphi_{d,e,f}^j), \quad D_4\Psi_{{\rm im}_1},(\varphi_{a,b,c}^j, \varphi_{a-1,b,c}^j, \varphi_{a,b-1,c}^j, \varphi_{d,e,f}^j),
\end{align*}
similarly for $\Psi_{{\rm im}_2}, \Psi_{{\rm im}_3}, \Psi_{{\rm im}_4}$, see Appendix \ref{3D_deriv_imp_const} for the expressions of these derivatives.

The resulting discrete action functional is
\begin{equation}\label{tilde_S_3D}{\fontsize{9pt}{13pt}\selectfont
\begin{aligned} 
&\widetilde{S}_d( \varphi _d ) \\
&= S( \varphi  _d ) - \sum_{j=0}^{N-1}\sum_{a=0}^{A-1}\sum_{b=0}^{B-1}\sum_{c=0}^{C-1} \Delta t \Delta s_1 \Delta s_2 \Delta s _3\,\Phi_{\rm in}\big(j^1 \varphi _d (\mbox{\mancube}^j_{a,b,c})\big) \\
& \;\;\; - \sum_{j=0}^{N-1}\!\sum_{(d,e,f)\, \in \,\mathbb{B}^\mathsf{f}}\!\Big[ \Phi_{{\rm im}_1}(\varphi_{a,b,c}^j, \varphi_{a-1,b,c}^j, \varphi _{a,b-1,c},\varphi_{d,e,f}^j)+ \Phi_{{\rm im}_2}(\varphi_{a,b,c}^j, \varphi_{a-1,b,c}^j, \varphi _{a,b-1,c},\varphi_{d,e,f}^j)\\
&\hspace{1cm}+\Phi_{{\rm im}_3}(\varphi_{a,b-1,c}^j, \varphi_{a-1,b-1,c}^j, \varphi _{a,b,c},\varphi_{d,e,f}^j)+ \Phi_{{\rm im}_4}(\varphi_{a-1,b-1,c}^j, \varphi_{a,b-1,c}^j, \varphi _{a-1,b,c},\varphi_{d,e,f}^j)\Big],
\end{aligned}}
\end{equation} 
where $S_d$ is given in \eqref{3D_discrete_action_sum}. The second term in \eqref{tilde_S_3D} is associated to the incompressibility of the elastic body. In the third term, the second sum is taken over the boundary nodes of the fluid and, given such a boundary node $ \varphi _{d,e,f}^j$, the nodes of the elastic body in the four terms $i=1,...,4$ are chosen in accordance with \eqref{3D_impenetrability_constraints} for each functions $ \Phi _{{\rm im}_i}$, $i=1,...,4$.

\paragraph{Symmetries and discrete Noether theorems.}
The discrete classical momentum map associated to the action of the special Euclidean group $SE(3)$ takes the form
\begin{equation}\label{3D_momentum_map_fluid}
\mathbf{J} _d ( \boldsymbol{\varphi} ^j, \boldsymbol{\varphi} ^{j+1}) =\sum_{a=0}^{A-1} \sum_{b=0}^{B-1}\sum_{c=0}^{C-1} \mathbf{J} ^\mathsf{e}\big(j^1 \varphi _d (\mbox{\mancube}_{a,b,c}^j)\big) +\sum_{d=0}^{D-1} \sum_{e=0}^{E-1}\sum_{f=0}^{F-1}  \mathbf{J} ^\mathsf{f}\big(j^1 \varphi _d (\mbox{\mancube}_{d,e,f}^j)\big)
\end{equation}
with
\begin{equation}
\mathbf{J} ^\mathsf{e}\big(j^1 \varphi _d (\mbox{\mancube}_{a,b,c}^j)\big)= \sum_{\alpha =a}^{a+1} \sum_{\beta =b}^{b+1}\sum_{\gamma  =c}^{c+1} \left( \varphi_{\alpha,\beta}^j \times \frac{M^\mathsf{e}}{8} v_{\alpha,\beta, \gamma }^j ,  \frac{M^\mathsf{e}}{8} v_{\alpha,\beta, \gamma }^j\right) 
\end{equation}
similarly for the fluid.

\subsection{Barotropic fluid flowing over a Mooney-Rivlin hyperelastic support}\label{3D_Mooney_Rivlin}

\begin{figure}[H] \centering 
\qquad \includegraphics[width=1.6 in]{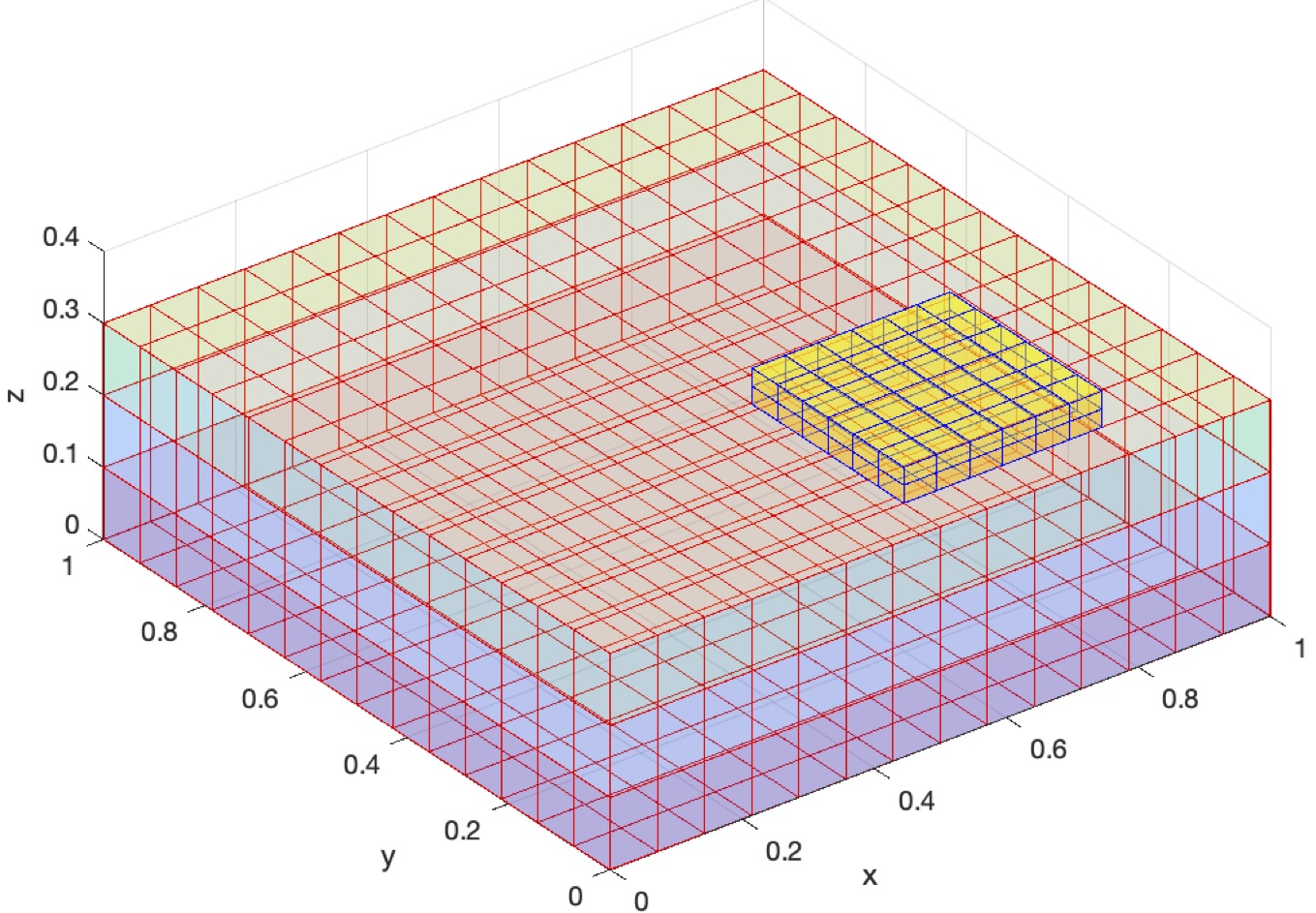} \hspace{0.15cm} \hspace{0.15cm} \qquad \qquad \hspace{1.7 in} \vspace{-1pt}
\\
\hspace{0.6cm} \includegraphics[width=2.35 in]{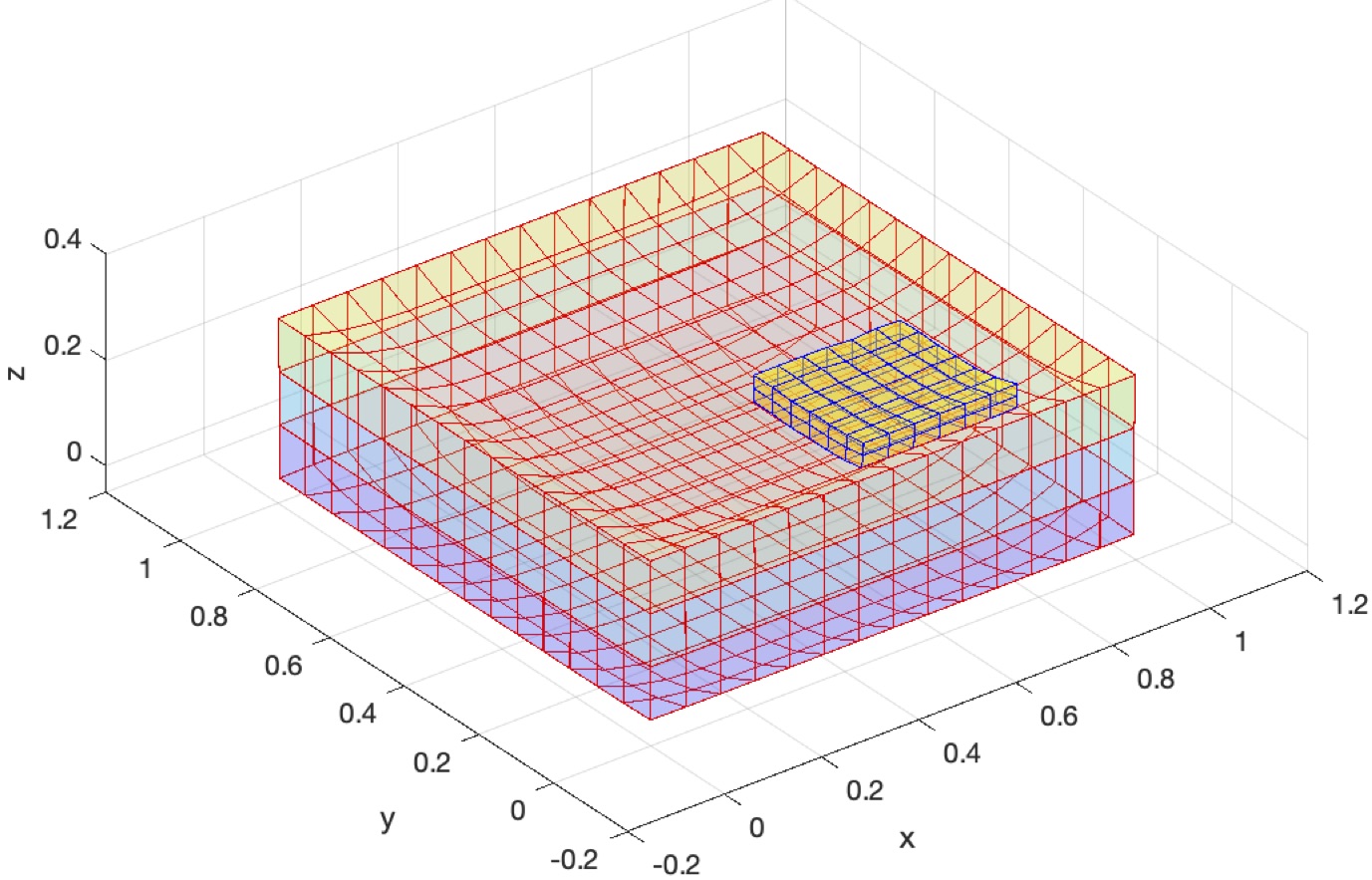}  \hspace{0.15cm} \qquad \qquad \includegraphics[width=1.7 in]{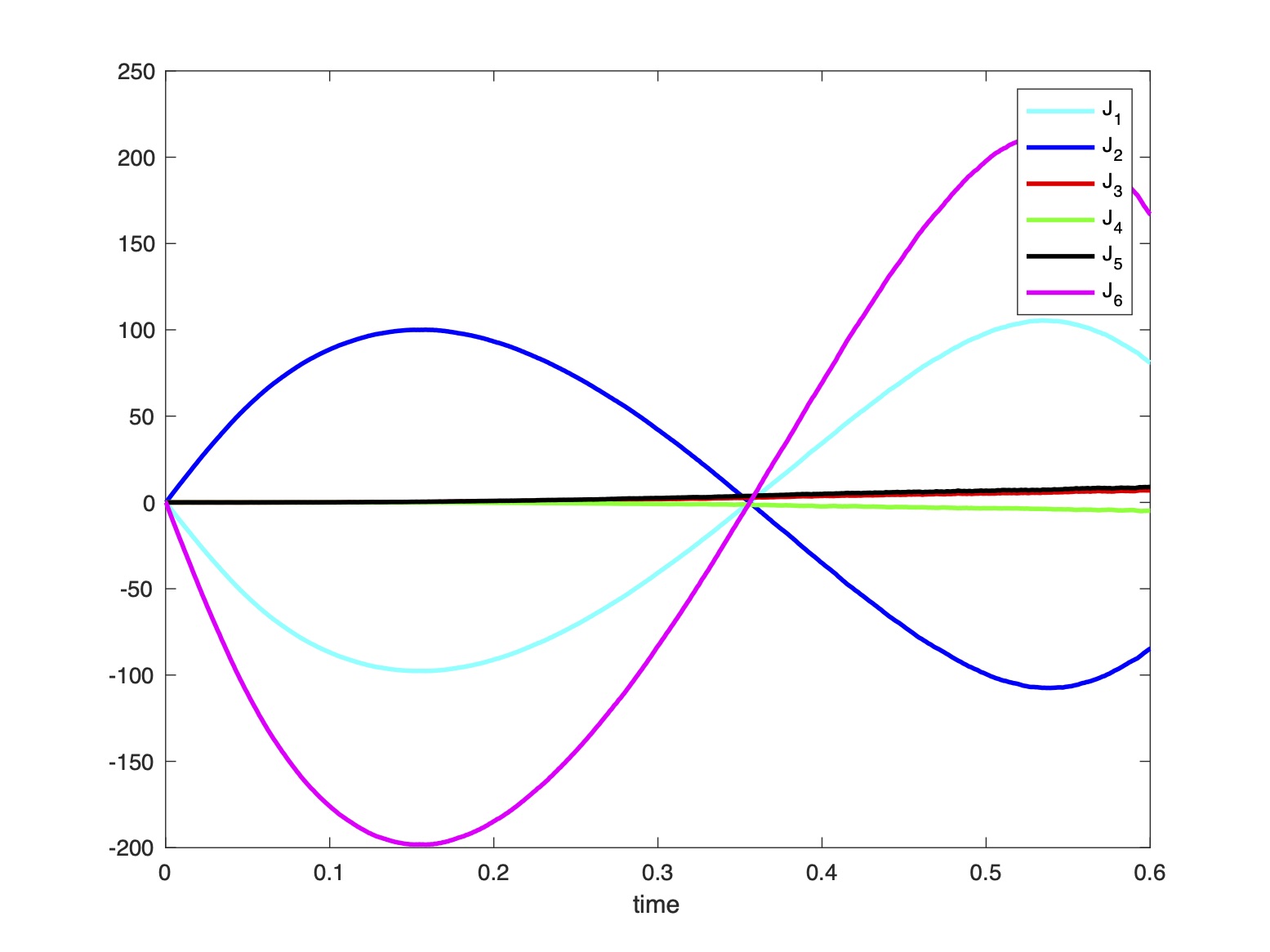} \vspace{-1pt}
\\
\hspace{1 cm} \includegraphics[width=1.95 in]{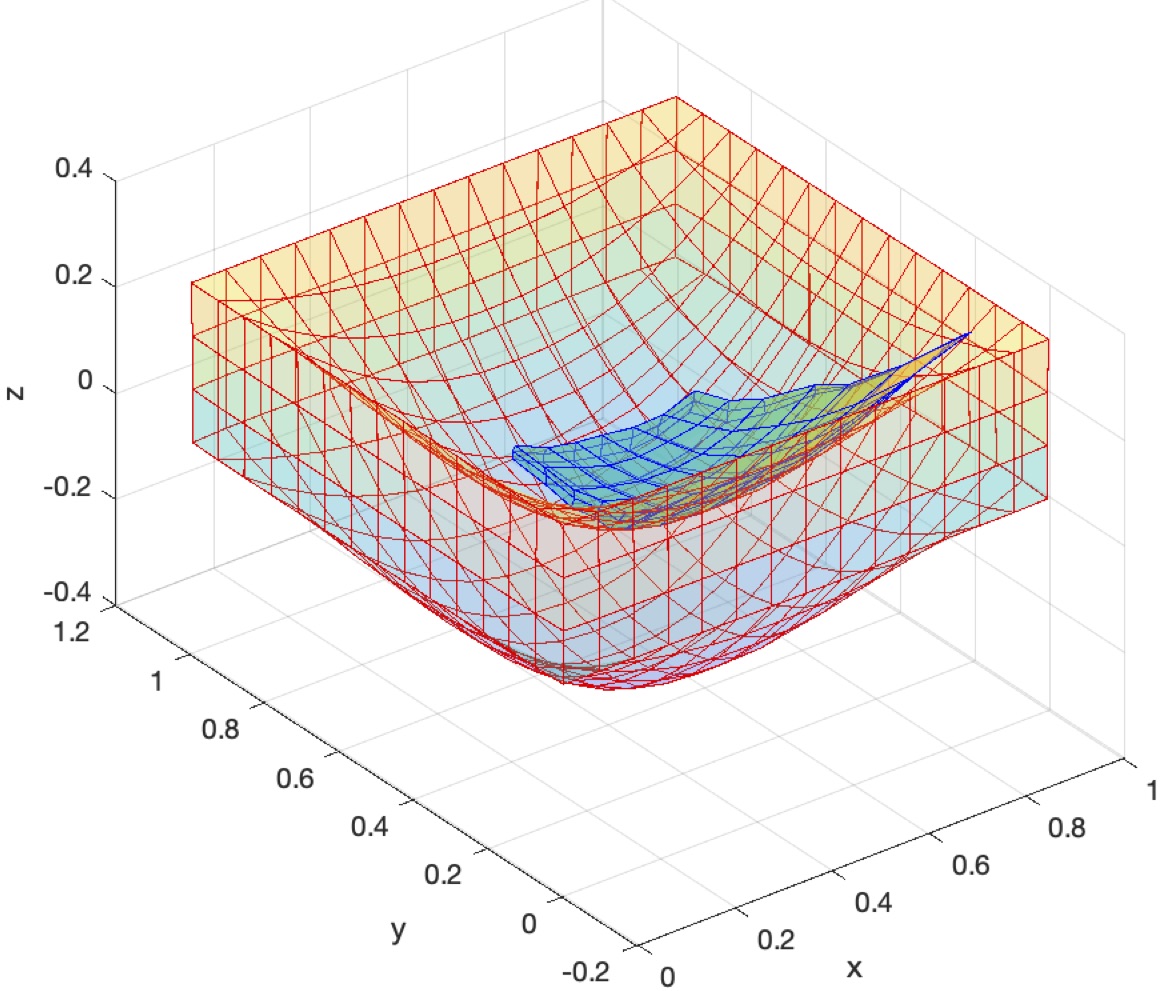} \hspace{0.9cm} \qquad \qquad \includegraphics[width=1.7 in]{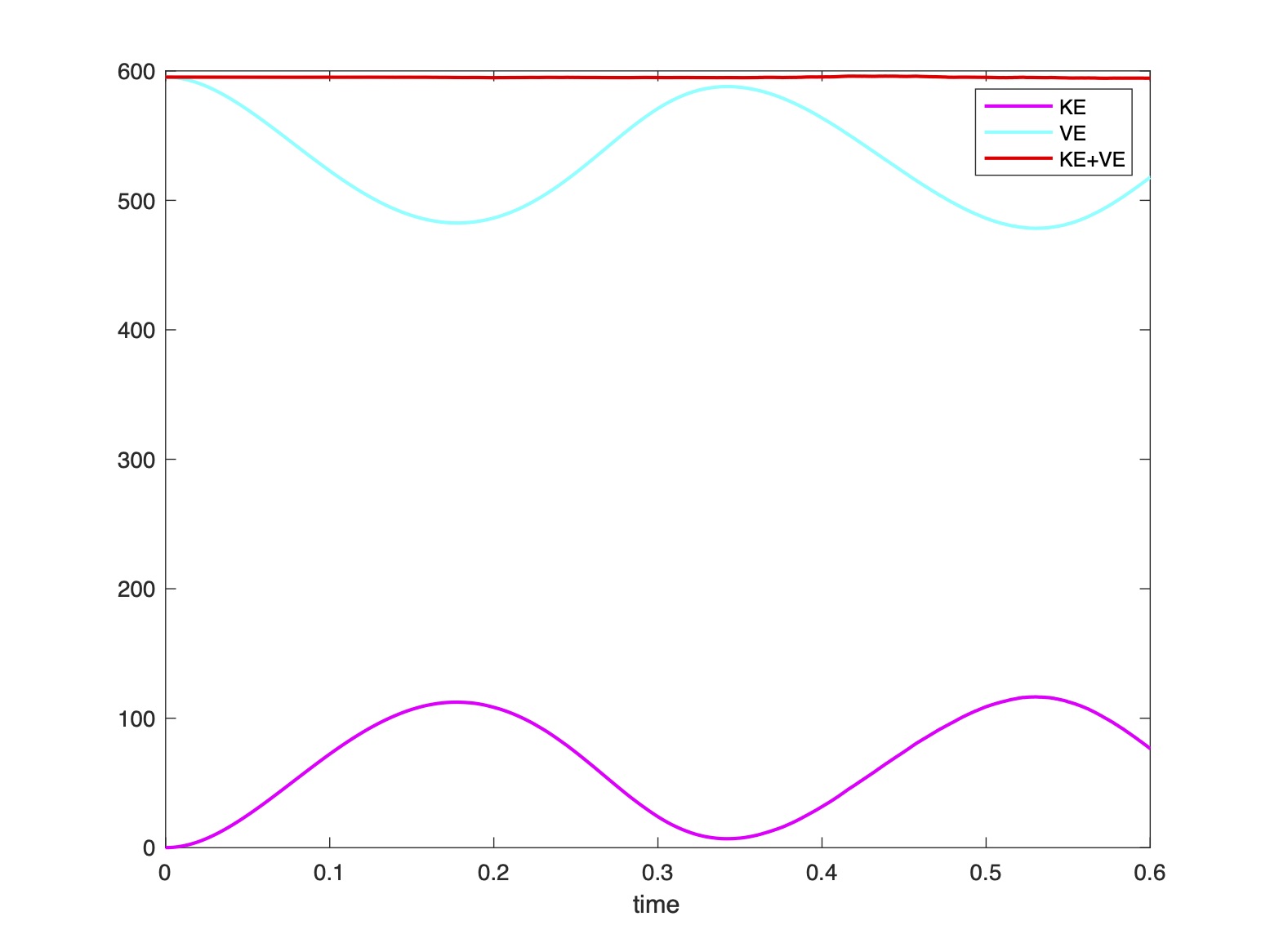} \vspace{-1pt}
\\
\hspace{0.6cm}  \includegraphics[width=2.4 in]{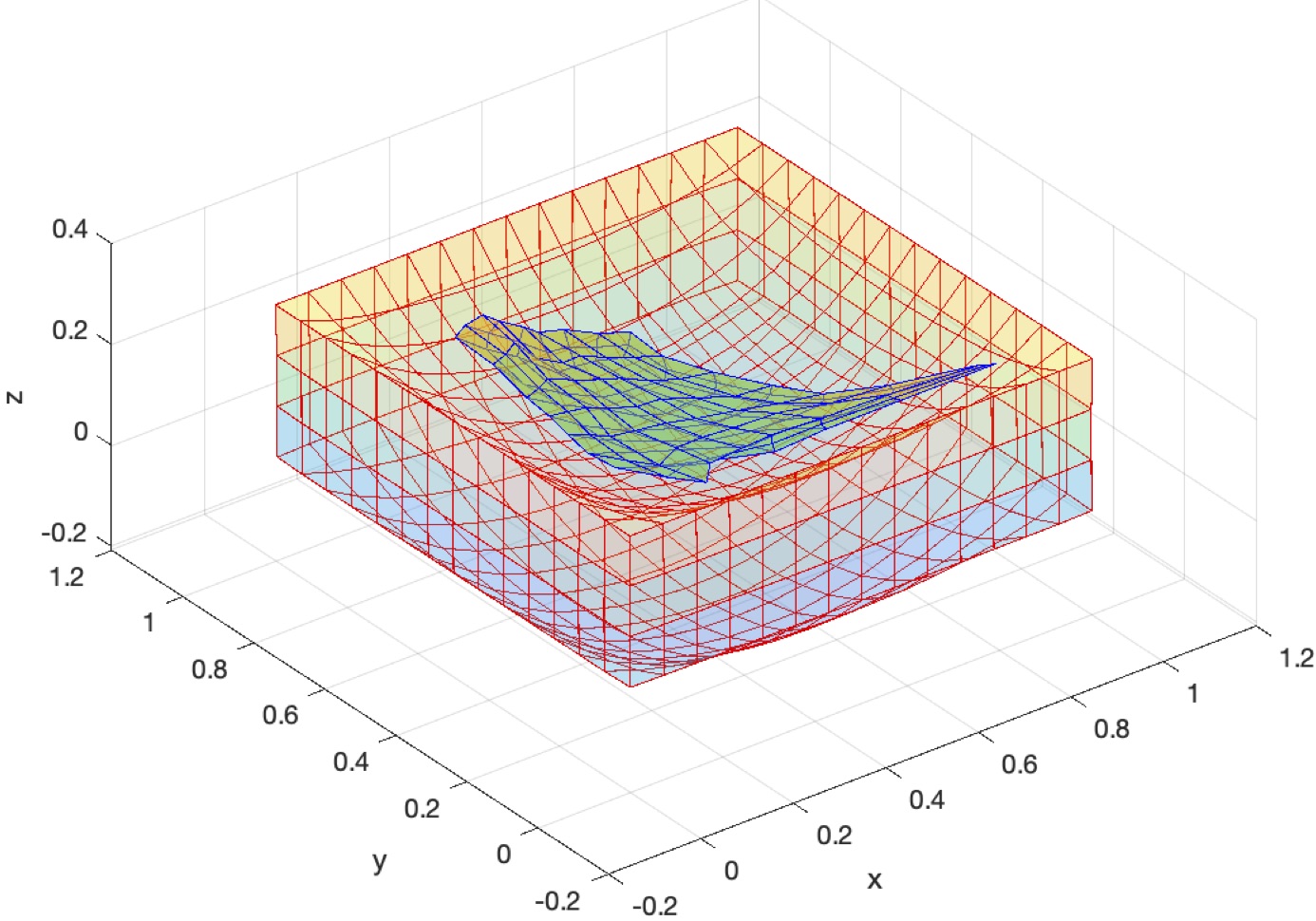} \hspace{0.1cm} \qquad \qquad  \includegraphics[width=1.7 in]{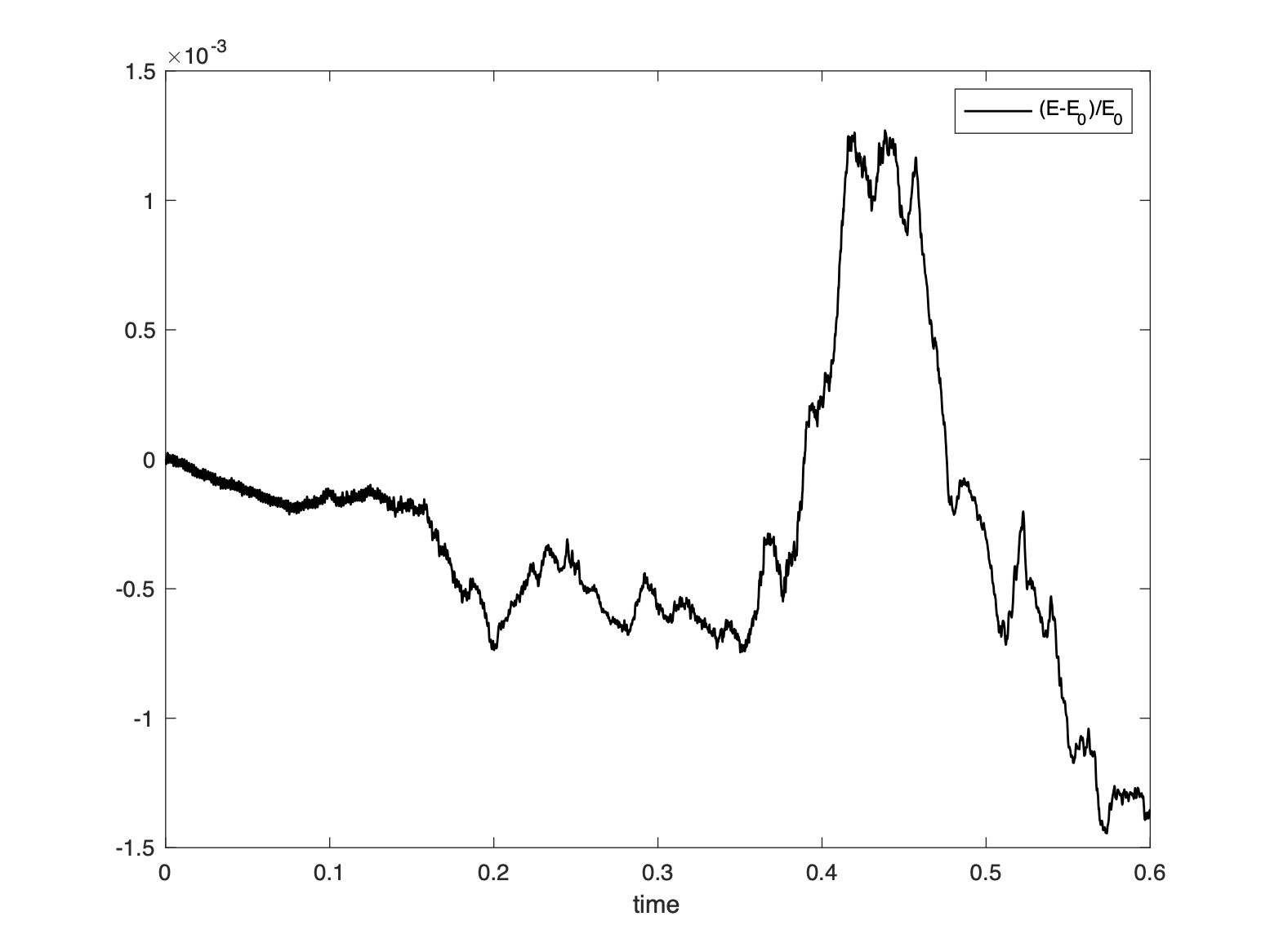} \vspace{-1pt}
\\
\hspace{0.6cm}  \includegraphics[width=2.2 in]{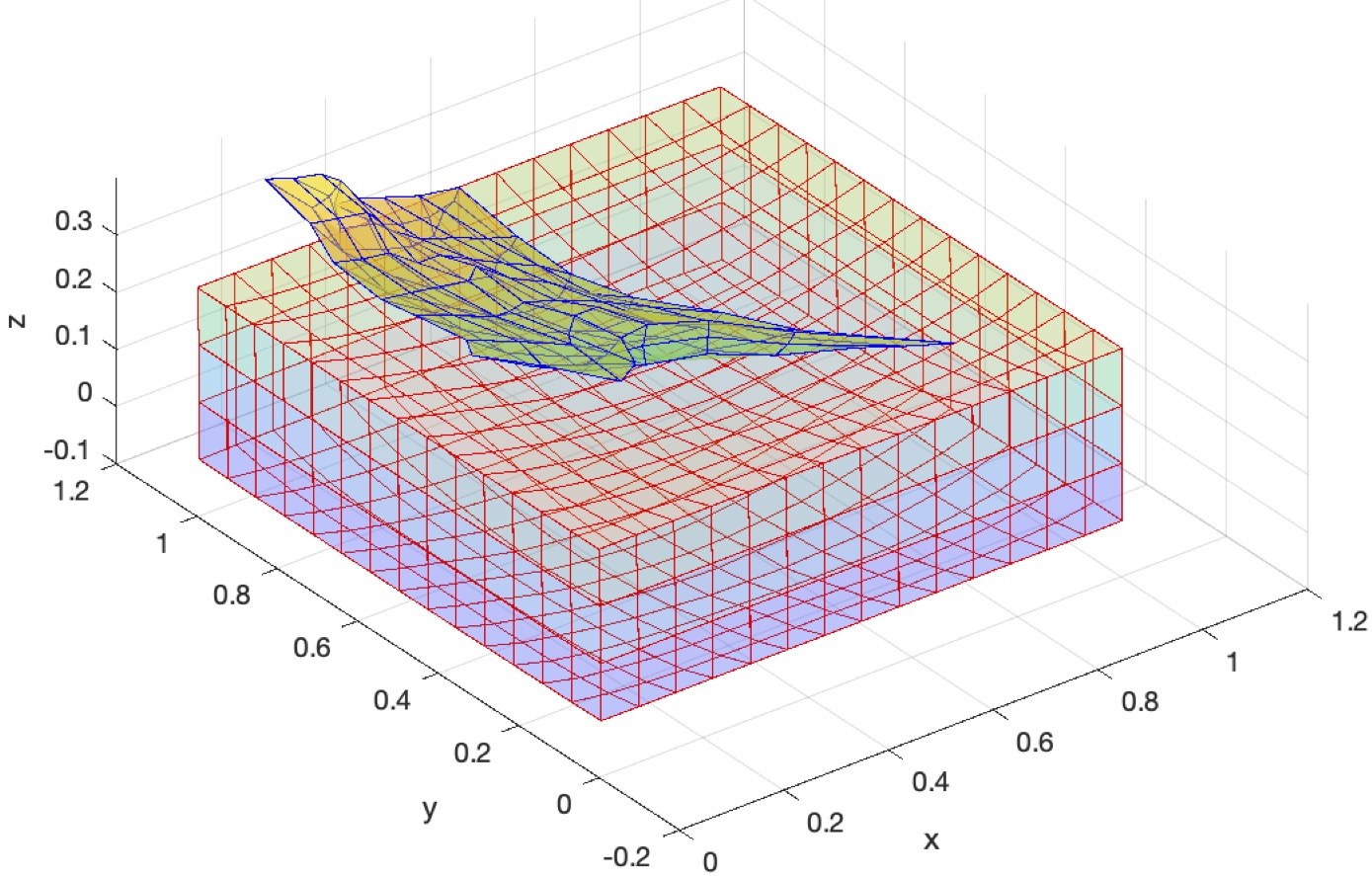} \hspace{0.15cm} \qquad \qquad \hspace{1.7 in}\vspace{-1pt}
\caption{\footnotesize \textit{Top to bottom on the left}: Configuration after $0.002$s, $0.1$s, $0.4$s, $0.5$s, $0.6$s. \textit{Top to bottom on the right}: Evolution of total (fluid+solid) momentum maps, total (fluid+solid) energy, and relative total energy during $0.5$s.}\label{fluid_elastic_body_3D}
\end{figure}

We consider an incompressible Mooney-Rivlin hyperelastic body, see \eqref{3D_disc_Mooney_Rivlin}, and a barotropic fluid described by the Tait equation with discrete internal energy term
\begin{equation}\label{Dis_Tait_eq_2D} 
W_d^\mathsf{f}\big( \rho  _0, j^1 \varphi _d( \mbox{\mancube})\big) =  \frac{1}{8}  \sum_{\ell=1}^8\left[ \frac{A}{ \gamma -1} \left(\frac{J_\ell (\mbox{\mancube})}{ \rho  _0}  \right)    ^{1- \gamma } + B  \left(\frac{J_\ell (\mbox{\mancube})}{ \rho  _0}  \right) \right].
\end{equation}
The integrator is found by computing the criticality condition for the discrete action functional which includes the incompressibility and impenetrability penalty terms. We simulate the following situation, see Fig.\,\ref{fluid_elastic_body_3D}:
\begin{enumerate}
\item The barotropic fluid has the properties $\rho_0= 997 \, \mathrm{kg/m}^2$, $\gamma =6 $, $A = \tilde{A} \rho_0^{-\gamma}$ with $\tilde A = 3.041\times 10^4$ Pa, and $B = 3.0397\times 10^4$ Pa. The size of the discrete reference configuration at time $t^0$ is $0.3\mathrm{m} \times 0.3 \mathrm{m} \times 0.05 \mathrm{m}$, with space-steps $\Delta s_1=\Delta s_2=0.05$, $\Delta s_3= 0.025$m.  

\item The Mooney-Rivlin model has the properties $\rho_0= 945 \, \mathrm{kg/m}^3$, $C_1= 1.848$, and $C_2= 0.264$ such that $C_1/C_2=7$. The model is nearly incompressible, with coefficient $r=10^4$ in the corresponding penalty term. The size of the discrete reference configuration at time $t^0$ is $1\mathrm{m} \times 1 \mathrm{m} \times 0.3 \mathrm{m}$ with space-steps $\Delta s_1=\Delta s_2= 0.0625$m, $\Delta s_3= 0.1$m.

\item The four upper edges of the Mooney-Rivlin incompressible model are fixed.
\item The time-step is $\Delta t=2\times10^{-4}$ and the test is carried out for $0.4s$.
\end{enumerate}
The results are reproduced in Fig.\,\ref{fluid_elastic_body_3D}, where we used the notation $\mathbf{J}_d=(\mathsf{J}_1,\mathsf{J}_2, \mathsf{J}_3,\mathsf{J}_4,\mathsf{J}_5,\mathsf{J}_6) \in \mathfrak{se}(3)^*$. Due to the presence of the gravitation potential, the discrete Lagrangian is only invariant under the subgroup $SE(2)$ of Euclidean transformations of the horizontal plane. This is why only the components $\mathsf{J}_3, \mathsf{J}_4, \mathsf{J}_5$ of the discrete momentum map $ \mathbf{J} _d$ are almost preserved. As in our 2D test, the fixing of the four upper edges prevents the exact conservation of these momentum maps.

\subsection{Convergence tests for the isotropic nonlinear hyperelastic model}

\begin{figure}[H] \centering 
\includegraphics[width=2.1 in]{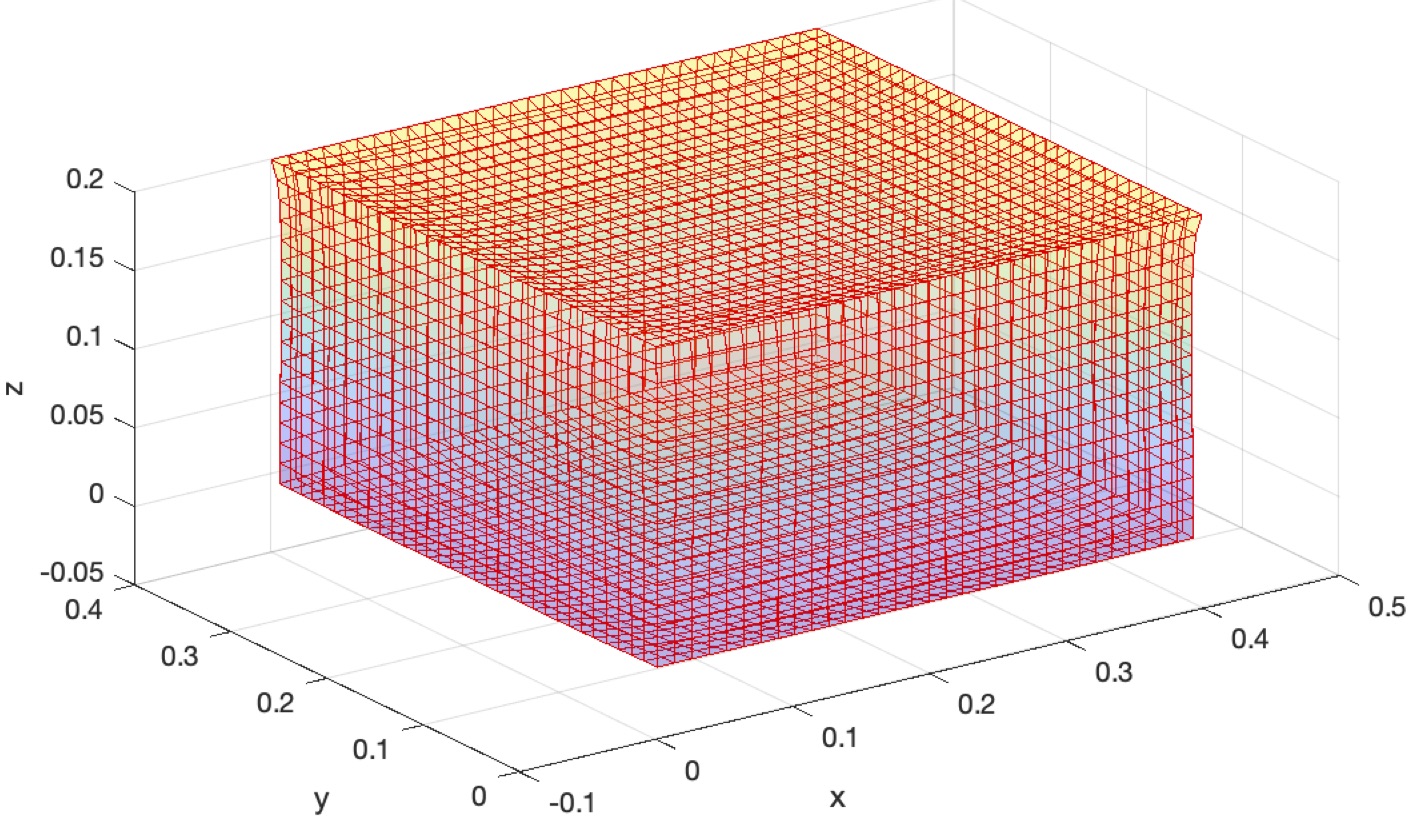} \vspace{-3pt}   \includegraphics[width=1.5 in]{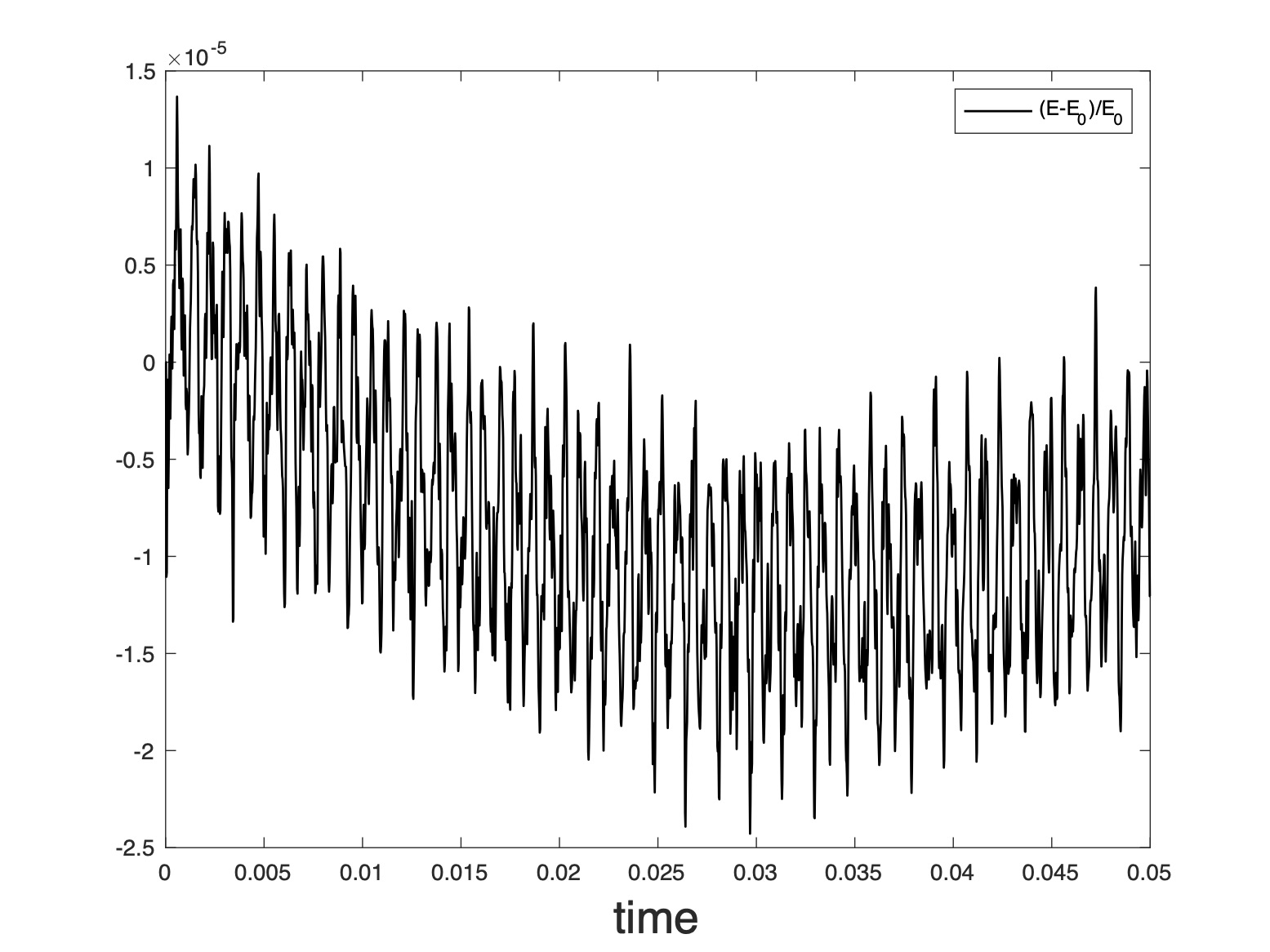} \vspace{-3pt}   \includegraphics[width=1.5 in]{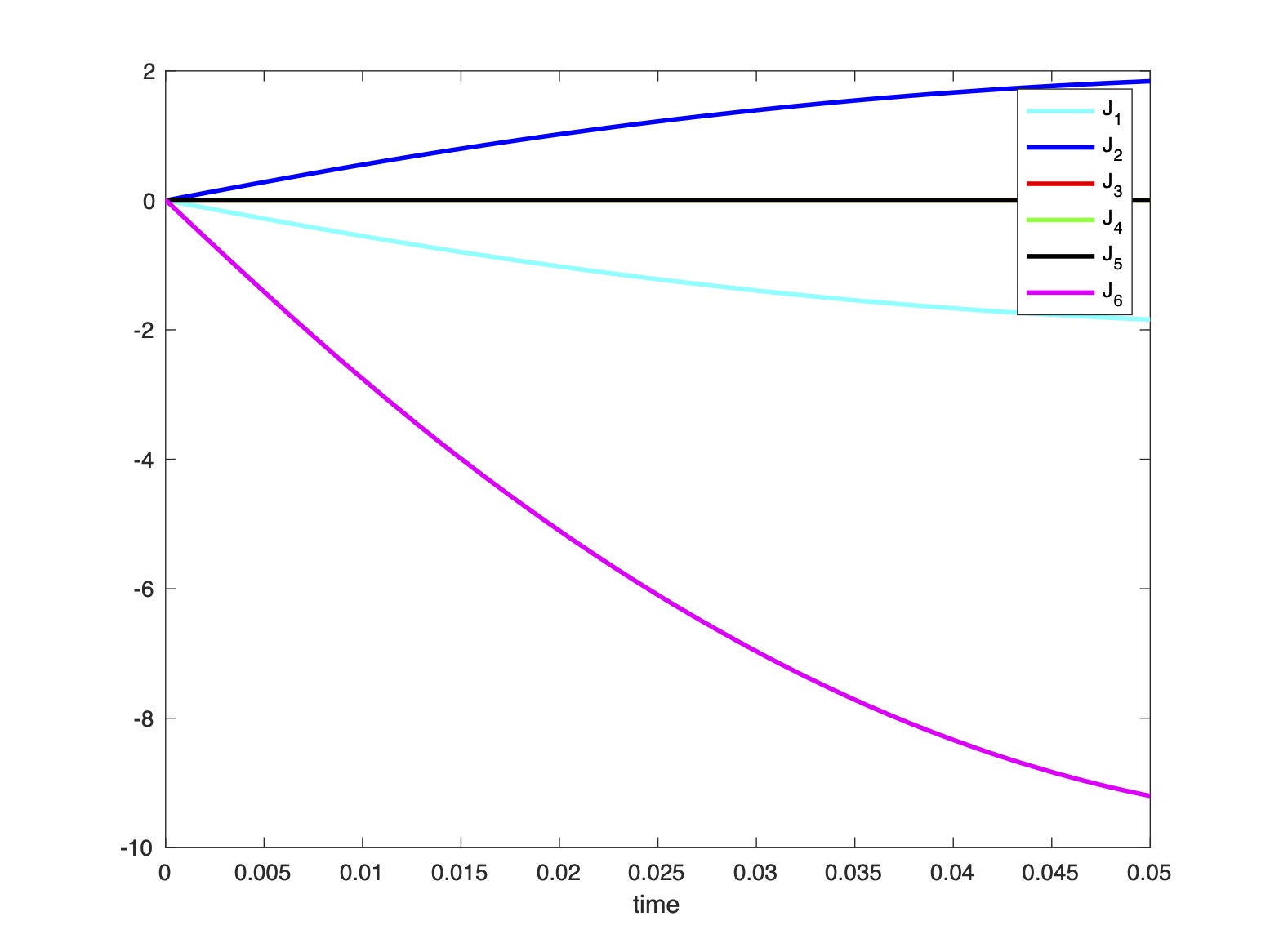} \vspace{-3pt}  
\caption{\footnotesize \textit{From left to right}: Configuration, relative energy, and momentum maps after $0.05$s, when $\Delta t=2\times 10^{-5}$ and $\Delta s_1=\Delta s_2=\Delta s_3=1.25\times 10^{-2}$. }\label{convergence_3D} 
\end{figure}

Consider the Mooney-Rivlin incompressible model described above, with properties $\rho_0= 945 \, \mathrm{kg/m}^3$, $C_1= 1.848$, and $C_2= 0.264$ with $C_1/C_2=7$. The model is incompressible, with coefficient $r=10^4$ in the penalty term.
The size of the discrete reference configuration at time $t^0$ is $0.4\mathrm{m} \times 0.4 \mathrm{m} \times 0.2 \mathrm{m}$. As earlier, the four upper edges of the Mooney-Rivlin incompressible model are fixed, see Fig.\,\ref{convergence_3D}. The motion is only induced by gravity forces. 

We consider the explicit integrator described in \eqref{concrete_EL_3D}, for the elastic body only, and study the convergence with respect to $\Delta t$ and $\Delta s_i$, $i=1,2,3$.

\begin{enumerate}
\item
 Given a fixed mesh, with values $\Delta s_1=\Delta s_2=\Delta s_3= 0.1$m, we vary the time-steps as $\Delta t \in \{ 2 \times 10^{-4}, \,  10^{-4}, \, 5 \times 10^{-5}, \, 2.5 \times 10^{-5} \}$.
We compute the $L^2$-errors in the position $\varphi_d$ at time $t^N=0.1$s, by comparing $\varphi_d$ with an ``exact solution'' obtained with the time-step $\Delta t_{\rm ref}=6.25 \times 10^{-6}$s.  That is, for each value of $\Delta t$ we calculate
\begin{equation}\label{L_2norm_3D3}
\| \varphi_d - \varphi_{\rm ref} \|_{L^2} = \left( \sum_a \sum_b\sum_c  \| \varphi_{a,b,c}^N - \varphi_{{\rm ref};a,b,c}^N \|^2 \right)^{1/2}.
\end{equation}
This yields the following convergence with respect to $\Delta t$ 
\begin{figure}[H] \centering 
\begin{tabular}{| c | c | c | c | c |}
\hline
$\Delta t$ & $2 \times 10^{-4}$ & $ 10^{-4}$ & $ 5 \times 10^{-5}$ & $2.5 \times 10^{-5}$  \\
\hline
$\| \varphi_d - \varphi_{\rm ref} \|_{L^2}$ &  $5.23 \times 10^{-4}$  & $2.55 \times 10^{-4}$  & $1.18\times 10^{-4}$  &  $5.035 \times 10^{-5}$ \\
\hline
$ \text{rate} $  &    & 1.036 & 1.12  & 1.29  \\
\hline
\end{tabular}
\end{figure}

\item Given a fixed time-step $\Delta t=2 \times 10^{-5}$s, we vary the space-steps as $\Delta s_1 =\Delta s_2=\Delta s_3$ $\in \{0.1,\, 0.05, \,0.025, \, 0.0125 \}$. The ``exact solution'' is chosen with $\Delta s_{1;\rm ref}=\Delta s_{2;\rm ref}  =\Delta s_{3;\rm ref}  = 0.00625$m.  We compute the $L^2$-errors in the position $\varphi_d$ at time $t^N=0.05$s. We get the following convergence with respect to $\Delta s_1=\Delta s_2=\Delta s_3$.
\begin{figure}[H] \centering 
\begin{tabular}{| c | c | c | c |c |c|c|}
\hline
$\Delta s_1=\Delta s_2 =\Delta s_3 $ & $0.1$  & $0.05$  & $0.025$    & $0.0125$\\
\hline
$\| \varphi_d - \varphi_{\rm ref} \|_{L^2}$ & $0.0375$   &  $0.0284 $ &  $0.0177$ &  $0.008$ \\
\hline
$ \text{rate} $  &  & 0.40  & 0.682 & 1.1457\\
\hline
\end{tabular}
\end{figure}
\end{enumerate}

In Fig.\,\ref{convergence_3D} on the right, the components $\mathsf{J}_3, \mathsf{J}_4, \mathsf{J}_5$ of the discrete momentum map are preserved, which corresponds to the symmetry subgroup $SE(2)$.

\section{Concluding remarks and future directions} \label{conclusion}

In this paper we introduced a structure preserving numerical integrator for hyperelastic solids, which is spacetime multisymplectic, symplectic in time, preserves exactly the momenta associated to symmetries, and nearly preserves total energy. The integrator is derived from a discrete variational formulation of hyperelasticity in the Lagrangian description. Thanks to its variational nature and its adequacy with a previously derived variational discretization for fluids, the scheme can naturally be extended to accommodate fluid-structure interaction problems via the inclusion of impenetrability penalty terms in the discrete action functional. Important steps in our approach were the definition of discrete deformation gradients, discrete Cauchy-Green deformation tensors and discrete Jacobians. The resulting integrator and its conservative properties were illustrated with numerical tests both in 2D and 3D. This paper only presented the very first steps towards the application of variational integrators to fluid-structure interaction problems. The discrete setting developed here opens several directions of research such as the inclusion of dynamic mesh update, inspired by arbitrary Lagrangian-Eulerian methods, and the treatment of h-adaptivity, by appropriate extensions of the variational discretization. One can also explore extensions of our model to take into account of the possible disconnection of fluid cells, which is useful to model the fragmentation of the flow after the impact. Finally, appropriate extensions of the discrete variational formulation need to be developed in order to treat other classes of continuum models, such as viscous-plastic materials.

\appendix 

\section{Appendix}\label{appendix}

\subsection{2D discrete hyperelastic models}\label{discrete_EL_derivation}

To compute the partial derivatives $D_k \mathcal{L} ^j_{a,b}$ of the discrete Lagrangian \eqref{Discrete_Lagrangian_2D_wave} we note
\[
D_k\left( \rho  _0^\mathsf{e} \frac{1}{4}  \sum_{\ell=1}^4 W^\mathsf{e}\big( \mathbf{C}_\ell(\mbox{\mancube}_{a,b}^j)\big) \right) \cdot \delta \varphi  =  \frac{1}{4} \sum_{\ell=1}^4 \operatorname{Tr} \left( \mathbf{S}_\ell(\mbox{\mancube}_{a,b}^j) \mathbf{F} _\ell(\mbox{\mancube}_{a,b}^j)^\mathsf{T} D_k  \mathbf{F} _\ell(\mbox{\mancube}_{a,b}^j) \cdot \delta \varphi \right) 
\]
for all $\ell=1,...,4$ and all $ k=1,3,5,7$. We observe that the 16 expressions $D_k \mathbf{F} _\ell(\mbox{\mancube}_{a,b}^j)$ are independent of the value of the field $ \varphi _d $ and on the location of the spacetime node (if the spacing $ \Delta  s_a$ and $ \Delta  s_b$ are constant), i.e., we can write $D_k \mathbf{F} _\ell(\mbox{\mancube}_{a,b}^j)=D_k \mathbf{F} _\ell$. For instance, for $\ell=1$, we have 
\begin{align*} 
D_1 \mathbf{F} _1 \cdot \delta \varphi &= \left[ - \frac{1}{ \Delta s_a} \delta \varphi  , - \frac{1}{ \Delta s_b} \delta \varphi \right] & D_3 \mathbf{F} _1 \cdot \delta \varphi &= \left[ \frac{1}{ \Delta s_a} \delta \varphi  , 0  \right]  \\
D_5 \mathbf{F} _1 \cdot \delta \varphi &= \left[ 0 ,  \frac{1}{ \Delta s_b} \delta \varphi \right] & D_7 \mathbf{F} _1 \cdot \delta \varphi &= \left[ 0, 0 \right] 
\end{align*} 
and for $\ell=2$, we have
\begin{align*} 
D_1 \mathbf{F} _2 \cdot \delta \varphi &= \left[ 0 ,  \frac{1}{ \Delta s_a} \delta \varphi \right] &D_3 \mathbf{F} _2 \cdot \delta \varphi &= \left[ - \frac{1}{ \Delta s_b} \delta \varphi  , - \frac{1}{ \Delta s_a} \delta \varphi   \right]\\
D_7 \mathbf{F} _2 \cdot \delta \varphi &= \left[  \frac{1}{ \Delta s_b} \delta \varphi  , 0 \right] &D_5 \mathbf{F} _2 \cdot \delta \varphi &= \left[ 0 , 0 \right],
\end{align*} 
similarly for $\ell=3,4$.
We can thus write
\begin{align*} 
D_k\left( \rho  _0^\mathsf{e} \frac{1}{4}  \sum_{\ell=1}^4 W^\mathsf{e}\big( \mathbf{C}_\ell(\mbox{\mancube}_{a,b}^j)\big) \right) \cdot \delta \varphi  &=\frac{1}{4} \sum_{\ell=1}^4  \left\langle   \mathbf{F} _\ell(\mbox{\mancube}_{a,b}^j) \mathbf{S}_\ell(\mbox{\mancube}_{a,b}^j) , D_k \mathbf{F} _\ell \cdot \delta \varphi  \right\rangle \\
&=\frac{1}{4} \sum_{\ell=1}^4  \mathfrak{D}_k^\ell  \left( \mathbf{F} _\ell(\mbox{\mancube}_{a,b}^j) \mathbf{S}_\ell(\mbox{\mancube}_{a,b}^j)  \right)  \cdot \delta \varphi ,
\end{align*} 
where in the last equality we defined the operator $\mathfrak{D}_k^\ell$, $k=1,3,5,7$, $\ell=1,2,3,4$, by duality.

From this, the partial derivatives $D_k \mathcal{L} _{a,b}^j$ of the discrete Lagrangian can be computed. For instance for $k=1$ one has
\[
\frac{1}{ \operatorname{vol}(\mbox{\mancube}) } D _1 \mathcal{L} ^{j}_{a,b}= - \rho  _0^\mathsf{e} \frac{1}{4 \Delta t}\mathbf{v}^j_{a,b} - \frac{1}{4}\sum_{\ell=1}^4\mathfrak{D}_1^\ell   \left( \mathbf{F} _\ell(\mbox{\mancube}_{a,b}^j) \mathbf{S}_\ell(\mbox{\mancube}_{a,b}^j)\right)   - \frac{1}{4\Delta t} \rho  _0^\mathsf{e} \nabla \Pi  ( \varphi ^j _{a,b}).
\]
Using this result, one then derives \eqref{concrete_EL} from \eqref{discrete_Euler_Lagrange}.  

\subsection{3D discrete Jacobian} \label{3D_Jacobian's}

{\small
\begin{align*}
J_2(\mbox{\mancube}_{a,b,c}^j) &=   (\mathbf{F}_{2;a+1,b,c}^j \times \mathbf{F}_{4;a+1,b,c}^j)\cdot \mathbf{F}_{3;a+1,b,c}^j ,\\
J_3(\mbox{\mancube}_{a,b,c}^j) &=   (\mathbf{F}_{5;a,b+1,c}^j \times \mathbf{F}_{1;a,b+1,c}^j)\cdot \mathbf{F}_{3;a,b+1,c}^j ,\\
J_4(\mbox{\mancube}_{a,b,c}^j) &= (\mathbf{F}_{2;a,b,c+1}^j \times \mathbf{F}_{1;a,b,c+1}^j)\cdot \mathbf{F}_{6;a,b,c+1}^j ,\\
J_5(\mbox{\mancube}_{a,b,c}^j) &=  (\mathbf{F}_{4;a+1,b+1,c}^j \times \mathbf{F}_{5;a+1,b+1,c}^j)\cdot \mathbf{F}_{3;a+1,b+1,c}^j,\\
J_6(\mbox{\mancube}_{a,b,c}^j) &= (\mathbf{F}_{1;a,b+1,c+1}^j \times \mathbf{F}_{5;a,b+1,c+1}^j)\cdot \mathbf{F}_{6;a,b+1,c+1}^j,\\ 
J_7(\mbox{\mancube}_{a,b,c}^j) & = (\mathbf{F}_{4;a+1,b,c+1}^j \times \mathbf{F}_{2;a+1,b,c+1}^j)\cdot \mathbf{F}_{6;a+1,b,c+1}^j,\\
J_8(\mbox{\mancube}_{a,b,c}^j) & = (\mathbf{F}_{5;a+1,b+1,c+1}^j \times \mathbf{F}_{4;a+1,b+1,c+1}^j)\cdot \mathbf{F}_{6;a+1,b+1,c+1}^j.
\end{align*}}

\subsection{3D discrete Cauchy-Green deformation tensor} \label{3D_Cauchy_Green}

{\footnotesize
\begin{align*}
 \mathbf{C}_2(\mbox{\mancube}_{a,b,c}^j) &=
 \begin{bmatrix}
 \langle\mathbf{F}_{2;a+1,b,c}^j,  \mathbf{F}_{2;a+1,b,c}^j\rangle & \; 
\langle \mathbf{F}_{2;a+1,b,c}^j,  \mathbf{F}_{4;a+1,b,c}^j \rangle & \; 
\langle \mathbf{F}_{2;a+1,b,c}^j , \mathbf{F}_{3;a+1,b,c}^j\rangle \\[6pt]
\langle \mathbf{F}_{4;a+1,b,c}^j,  \mathbf{F}_{2;a+1,b,c}^j \rangle  & \; 
\langle\mathbf{F}_{4;a+1,b,c}^j , \mathbf{F}_{4;a+1,b,c}^j \rangle & \; 
\langle \mathbf{F}_{4;a+1,b,c}^j,  \mathbf{F}_{3;a+1,b,c}^j\rangle \\[6pt]
\langle \mathbf{F}_{3;a+1,b,c}^j , \mathbf{F}_{2;a+1,b,c}^j\rangle  & \; 
\langle\mathbf{F}_{3;a+1,b,c}^j,  \mathbf{F}_{4;a+1,b,c}^j\rangle & \; 
\langle\mathbf{F}_{3;a+1,b,c}^j,  \mathbf{F}_{3;a+1,b,c}^j\rangle   \end{bmatrix}.\\
\mathbf{C}_3(\mbox{\mancube}_{a,b,c}^j) &=
\begin{bmatrix}
\langle\mathbf{F}_{5;a,b+1,c}^j , \mathbf{F}_{5;a,b+1,c}^j \rangle & \; 
\langle \mathbf{F}_{5;a,b+1,c}^j , \mathbf{F}_{1;a,b+1,c}^j\rangle & \; 
\langle \mathbf{F}_{5;a,b+1,c}^j , \mathbf{F}_{3;a,b+1,c}^j\rangle \\[6pt]
\langle \mathbf{F}_{1;a,b+1,c}^j , \mathbf{F}_{5;a,b+1,c}^j\rangle  & \; 
\langle\mathbf{F}_{1;a,b+1,c}^j,  \mathbf{F}_{1;a,b+1,c}^j\rangle & \; 
\langle\mathbf{F}_{1;a,b+1,c}^j,  \mathbf{F}_{3;a,b+1,c}^j\rangle \\[6pt]
\langle \mathbf{F}_{3;a,b+1,c}^j , \mathbf{F}_{5;a,b+1,c}^j\rangle  & \; 
\langle\mathbf{F}_{3;a,b+1,c}^j,  \mathbf{F}_{1;a,b+1,c}^j\rangle & \; 
\langle\mathbf{F}_{3;a,b+1,c}^j,  \mathbf{F}_{3;a,b+1,c}^j\rangle   \end{bmatrix}.\\
\mathbf{C}_4(\mbox{\mancube}_{a,b,c}^j) & =
\begin{bmatrix}
\langle\mathbf{F}_{2;a,b,c+1}^j , \mathbf{F}_{2;a,b,c+1}^j \rangle & \; 
\langle \mathbf{F}_{2;a,b,c+1}^j , \mathbf{F}_{1;a,b,c+1}^j\rangle & \; 
\langle \mathbf{F}_{2;a,b,c+1}^j , \mathbf{F}_{6;a,b,c+1}^j\rangle \\[6pt]
\langle \mathbf{F}_{1;a,b,c+1}^j , \mathbf{F}_{2;a,b,c+1}^j\rangle  & \; 
\langle\mathbf{F}_{1;a,b,c+1}^j,  \mathbf{F}_{1;a,b,c+1}^j\rangle & \; 
\langle\mathbf{F}_{1;a,b,c+1}^j,  \mathbf{F}_{6;a,b,c+1}^j\rangle \\[6pt]
\langle \mathbf{F}_{6;a,b,c+1}^j , \mathbf{F}_{2;a,b,c+1}^j\rangle  & \; 
\langle\mathbf{F}_{6;a,b,c+1}^j,  \mathbf{F}_{1;a,b,c+1}^j\rangle & \; 
\langle\mathbf{F}_{6;a,b,c+1}^j,  \mathbf{F}_{6;a,b,c+1}^j\rangle   \end{bmatrix}.\\
\mathbf{C}_5(\mbox{\mancube}_{a,b,c}^j) &=
\begin{bmatrix}
\langle\mathbf{F}_{4;a+1,b+1,c}^j , \mathbf{F}_{4;a+1,b+1,c}^j \rangle & \; \langle \mathbf{F}_{4;a+1,b+1,c}^j , \mathbf{F}_{5;a+1,b+1,c}^j\rangle & \; \langle \mathbf{F}_{4;a+1,b+1,c}^j , \mathbf{F}_{3;a+1,b+1,c}^j\rangle \\[6pt]
\langle \mathbf{F}_{5;a+1,b+1,c}^j , \mathbf{F}_{4;a+1,b+1,c}^j\rangle  & \; \langle\mathbf{F}_{5;a+1,b+1,c}^j,  \mathbf{F}_{5;a+1,b+1,c}^j\rangle & \; \langle\mathbf{F}_{5;a+1,b+1,c}^j,  \mathbf{F}_{3;a+1,b+1,c}^j\rangle \\[6pt]
\langle \mathbf{F}_{3;a+1,b+1,c}^j , \mathbf{F}_{4;a+1,b+1,c}^j\rangle  & \; \langle\mathbf{F}_{3;a+1,b+1,c}^j,  \mathbf{F}_{5;a+1,b+1,c}^j\rangle & \; \langle\mathbf{F}_{3;a+1,b+1,c}^j,  \mathbf{F}_{3;a+1,b+1,c}^j\rangle   \end{bmatrix}.\\
\mathbf{C}_6(\mbox{\mancube}_{a,b,c}^j) &=
\begin{bmatrix}
\langle\mathbf{F}_{1;a,b+1,c+1}^j , \mathbf{F}_{1;a,b+1,c+1}^j \rangle & \; \langle \mathbf{F}_{1;a,b+1,c+1}^j , \mathbf{F}_{5;a,b+1,c+1}^j\rangle & \; \langle \mathbf{F}_{1;a,b+1,c+1}^j , \mathbf{F}_{6;a,b+1,c+1}^j\rangle \\[6pt]
\langle \mathbf{F}_{5;a,b+1,c+1}^j , \mathbf{F}_{1;a,b+1,c+1}^j\rangle  & \; \langle\mathbf{F}_{5;a,b+1,c+1}^j,  \mathbf{F}_{5;a,b+1,c+1}^j\rangle & \; \langle\mathbf{F}_{5;a,b+1,c+1}^j,  \mathbf{F}_{6;a,b+1,c+1}^j\rangle \\[6pt]
\langle \mathbf{F}_{6;a,b+1,c+1}^j , \mathbf{F}_{1;a,b+1,c+1}^j\rangle  & \; \langle\mathbf{F}_{6;a,b+1,c+1}^j,  \mathbf{F}_{5;a,b+1,c+1}^j\rangle & \; \langle\mathbf{F}_{6;a,b+1,c+1}^j,  \mathbf{F}_{6;a,b+1,c+1}^j\rangle   \end{bmatrix}.\\
\mathbf{C}_7(\mbox{\mancube}_{a,b,c}^j) &=
\begin{bmatrix}
\langle\mathbf{F}_{4;a+1,b,c+1}^j , \mathbf{F}_{4;a+1,b,c+1}^j \rangle & \; \langle \mathbf{F}_{4;a+1,b,c+1}^j , \mathbf{F}_{2;a+1,b,c+1}^j\rangle & \; \langle \mathbf{F}_{4;a+1,b,c+1}^j , \mathbf{F}_{6;a+1,b,c+1}^j\rangle \\[6pt]
\langle \mathbf{F}_{2;a+1,b,c+1}^j , \mathbf{F}_{4;a+1,b,c+1}^j\rangle  & \; \langle\mathbf{F}_{2;a+1,b,c+1}^j,  \mathbf{F}_{2;a+1,b,c+1}^j\rangle & \; \langle\mathbf{F}_{2;a+1,b,c+1}^j,  \mathbf{F}_{6;a+1,b,c+1}^j\rangle \\[6pt]
\langle \mathbf{F}_{6;a+1,b,c+1}^j , \mathbf{F}_{4;a+1,b,c+1}^j\rangle  & \; \langle\mathbf{F}_{6;a+1,b,c+1}^j,  \mathbf{F}_{2;a+1,b,c+1}^j\rangle & \; \langle\mathbf{F}_{6;a+1,b,c+1}^j,  \mathbf{F}_{6;a+1,b,c+1}^j\rangle   \end{bmatrix}.\\
\mathbf{C}_8(\mbox{\mancube}_{a,b,c}^j) &=
\begin{bmatrix}
\langle\mathbf{F}_{5;a+1,b+1,c+1}^j , \mathbf{F}_{5;a+1,b+1,c+1}^j \rangle &  
\langle \mathbf{F}_{5;a+1,b+1,c+1}^j , \mathbf{F}_{4;a+1,b+1,c+1}^j\rangle &  
\langle \mathbf{F}_{5;a+1,b+1,c+1}^j , \mathbf{F}_{6;a+1,b+1,c+1}^j\rangle \\[6pt]
\langle \mathbf{F}_{4;a+1,b+1,c+1}^j , \mathbf{F}_{5;a+1,b+1,c+1}^j\rangle  & 
\langle\mathbf{F}_{4;a+1,b+1,c+1}^j,  \mathbf{F}_{4;a+1,b+1,c+1}^j\rangle &  
\langle\mathbf{F}_{4;a+1,b+1,c+1}^j,  \mathbf{F}_{6;a+1,b+1,c+1}^j\rangle \\[6pt]
\langle \mathbf{F}_{6;a+1,b+1,c+1}^j , \mathbf{F}_{5;a+1,b+1,c+1}^j\rangle  &  
\langle\mathbf{F}_{6;a+1,b+1,c+1}^j,  \mathbf{F}_{4;a+1,b+1,c+1}^j\rangle &  
\langle\mathbf{F}_{6;a+1,b+1,c+1}^j,  \mathbf{F}_{6;a+1,b+1,c+1}^j\rangle   \end{bmatrix}.
\end{align*}}

\subsection{3D derivatives of impenetrability constraints \eqref{3D_impenetrability_constraints}}\label{3D_deriv_imp_const}

\begin{equation*}{\small
\begin{aligned}
D_1\Psi_{\rm im_1} & = -(\varphi_{a-1,b,c}^j -\varphi_{a,b,c}^j) \times (\varphi_{a,b-1,c}^j -\varphi_{a,b,c}^j) - (\varphi_{a,b-1,c}^j -\varphi_{a,b,c}^j) \times (\varphi_{d,e,f}^j - \varphi_{a,b,c}^j )
\\
& \quad - (\varphi_{d,e,f}^j - \varphi_{a,b,c}^j) \times (\varphi_{a-1,b,c}^j -\varphi_{a,b,c}^j),
\\
D_2\Psi_{\rm im_1} & = (\varphi_{a,b-1,c}^j -\varphi_{a,b,c}^j) \times (\varphi_{d,e,f}^j - \varphi_{a,b,c}^j ),
\\
D_3\Psi_{\rm im_1} & = (\varphi_{d,e,f}^j - \varphi_{a,b,c}^j ) \times (\varphi_{a-1,b,c}^j -\varphi_{a,b,c}^j),
\\
D_4 \Psi_{\rm im_1}& =  (\varphi_{a-1,b,c}^j -\varphi_{a,b,c}^j) \times (\varphi_{a,b-1,c}^j -\varphi_{a,b,c}^j),
\\
D_1\Psi_{\rm im_2} & =  -  (\varphi_{a,b-1,c}^j -\varphi_{a,b,c}^j) \times (\varphi_{a+1,b,c}^j -\varphi_{a,b,c}^j) - (\varphi_{d+1,e,f}^j - \varphi_{a,b,c}^j ) \times (\varphi_{a,b-1,c}^j -\varphi_{a,b,c}^j)
\\
& \quad - (\varphi_{a+1,b,c}^j -\varphi_{a,b,c}^j) \times  (\varphi_{d+1,e,f}^j - \varphi_{a,b,c}^j ),
\\
D_2\Psi_{\rm im_2} & =  (\varphi_{d+1,e,f}^j - \varphi_{a,b,c}^j) \times (\varphi_{a,b-1,c}^j -\varphi_{a,b,c}^j),
\\
D_3\Psi_{\rm im_2} & = (\varphi_{a+1,b,c}^j -\varphi_{a,b,c}^j) \times  (\varphi_{d+1,e,f}^j - \varphi_{a,b,c}^j),
\\
D_4\Psi_{\rm im_2} & =  (\varphi_{a,b-1,c}^j -\varphi_{a,b,c}^j) \times (\varphi_{a+1,b,c}^j -\varphi_{a,b,c}^j), 
\\
D_1\Psi_{\rm im_3} & = - (\varphi_{a,b+1,c}^j -\varphi_{a,b,c}^j) \times (\varphi_{a-1,b,c}^j -\varphi_{a,b,c}^j) - (\varphi_{d,e+1,f}^j - \varphi_{a,b,c}^j ) \times (\varphi_{a,b+1,c}^j -\varphi_{a,b,c}^j),
\\ & \quad - (\varphi_{a-1,b,c}^j -\varphi_{a,b,c}^j) \times (\varphi_{d,e+1,f}^j - \varphi_{a,b,c}^j ),
\\
D_2\Psi_{\rm im_3} & = (\varphi_{d,e+1,f}^j - \varphi_{a,b,c}^j ) \times (\varphi_{a,b+1,c}^j -\varphi_{a,b,c}^j),
\\
D_3\Psi_{\rm im_3} & =(\varphi_{a-1,b,c}^j -\varphi_{a,b,c}^j) \times (\varphi_{d,e+1,f}^j - \varphi_{a,b,c}^j ),
\\
D_4\Psi_{\rm im_3} & =  (\varphi_{a,b+1,c}^j -\varphi_{a,b,c}^j) \times (\varphi_{a-1,b,c}^j -\varphi_{a,b,c}^j),
\\
D_1\Psi_{\rm im_4} & = - (\varphi_{a+1,b,c}^j -\varphi_{a,b,c}^j) \times (\varphi_{a,b+1,c}^j -\varphi_{a,b,c}^j) -  (\varphi_{d+1,e+1,f}^j - \varphi_{a,b,c}^j ) \times (\varphi_{a+1,b,c}^j -\varphi_{a,b,c}^j)
\\
& \quad - (\varphi_{a,b+1,c}^j -\varphi_{a,b,c}^j) \times (\varphi_{d+1,e+1,f}^j - \varphi_{a,b,c}^j),
\\
D_2\Psi_{\rm im_4} & = (\varphi_{a,b+1,c}^j -\varphi_{a,b,c}^j) \times (\varphi_{d+1,e+1,f}^j - \varphi_{a,b,c}^j),
\\
D_3\Psi_{\rm im_4} & = (\varphi_{d+1,e+1,f}^j - \varphi_{a,b,c}^j) \times (\varphi_{a+1,b,c}^j -\varphi_{a,b,c}^j),
\\
D_4\Psi_{\rm im_4} & =  (\varphi_{a+1,b,c}^j -\varphi_{a,b,c}^j) \times (\varphi_{a,b+1,c}^j -\varphi_{a,b,c}^j).
\end{aligned}}
\end{equation*}

{\footnotesize

\bibliographystyle{new}

}

\end{document}